\documentclass[mathematics,article,accept,moreauthors,pdftex]{Definitions/mdpi} 

\usepackage{adjustbox}

\firstpage{1} 
\pubvolume{1}
\issuenum{1}
\articlenumber{0}
\pubyear{2021}
\copyrightyear{2021}
\externaleditor{Academic Editor: Begoña Cano and Mechthild Thalhammer} 

\datereceived{4 June 2021} 
\dateaccepted{25 July 2021} 
\datepublished{} 
\hreflink{https://doi.org/} 

\usepackage{comment}
\usepackage{multirow}
\usepackage{float}
\usepackage{amsmath}
\usepackage{listings}
\usepackage{amsfonts}
\usepackage{pgf,tikz,pgfplots}
\pgfplotsset{compat=1.15}
\usepackage{mathrsfs}
\usepackage{comment}
\usetikzlibrary{arrows}
\usepackage[normalem]{ulem}

\definecolor{darkgreen}{RGB}{0, 150, 0}

\newcommand{\step}{\\[2mm]}

\usepackage[labelformat=simple]{subfig}

\newtheorem{alg}{Algorithm}

\Title{ Collocation Methods for High-Order Well-Balanced Methods for Systems of Balance Laws}

\TitleCitation{Collocation Methods for High-Order Well-Balanced Methods for Systems of Balance Laws}

\Author{\hl{Irene G}ómez-Bueno $^{1,}$*\orcidA{}, Manuel Jesús Castr\hl{o D}íaz $^{1}$, \hl{Carlos P}arés $^{1}$ and Giovanni Russo $^{2}$}

\AuthorNames{Irene Gómez-Bueno, Manuel Jesús Castro Díaz, Carlos Parés, and Giovanni Russo}

\AuthorCitation{Gómez-Bueno, I.; Castro Díaz, M.J.; Parés, C.; Russo, G.}

\address{%
$^{1}$ \quad \hl{University} of Málaga; mjcastro@uma.es (M.J.C.D.); pares@uma.es (C.P.)\\
$^{2}$ \quad \hl{University} of Catania; russo@dmi.unict.it}

\corres{Correspondence: igomezbueno@uma.es}

\abstract{In some previous works, two of the authors introduced a technique to design high-order numerical methods for one-dimensional balance laws that preserve all their stationary solutions. The basis of these methods is a well-balanced reconstruction operator. Moreover, they introduced a procedure to modify any standard reconstruction operator, like MUSCL, ENO, CWENO, etc., in order to be well-balanced. This strategy involves a non-linear problem at every cell at every time step that consists in finding the stationary solution whose average is the given cell value. In a recent \highlighting{paper} \cite{GomezCastroPares2020}
, a fully well-balanced method is presented where the non-linear problems to be solved in the reconstruction procedure are interpreted as control problems. The goal of this paper is 
 to introduce a new technique to solve these local non-linear problems based on the application of the collocation RK methods. Special care is put to analyze the effects of computing the averages and the source terms using quadrature formulas. A general technique which allows us to deal with resonant problems is also introduced. To check the efficiency of the methods and their well-balance property, they have been applied to a number of tests, ranging from easy academic systems of balance laws consisting of Burgers equation with some non-linear source terms to the shallow water equations---without and with Manning friction---or Euler equations of gas dynamics with gravity effects.
 }


\begin{document}

\vspace{6pt} 




\section{Introduction}
Let us consider a PDE system of the form:
\begin{equation} \label{sle}
U_t(x,t)+ f(U(x,t))_x= S(U(x,t)) H_x(x), \quad x\in\mathbb{R}, \, t>0,
\end{equation}
where $U(x,t)$ takes values on an open convex set $\Omega \subset \mathbb{R}^N$, $f:\Omega \longrightarrow \mathbb{R}^N$ is the flux function, $S:\Omega \longrightarrow \mathbb{R}^N$, and~$H$ is a continuous known function from $\mathbb{R}$ to $\mathbb{R}$ (possibly the identity function $H(x)=x$). It is supposed that system (\ref{sle}) is strictly hyperbolic, that is,
$D_f(U)= \displaystyle \frac{\partial f}{\partial U} (U)$ has $N$ real distinct eigenvalues $r_1 (U), \cdots, r_N (U) $ and associated eigenvectors $v_1, \cdots, v_N$. 

\textls[-18]{Systems of the form \eqref{sle} have non-trivial stationary solutions that satisfy the ODE~system:}
\begin{equation}\label{sblst}
f(U)_x = S(U)H_x.
\end{equation} 
\hl{A} numerical method is said to be well-balanced if it has discrete stationary solutions that approximate all the stationary solutions of the system or, at~least, a~relevant family of them. The~design of numerical methods that have the previous property is of great importance in problems where small perturbations of a stationary solution are going to be considered, like in many geophysical problems. Methods that are well-balanced have been proposed by many authors: for example, References
 \cite{Audusse,Bouchut,bicaper,EFN-03,EFN-CRAS,Chandrashekar15,Chandrashekar17, Desveaux16,desveauxRipa2016,Gosse00,Gosse01,Gosse02,Kappeli14,LR96,LR97,LeVeque98, Lukacova,Noelle06,Noelle07,Pelanti,PerSi1,PerSi2,Russo,Tang04,touma2016,Xing06, Klingenberg19,Klingenberg20,Kurganov2016,Kurganov2018,Chertok2015,handbook,
 franck2016finite,grosheintz2020well,kappeli2016wellb} and their~references. 
 
 The aim of this work is to introduce a general methodology to build well-balanced high-order finite-volume numerical methods for \eqref{sle} of the form:
\begin{equation} \label{met_num}
\frac{d \widetilde U_i}{dt}= -\frac{1}{\Delta x} \left( F_{i+\frac{1}{2}}(t) - F_{i-\frac{1}{2}}(t) \right) + \frac{1}{\Delta x} S_i,
\end{equation}
where: 
\begin{itemize}
\item  $I_i=\left[x_{i-\frac{1}{2}},x_{i+\frac{1}{2}}\right]$ are the computational cells, whose length $\Delta x$ is supposed to be constant for simplicity;
\item $\widetilde U_i(t)$ is the approximation of the average of the exact solution at the $i$th cell at time $t$, that is, 
\begin{equation*}
\widetilde U_i(t) \cong \frac{1}{\Delta x} \int_{x_{i-\frac{1}{2}}}^{x_{i+\frac{1}{2}}}  U(x,t) \, dx;
\end{equation*}
\item $ F_{i+\frac{1}{2}} = \mathbb{F}(U_{i+\frac{1}{2}}^{t, -}, U_{i+\frac{1}{2}}^{t, +})$, where $\mathbb{F}$ is a consistent numerical flux, and $U_{i+\frac{1}{2}}^{t, \pm}$ are the reconstructed states at the intercells, i.e.,
\begin{equation*}
U_{i+\frac{1}{2}}^{t, -}=P_i^t(x_{i+\frac{1}{2}}), \quad U_{i+\frac{1}{2}}^{t, +}=P_{i+1}^t(x_{i+\frac{1}{2}}),
\end{equation*}
where $P_i^t (x)$ is a reconstruction of order $p$ of the solution at the given $i$th cell, computed from the sequence $\{\widetilde U_i(t) \}$:
\begin{equation*}
P_i^t(x)= P_i(x; \{\widetilde U_j(t)\}_{j \in \mathcal{S}_i}),
\end{equation*}
where $\mathcal{S}_i$ is the set of cell indices belonging to the stencil associated to $P_i^t(x)$;
\item finally,
\begin{equation} \label{eq:S_i}
S_i \approx \int_{x_{i-\frac{1}{2}}}^{x_{i+\frac{1}{2}}} S(P_i^t (x)) H_x(x) \, dx.
\end{equation}
\end{itemize}
\hl{Given} a function $U$, the~following notation will be used to represent its cell averages and the approximations to its cell averages or its point values at the intercells:
$$
\bar U_i = \frac{1}{\Delta x} \int_{x_{i-1/2}}^{x_{i+1/2}} U(x) \,dx, \quad \widetilde U_i \approx \bar U_i, \quad U^{i+1/2} \approx U(x_{i+1/2}), \quad \forall i.
$$

Following Reference~\cite{sinum2008,CastroPares2019}, the~well-balanced property of the methods can be transferred to the reconstruction operator: 
\begin{Definition}\label{def:wbro}
A reconstruction operator $P_i(x)$ is said to be well-balanced for a function $U^*$ if
\begin{equation}
P_i(x)=U^*(x), \quad  \forall x \in [x_{i-\frac{1}{2}}, x_{i+\frac{1}{2}}], \, \forall i,
\end{equation}
where $P_i$ is the approximation of $U^*$ computed from the vector of cell-averages $\{ \bar U^*_i \}$ of the function~$U^*$.
\end{Definition}
In Reference~\cite{CastroPares2019}, it has been proved that the numerical method \eqref{met_num} is \textit{\hl{exactly well-balanced}} according to Definition \ref{def:ewb} if the reconstruction operator is well-balanced for every stationary solution $U^*$ and
\begin{equation} 
S_i  = \int_{x_{i-\frac{1}{2}}}^{x_{i+\frac{1}{2}}} S(P_i^t (x)) H_x(x) \, dx.
\end{equation}

\begin{Definition}\label{def:ewb}
The numerical method \eqref{met_num} is said to be exactly well-balanced if the sequence of cell-averages $\{ \bar U^*_i \}$ (or the sequence of their approximations $\{\widetilde U^*_i \}$ if a quadrature formula is used to compute them) of any stationary solution $U^*$ of \eqref{sle} is an equilibrium of the ODE system (\ref{met_num}).
\end{Definition}

Furthermore, in Reference~\cite{sinum2008,CastroPares2019}, the authors propose a method to build a well-balanced reconstruction operator from any standard one that will be recalled in the next section. As~we shall see, the~main difficulty in the application of this technique comes from the fact that, at~every cell and at every time step, one has to solve a non-linear problem of the form:

Find $U$ such that
\begin{equation}\label{ODEave}
\left\{ 
\begin{array}{l}
\displaystyle f(U)_x = S(U)H_x, \\[0.2cm]
\smallskip
\displaystyle \frac{1}{\Delta x}\int_{x_{i-1/2}}^{x_{i+1/2}} U (x) \,dx = \widetilde U_i,
\end{array}
\right.
\end{equation}
where $\widetilde U_i$ is an approximation of the average at the $i$-th cell of the sought solution of \eqref{sle}. Once this problem with given average has been solved at the cell, the~solution of \eqref{ODEave}, that will be denoted by $U_i^*$, has to be extended to the whole stencil by solving two Cauchy problems of the form
\begin{equation}\label{cauchy_stencil}
 \left\{ 
 \begin{array}{l}
 \displaystyle f(U)_x = S(U)H_x, \\
 U(x_{i + 1/2}) = U^{i + 1/2}.
 \end{array}
 \right.
 \end{equation}
\hl{More} precisely, \eqref{cauchy_stencil} with initial condition $U(x_{i+1/2}) = U^*(x_{i+1/2})$ (resp. with final condition $U(x_{i-1/2}) = U^*(x_{i-1/2})$ ) has to be solved forward in space (resp. backward in space) in the cells at the right (resp. at the left) of the $i$th cell in its~stencil.

The previous procedure has been applied to balance laws whose stationary solutions are known in implicit or explicit form, so that problems \eqref{ODEave} and \eqref{cauchy_stencil} can be easily solved: see, for example, Reference~\cite{lopez2013, handbook}, and~the references therein. However, when solving these problems by analytic procedures is not possible, they have to be solved numerically. Please note that the two main difficulties come from the condition on the average in \eqref{ODEave} and from the fact that the ODE is not in normal form: special care has to be taken in \textit{\hl{resonant problems}} i.e.,~problems in which the Jacobian matrix of the flux,
$D_f(U)$ becomes singular. When the Jacobian is regular, the~ODE system \eqref{sblst} can be written in normal form
\begin{equation}\label{sblstnormal}
U_x = D_f(U)^{-1} S(U) H_x,
\end{equation}
and the solution of the problems is easier. At~a sonic point, i.e., a point in which one of the eigenvalues vanishes, the~ODE system may have no solution or may have more than~one. 

In recent papers~\cite{GomezCastroPares2020, gomez2021well}, a strategy based on the interpretation of \eqref{ODEave} as a control problem has been proposed to deal with the first difficulty: the value of $U^*_i$ at the leftmost intercell of the stencil is the control that has to be selected so that the condition on the average is satisfied. Newton's method is used to find the control, where the gradient is computed by solving the adjoint problem. A~standard ODE solver, namely RK4, is used to solve the state and the adjoint equations, and~to extend the solution to the whole~stencil. 

 The first goal of this paper is to present a general framework to design and analyze well-balanced high-order numerical methods in which problems \eqref{ODEave} and \eqref{cauchy_stencil} are numerically solved: the definition of the well-balanced property and a general result that allows one to check it will be stated. The~second goal is then to introduce a methodology based on the resolution of \eqref{ODEave} and \eqref{cauchy_stencil} using RK collocation methods that leads to numerical methods that are well-balanced according to the given definition. As~we shall see, the~condition on the average in \eqref{ODEave} will naturally lead to non-linear problems that will be solved by an iterative procedure. Finally, the~third goal is to propose a strategy to deal with resonant problems that is general, albeit problem-dependent: this strategy will be illustrated in the particular case of the shallow water~system.

 The organization of the article is as follows: Section~\ref{sec2} is devoted to state the general framework. Definitions of well-balanced methods and reconstruction operators adapted to the numerical solution of the local problems will be introduced and a general result will be stated showing that well-balanced reconstruction operators lead to well-balanced methods. In~Section~\ref{Section_collocation}, we introduce a strategy based on the use of RK collocation methods to solve problems \eqref{ODEave} and \eqref{cauchy_stencil}, and we show that it leads to numerical schemes that are well-balanced. Section~\ref{Section_collocation} ends with the introduction of the technique that allows us to deal with resonant~problems.

In Section~\ref{numerical_experiments}, a~number of numerical tests are presented to analyze the performance of the methods and to check their accuracy and well-balancedness. Some numerical tests for scalar and systems of balance laws are considered: Burgers' equation with a non-linear source term, the~shallow water system---both with and without Manning friction---or the Euler equations of gas dynamics including gravity effects are considered. In~addition, we also present some numerical results to check that such numerical methods introduced in this work behave correctly in the presence of critical states. Moreover, we show that the numerical methods are able to preserve subsonic and supersonic moving stationary solutions for the compressible Euler equations with gravitational force. Finally, some conclusions are drawn in Section~\ref{sec5}, and further developments are also~discussed.
 
 \section{Well-Balanced~Methods}
 \label{sec2}

The following strategy to obtain a reconstruction operator 
$$
P_i(x) = P_i(x;\{\bar U_j\}_{j \in \mathcal{S}_i})
$$
that is well-balanced according to Definition \ref{def:wbro}, from~a standard one
$$
Q_i(x) = Q_i(x;\{\bar U_j\}_{j \in \mathcal{S}_i}),
$$
was introduced in Reference~\cite{sinum2008}:

\begin{alg}  \label{alg:ewbrec} \hl{Given a family} of cell values $\{\bar U_i\}$, at~every cell $I_i$:
\begin{enumerate}
\item Find, if~possible, a~stationary solution $U_i^*(x)$ defined in the stencil of cell $I_i$, 
$\displaystyle \cup_{j\in\mathcal{S}_{i}}I_j$, such~that:
\begin{equation} \label{step1}
\frac{1}{\Delta x} \int_{x_{i-\frac{1}{2}}}^{x_{i+\frac{1}{2}}} U_i^* (x) \, dx = \bar U_i.
\end{equation}
 Otherwise, take $U_i^* \equiv 0$.
 
\item Compute the \textit{fluctuations} $\{V_j\}_{j \in S_i}$ given by
\begin{equation*}
V_j= \bar U_j - \frac{1}{\Delta x} \int_{x_{j-\frac{1}{2}}}^{x_{j+\frac{1}{2}}} U_i^* (x) \, dx ,\quad j \in \mathcal{S}_i, 
\end{equation*}
and compute the reconstruction operator:
\begin{equation*}
Q_i(x)=Q_i(x;\{V_j\}_{j \in \mathcal{S}_i}).
\end{equation*}
\item Finally, define
\begin{equation} \label{step3}
P_i(x)=U_i^*(x)+Q_i(x).
\end{equation}

\end{enumerate}
\end{alg}

According to Reference~\cite{sinum2008}, the~reconstruction operator $P_i$ is well-balanced if $Q_i$ is exact for the null function, conservative if $Q_i$ is conservative \begin{equation*}
\frac{1}{\Delta x} \int_{x_{i-\frac{1}{2}}}^{x_{i+\frac{1}{2}}} P_i (x) \, dx = \bar U_i, \, \forall i,
\end{equation*} 
and $P_i$ is of the same order of accuracy $p$ of $Q_i$ provided that the stationary solutions are~smooth.

Observe that, if~it is impossible to find a stationary solution defined in the stencil that satisfies \eqref{step1}, then the standard reconstruction is used. Please note that this choice does not spoil the well-balanced character of the numerical method: in this case, the~cell values in the stencil cannot be the averages of a stationary solution (otherwise, there would be at least one solution $U^*_i$); thus, there is no local equilibrium to preserve. On~the other hand, if~there is more than one stationary solution defined on the stencil that satisfies \eqref{step1}, a~criterion is needed to select one of them. This criterion depends on the particular problem and it will be discussed~later.

Notice that the use of quadrature formulas for the computation of the integral \eqref{met_num} appearing in the right-hand side of \eqref{eq:S_i} may destroy the well-balanced property of the method. In Reference~\cite{CastroPares2019}, the authors proposed to write the source term as follows:
\begin{equation}
\begin{split}
\int_{x_{i-\frac{1}{2}}}^{x_{i+\frac{1}{2}}} S(P_i^t (x)) H_x(x) \, dx = f\left( U_i^{t,*}(x_{i+\frac{1}{2}})\right) - f\left( U_i^{t,*}(x_{i-\frac{1}{2}})\right) \\
+\int_{x_{i-\frac{1}{2}}}^{x_{i+\frac{1}{2}}} \left((S(P_i^t (x)) -S(U_i^{t,*}(x)) \right) H_x(x) \, dx,
\end{split}
\end{equation}
where $U_i^{t,*}$ is the stationary solution found in (\ref{step1}). Finally, the~integral of the source term is approximated as follows:
\begin{equation}\label{Siquad}
    S_i =  f\left( U_i^{t,*}(x_{i+\frac{1}{2}})\right) - f\left( U_i^{t,*}(x_{i-\frac{1}{2}})\right) 
+ \Delta x \sum_{m=1}^M b_m \left( S(P_i^t (x^m_i)) -S(U_i^{t,*}(x^m_i)) \right) H_x(x^m_i), 
\end{equation}
where $x^m_i$, $b_m$, $m= 1, \dots, M$ are, respectively, the nodes and the weights of the quadrature formula chosen in the cell 
$I_i$, whose order of accuracy $s$ is larger or equal than $p$. It can be checked that, if the reconstruction operator is well-balanced for $U^*$, then the numerical method \eqref{met_num}--\eqref{Siquad} is exactly well-balanced for $U^*$. 

If the quadrature formula is used, as well, to compute cell averages, i.e.,
\begin{equation}\label{appave}
\bar U^*_i \approx \widetilde U^*_i = \sum_{m=1}^M b_m U^*(x^m_i),
\end{equation}
the reconstruction algorithm has to be modified as follows (see Reference~\cite{CastroPares2019}):
\begin{alg}\label{alg:ewbnirec} \hl{Given a} family of cell values $\{\widetilde U_i\}$, at~every cell $I_i$: 
\begin{enumerate}
\item Find, if~possible, a~stationary solution $U_i^*(x)$ defined in the stencil of cell $I_i$, (\/$\cup_{j\in\mathcal{S}_{i}}I_j$) such that
\begin{equation} \label{step1qf}
\sum_{m=1}^M b_m U_i^* (x^m_i)  = \widetilde U_i.
\end{equation}
Otherwise, take $U_i^* \equiv 0$.
\item Compute the \textit{fluctuations} $\{V_j\}_{j \in S_i}$ given by
\begin{equation*}
V_j= \widetilde U_j -  \sum_{m=1}^M b_m U_i^* (x^m_j),\quad j \in \mathcal{S}_i, 
\end{equation*}
and compute the reconstruction operator:
\begin{equation*}
Q_i(x)=Q_i(x;\{V_j\}_{j \in \mathcal{S}_i}).
\end{equation*}
\item Define
\begin{equation*} 
P_i(x)=U_i^*(x)+Q_i(x).
\end{equation*}
\end{enumerate}
\end{alg}
Again, it can be proved that this reconstruction operator is well-balanced and leads to numerical methods that are well-balanced according to Definition \ref{def:ewb}.

The main difficulty to implement the well-balanced reconstruction procedure given by Algorithm \ref{alg:ewbnirec} comes from the first step in which the values at the quadrature points and the intercells of the stationary solution that satisfies \eqref{step1qf} have to be computed in the stencil. When it is not possible to solve these problems by analytical procedure, a~numerical one is required. In~this case, the~generic local discrete problem to be solved is the following:
\begin{Problem}[Local problem (LP)]
Given an index $i$ and a state $\widetilde W \in \Omega$, find approximations 
$$
U^{*,m}_{i, j}, \  m=1, \dots, M, \ j \in \mathcal{S}_i; \quad U_{i}^{*, i\pm 1/2},
$$
of the values
$$
U^*_i(x^m_j), \  m=1, \dots, M, \ j \in \mathcal{S}_i; \quad U^*_i(x_{i \pm 1/2}),
$$
where $U^*_i$ is the stationary solution that satisfies
\begin{equation} \label{localproblem}
\sum_{m=1}^M b_m U_i^* (x^m_i)  = \widetilde W.
\end{equation}
\end{Problem}

Although the approximations to be found depend in general on $\Delta x$ and $\widetilde W$, this dependency will not be explicitly written to avoid an excess of~notation.

 Once a numerical solver has been selected for the local discrete problems, the~approximations of the stationary solution satisfying \eqref{localproblem} at the quadrature points $U^{*,m}_{i,j}$ and at the intercells $U_{i}^{*, i\pm 1/2}$ are used to define the reconstruction operator as follows:
\begin{alg} \label{alg:wbrec} \hl{Given a} sequence of cell values $\{\widetilde U_i\}$, at~every cell $I_i$: 
\begin{enumerate}
\item Apply, if~possible, the~local solver at the $i$-th cell with $\widetilde W = \widetilde U_i$ to obtain 
$$
U^{*,m}_{i, j}, \quad m=1, \dots, M, \ j\in \mathcal{S}_i; \quad U_{i}^{*, i\pm 1/2}.
$$
 Otherwise, take 
$$ 
U^{*,m}_{i, j} = 0, \quad m=1, \dots, M, \ j\in \mathcal{S}_i; \quad U_{i}^{*, i\pm 1/2}=0.
$$
\item Compute the fluctuations $\{V_j\}_{j \in S_i}$ given by
\begin{equation*}
V_j=\widetilde U_j - \sum_{m=1}^M b_m  U^{*,m}_{i, j}, \quad j \in \mathcal{S}_i, 
\end{equation*}
and compute:
\begin{equation*}
Q_i(x)=Q_i(x;\{V_j\}_{j \in \mathcal{S}_i}).
\end{equation*}
\item Define
\begin{eqnarray} \label{step3bisa}
& & P^m_i = U^{*,m}_{i,i} + Q_i(x^m_i), \quad m=1, \dots, M, \\
& & U^+_{i - 1/2}  =  U_{i}^{*, i - 1/2} + Q_i(x_{i - 1/2}), \\
& & U^-_{i + 1/2}  =  U_{i}^{*, i+ 1/2} + Q_i(x_{i + 1/2}). 
\end{eqnarray}
\end{enumerate}
\end{alg}

Then, the~numerical method is defined by \eqref{met_num} with
\begin{equation}\label{Siquadbis}
    S_i =  f\left( U_{i}^{*, i + 1/2}\right) - f\left( U_{i}^{*, i- 1/2}\right) 
+ \Delta x \sum_{m=1}^M b_m \left( S(P^m_{i}) -S(U_{i,i}^{*,m}) \right) H_x(x^m_i), 
\end{equation}
where the dependency on $t$ has not been explicitly written 
to avoid an excess of~indices. 

 The question is now what is the well-balanced property satisfied by the previous method. The~following definition is considered:
\begin{Definition}\label{defwb}
The numerical method \eqref{met_num} is said to be well-balanced with order $q \geq p$ if for every stationary solution $U^*$ of \eqref{sle}, and, for every $\Delta x$, there exists
a discrete stationary solution $\{\widetilde U^*_{\Delta x, i} \}$, i.e.,~an equilibrium of \eqref{met_num}, such that
\begin{equation}\label{orderq}
    \bar U^*_i = \widetilde U^*_{\Delta x, i}  + O(\Delta x^q),\quad \forall i.
\end{equation}
\end{Definition}

In order to illustrate the definition, let us consider some examples. We consider first the equation:
\begin{equation}\label{lintrans}
u_t + u_x = u, 
\end{equation}
whose stationary solutions are
$$
u(x) = Ce^x, \quad C \in \mathbb{R}.
$$
\hl{The} cell-averages will be computed by the mid-point value so that
$$
\frac{1}{\Delta x}\int_{x_{i-1/2}}^{x_{i+1/2}} u_0 (x) \,dx \approx u_0(x_i).
$$
\hl{We} consider the upwind numerical flux and different discretizations of the source term:

\begin{itemize}

\item Method 1:
\begin{equation}\label{wb:ex1}
\frac{d u_i}{dt} = \frac{1}{\Delta x}\left(  e^{\Delta x}u_{i-1} - u_i \right).
\end{equation}
This is the first-order well-balanced numerical method based on the well-balanced piecewise constant reconstruction operator described in Algorithm \ref{alg:ewbnirec}. The~discrete stationary solutions satisfy
$$
u_i =  e^{\Delta x}u_{i-1}, \quad i=0,1,\dots
$$
It is clear then that the set of discrete stationary solutions is given by the point values of the continuous ones:
$$
\{ C e^{x_i} \}, \quad C \in \mathbb{R},
$$
and the numerical method is the exactly well-balanced according to Definition \eqref{def:ewb}.

\item Method 2:
\begin{equation}\label{wb:ex2}
\frac{d u_i}{dt}  =  \frac{1}{\Delta x}\left(u_i - u_{i-1}  \right) + \frac{1}{2}(u_{i-1} + u_i).
\end{equation}
This is the first-order well-balanced numerical method based on the RK2 collocation method that will be introduced in Section~\ref{Section_collocation}.
The discrete stationary solutions satisfy
$$
u_i  =  u_{i-1} + \frac{\Delta x}{2}(u_{i-1} + u_i),
$$
which leads to
$$
u_i = \frac{1 + \Delta x/2}{1-\Delta x /2} u_{i-1}, \quad i = 1, 2, \dots
$$
The set of discrete stationary solutions is given by:
$$
\left\{ C\left( \frac{1 + \Delta x/2}{1-\Delta x /2} \right)^i \right \}, \quad C \in \mathbb{R}.
$$
One has for every $C$:
$$
u(x_i) = C e^{i \Delta x} =C \left( \frac{1 + \Delta x/2}{1-\Delta x /2} \right)^i + O(\Delta x^2), \quad \forall i,
$$
so that the numerical method is well-balanced with order~2.

\item Method 3:
\begin{equation}\label{wb:ex3}
\frac{d u_i}{dt}  =  \frac{1}{\Delta x}\left( u_{i-1} - u_i \right) +  u_{i-1}.
\end{equation}
This is the numerical method that corresponds to the upwind treatment of the source term (see Reference~\cite{bermudez1994upwind}). The~discrete stationary solutions satisfy
$$
u_i  =  u_{i-1} + \Delta x u_{i-1},
$$
which leads to
$$
u_i = (1 + \Delta x)u_{i-1}, \quad i = 1, 2, \dots
$$
The set of discrete stationary solutions is given by:
$$
\left\{ C(1 + \Delta x)^i \right \}, \quad C \in \mathbb{R}.
$$
One has for every $C$:
$$
u(x_i) = C e^{i \Delta x} = C(1 + \Delta x)^i + O(\Delta x), \quad \forall i,
$$
so that the numerical method is well-balanced with order~1.

\item Consider now the equation
$$
u_t + \left( \frac{u^2}{2} \right)_x =  u,
$$
whose stationary solutions are
$$
u = x + constant,
$$
and $u \equiv 0$. Let us consider the numerical method:
\begin{equation}\label{wb:ex4}
\frac{du_i}{dt} = \frac{1}{\Delta x} \left(F_{i-1/2} - F_{i+1/2} \right) + \bar u_{i-1/2},
\end{equation}
where
$$
\bar u_{i+1/2} = \frac{u_i + u_{i+1}}{2}, \quad F_{i+1/2} = \begin{cases} 
{u_i^2}/{2} & \text{if $\bar u_{i+1/2} > 0$;}\\
0 & \text{if $\bar u_{i+1/2} = 0$;}\\
{u_{i+1}^2}/{2} & \text{if $\bar u_{i+1/2} < 0$.}
\end{cases}
$$
The discrete stationary solutions satisfy
$$
F_{i-1/2} - F_{i+1/2} + \Delta x \bar u_{i-1/2} = 0.
$$

Let us see that there are no discrete stationary solutions that change their sign from negative to positive (as it happens with all the stationary solutions)

$$
u^*(x) = C + x, \quad C < 0
$$
if the interval contains point $x_0=-C$.

\begin{itemize}
    \item Let us suppose that there exists a discrete stationary solution such that $\bar u_{i-1/2} < 0$, $\bar u_{i+1/2} > 0$. Then, at~the $i$-th cell, we would have:
    $$
    0 = \frac{u_i^2}{2} - \frac{u_i^2}{2} + \Delta x \bar u_{i-1/2}  \implies \bar u_{i-1/2} = 0,
    $$
    which is a~contradiction.
 
    \item Let us suppose now that $\bar u_{i-3/2} < 0$, $\bar u_{i-1/2} = 0 $, $\bar u_{i+1/2} > 0$. Then, at~the $i$-th cell, we have
    $$
  0 =   - \frac{u_i^2}{2} + \Delta x\, \bar u_{i-1/2} = - \frac{u_i^2}{2} \implies u_i = 0.
  $$
 Since $\bar u_{i-1/2} = 0$, we have that $u_{i-1} = 0$, as well. Therefore, at~the $(i-1)$-th cell, we have
  $$
  0 = \frac{u_{i-1}^2}{2} + \Delta x\, \bar u_{i-3/2} = \Delta x\, \bar u_{i-3/2}  \implies \bar u_{i-3/2}  = 0,
$$
which is again a contradiction.
\end{itemize}

Therefore, it is not possible to find discrete stationary solutions that change from negative to positive. The~numerical method is not well-balanced in spite of the fact that it is consistent and~stable.

\end{itemize}

\begin{Remark}
Some authors impose $q > p$ in the definition of well-balanced, i.e.,~stationary solutions have to be approximated with enhanced accuracy. Please note that, although~any standard numerical method of order $p$ will provide approximations of any smooth stationary solution with order of accuracy $p$, these approximations are not necessarily discrete stationary solutions. This is the case of the fourth example above: although the numerical method will provide first order numerical approximations of the stationary solutions of the form $u^*(x) = x + C$ with $C <0$, these approximations cannot be discrete stationary solutions. 
\end{Remark}

\begin{Remark} \label{rem:wpic}
From the practical point of view, it is important to know how to compute discrete stationary solutions in order to 'prepare well' the initial conditions. Think, for instance, of a numerical experiment for Equation \eqref{lintrans} with initial condition 
$$
u_0(x) = u^*(x) + \delta(x),
$$
where $u^*(x)$ is a stationary solution, say $u^*(x) = e^x$, and~$\delta(x)$ a small perturbation with compact support. If~the initial condition is computed by 
$$
u^0_i = e^{x_i} + \delta(x_i), \quad i = 1, 2\dots,
$$
and Method 1 is used, the~stationary solution will remain unperturbed in regions where the wave generated by the initial perturbation is not arrived. Nevertheless, if~Method 2 is used, the~numerical solution will move everywhere to fit a discrete stationary solution of the problem. If~the amplitude of the initial perturbation is of order $O(\Delta x)$, the wave can be lost among the numerical errors at the beginning of the experiment. If, instead, the~initial condition is computed using the discrete stationary solution that approximates $u^*$,
$$
u^0_i =  \left( \frac{1 + \Delta x/2}{1-\Delta x /2} \right)^i + \delta(x_i), \quad i = 1, 2\dots
$$
the numerical solution will remain again unperturbed until the arrival of the~wave. 

\end{Remark}

According to Definition \ref{defwb}, in~order to be well-balanced, a numerical method has to have discrete stationary solutions. Let us discuss what are the equilibria of \eqref{met_num}--\eqref{Siquadbis} when the well-balanced reconstruction operator given by Algorithm \ref{alg:wbrec} is used. To~do this, let us first introduce the following discrete version of Definition \ref{def:wbro}:

\begin{Definition} The reconstruction operator defined by Algorithm \ref{alg:wbrec} is said to be well-balanced for a sequence of cell values $\{ \widetilde U_{i} \}$ if, for~every $i$,
\begin{eqnarray}
& & \sum_{m=1}^M b_m U^{*,m}_{i,j} = \widetilde U_{j}, \quad j \in \mathcal{S}_i, \label{wbseq1} \\ 
& &  U_{i+1/2}^-=U_{i}^{*, i + 1/2} =  U_{i+1}^{*, i +  1/2}=U_{i+1/2}^+, \label{wbseq2}
\end{eqnarray}
where $U^{*,m}_{i,j}$, $U_{i}^{*, i\pm 1/2}$ is the output given by the solver of local problems 
for $\widetilde W = \widetilde U_{i} $ at the $i$th cell.
\end{Definition}

In other words, the~reconstruction procedure defined by Algorithm \ref{alg:wbrec} is well-balanced for a sequence $\{ \widetilde U_{i} \}$ if (a) the cell values are recovered when the quadrature formula is applied to the approximations at the points given by the local solver and (b) the approximations at the intercells do not depend on the~stencil.

Like in the continuous case, there is a direct relation between discrete stationary solutions and well-balanced reconstruction operators:

\begin{Theorem}\label{th:wb1}
If the reconstruction operator described in Algorithm \ref{alg:wbrec} is well-balanced for a sequence of cell values
$\{ \widetilde U^*_{ i} \}$, then the sequence is a discrete stationary solution, i.e., an equilibrium of the ODE system given by \eqref{met_num}--\eqref{Siquadbis}.

\end{Theorem}
\begin{proof}
Let us show that the sequence 
$\{ \widetilde U^*_{i} \}$ is an equilibrium of \eqref{met_num}--\eqref{Siquadbis}. Due to the well-balancedness of the reconstruction operator, at~the $i$th cell, we have:
$$
V_j= \widetilde U^*_{j} - \sum_{j=1}^M b_m  U^{*,m}_{i,j} = 0, \quad j \in \mathcal{S}_i, 
$$
so that $Q_i(x) = 0$. Therefore,
\begin{eqnarray*} 
& & P^m_i = U^{*,m}_{i,i}, \quad m=1, \dots, M, \\
& & U^+_{i - 1/2}  = U_{i}^{*, i - 1/2}, \\
& & U^-_{i + 1/2}  = U_{i}^{*, i + 1/2}. 
\end{eqnarray*}
As a consequence, \eqref{Siquadbis} reduces to
$$
S_i =  f\left( U_{i}^{*, i + 1/2}\right) - f\left( U_{i}^{*, i - 1/2} \right).
$$
On the other hand:
\begin{eqnarray*}
F_{i+\frac{1}{2}} - F_{i+\frac{1}{2}}  & = &  \mathbb{F}(U_{i+\frac{1}{2}}^{ -}, U_{i+\frac{1}{2}}^{ +})
- \mathbb{F}(U_{i-\frac{1}{2}}^{ -}, U_{i-\frac{1}{2}}^{ +}) \\
& = &  \mathbb{F}(U_{i}^{*, i + 1/2}, U_{i+1}^{*, i + 1/2})
- \mathbb{F}(U_{i-1}^{*, i - 1/2},U_{i}^{*, i - 1/2})\\
& = &  f\left( U_{i}^{*, i + 1/2}\right) - f\left(U_{i}^{*, i - 1/2}\right).
\end{eqnarray*}
In the previous expression, the~consistency of the numerical flux and \eqref{wbseq2} have been used. Therefore, the~right-hand side of \eqref{met_num} vanishes.
\end{proof}

The previous result gives the condition under which a sequence of cell values defines a discrete stationary solution. In~order to prove that the method is well-balanced, one has to prove that, given any stationary solution $U^*$ and $\Delta x$, it is possible to construct a discrete one that approximates $U^*$ with order of accuracy greater or equal than $p$. If~we observe the three first examples above, we see that the discrete stationary solutions satisfy a consistent and well-defined ODE solver of the equation satisfied by the continuous ones
$$
u_x = u,
$$
the exact solver for Method 1, the~trapezoidal scheme for Method 2, and~forward Euler scheme for Method 3. This is not the case for the fourth example. The~following result states that, if~there exists an ODE solver of \eqref{sblst} with order of accuracy greater or equal than $p$ that is \textit{consistent} in some sense with the local problems solver, then the numerical method is well-balanced:

\begin{Theorem}\label{th:wb2}
Let us suppose that there exists an ODE solver of \eqref{sblst} with order of accuracy $q \geq p$ that, given any stationary solution $U^*$ and any value of $\Delta x$, provides approximations
$$
U^{*, m}_{\Delta x, i}, \  m=1, \dots, M; \quad U_{\Delta x}^{*, i + 1/2}, \quad \forall i,
$$
of the values of $U^*$ at the quadrature points and the intercells
$$
U^*(x^m_i), \  m=1, \dots, M; \quad U^*(x_{i +  1/2}),  \quad \forall i,
$$
in such a way that, at~every cell $I_i$, the~solution of the local discrete problem \eqref{localproblem} with $\widetilde W$ equal to
\begin{equation}
\label{discstat}
\widetilde U^*_{\Delta x, i} := \sum_{m=1}^M b_m U^{*,m}_{\Delta x, i},
\end{equation}
is given by
\begin{eqnarray}
& & U^{*,m}_{i, j} = U^{*,m}_{\Delta x, j}, \quad m=1, \dots, M, \ j\in \mathcal{S}_i;\\
& &  U_{i}^{*, i \pm 1/2} = U_{\Delta x}^{*, i \pm 1/2}.
\end{eqnarray}
Then, the~numerical method
\eqref{met_num}--\eqref{Siquadbis} is well-balanced. 
\end{Theorem}

\begin{proof}
The sequence $\{ \widetilde U^*_{\Delta x, i} \}$ given by \eqref{discstat} approaches the cell values of $U^*$ with order at least $p$ since both the orders of the ODE solver and the quadrature formula are greater or equal than $p$. 
On the other hand, it is obviously well-balanced for the reconstruction operator since
$$
 \sum_{m=1}^M b_m U^{*,m}_{i,j}  =  \sum_{m=1}^M b_m U^{*,m}_{\Delta x,j} = \widetilde U^*_{\Delta x, j}
 $$
 and, furthermore
 $$
 U_{i}^{*,i+1/2} = U_{i+1}^{*,i+1/2} = U_{\Delta x}^{*, i + 1/2}
 $$
 so that \eqref{wbseq1} and \eqref{wbseq2} are satisfied. As~a consequence of Theorem \ref{th:wb1},
 $\{ \widetilde U^*_{\Delta x, i} \}$ is a discrete stationary solution, which completes the proof.
 \end{proof}

 In a recent paper~\cite{GomezCastroPares2020}, the~ODE solver chosen for \eqref{sblst} was RK4. The~local problems were interpreted as control problems in which both the state and the adjoint local problems were solved using RK4, as well, which gave the consistency of both solvers and, thus, the well-balanced property of the~method.

\section{Collocation~Methods} \label{Section_collocation}
In this paper, we propose to build well-balanced numerical methods for \eqref{sle} on the basis of collocation RK methods. Remember that these methods admit a double interpretation: on the one hand, they are standard RK methods based on a Butcher
tableau 
$$
 \begin{array}{c|ccc}
 c_1 & a_{1,1} &\dots & a_{1,s} \\
 c_2 & a_{2,1} & \dots &  a_{2,s} \\ 
  \vdots &\vdots & \ddots & \vdots \\
 c_s & a_{s,1} & \dots & a_{s,s}\\
 \hline
  & b_1 &  \dots & b_s.
  \end{array}.
 $$

Given an index $i_0$, when these methods are applied to a Cauchy problem
\begin{equation}\label{Cauchynormal}
     \left\{
     \begin{array}{l}
     \displaystyle U_x = G(x, U), \\
     \displaystyle U(x_{i_0-1/2}) = U^{i_0-1/2},
     \end{array}
     \right.
 \end{equation}
 in a uniform 
 mesh of nodes $x_{i+1/2} = x_{i-1/2} +   \Delta  x$, $i=i_0,i_0+1,\dots$ , the~numerical solutions are updated as follows:
\begin{equation}
U^{i+1/2}  = U^{i-1/2} + \Delta x\, \Phi_{\Delta x}(U^{i-1/2}),\quad i = i_0,i_0+1,\dots
 \end{equation}
 where 
 $$
 \Phi_{\Delta x} (U^{i-1/2}) = \sum_{j = 1}^s b_j K_i^j.
 $$
 Here, $K_i^1, \dots, K_i^s$ solve the non-linear system
\begin{equation}\label{stages}
 K_i^j = G\left(x_i^j, U^{i-1/2} + \Delta x \sum_{l=1}^s a_{j,l} K_i^l \right), \quad j= 1, \dots, s,
 \end{equation}
 where
\begin{equation} \label{Gausspoints} 
 x_i^j = x_{i-1/2} + c_j\Delta x, \quad j= 1, \dots, s.
 \end{equation}
 
 On the other hand, they can be interpreted as follows: 
 $$
 U^{i+1/2} = P_i (x_{i+1/2}),
 $$
 where $P_i$ is the only polynomial of degree $s$ that satisfies:
\begin{equation}
     \left\{
 \begin{array}{l}
\displaystyle P_i(x_{i-1/2})  =  U^{i-1/2}, \\\smallskip
\displaystyle P'_i(x_i^j)  = G(x_i^j, P_i(x_i^j)), \quad j= 1, \dots, s. 
 \end{array}
 \right.
 \end{equation}
 
 We will consider here Gauss-Legendre methods, in~which $x_i^1,\dots, x_i^s$ and $b_1, \dots, b_s$ are, respectively, the quadrature points and the weights of the Gauss quadrature formula in the interval $[x_{i-1/2}, x_{i+1/2}]$. This quadrature formula will be used to compute the averages at the cells.
 The order of accuracy of these methods is $2s$. 
 
 Because of this double interpretation, it can be easily shown that Gauss methods are symmetric or reversible in the following sense (see Reference~\cite{haier2006geometric}):
\begin{equation}
\Phi_{\Delta x} \circ \Phi_{-\Delta x }=Id, \quad \text{or equivalently} \quad \Phi_{\Delta x}=\Phi_{-\Delta x}^{-1}.
\end{equation}

\subsection{Cauchy Problem~Solver}\label{ss:cauchy}

 Let us describe how the collocation RK methods will be adapted to solve problems of the form \eqref{cauchy_stencil}. The~following algorithm will be used to approximate the solution of a Cauchy problem with initial condition
$$U^*(x_{i_0-1/2}) = U^{i_0-1/2}$$
forward and backward in space:

\begin{alg} \label{alg:gp} \hl{Numerical solver} for the Cauchy problem \eqref{cauchy_stencil} using collocation RK~methods. 
\begin{itemize}

\item For $i = i_0, i_0 + 1 ,\dots$
\begin{itemize}
    \item Compute $K_i^m$, $m=1, \dots, s$ by solving the non-linear system:
    \begingroup\makeatletter\def\f@size{8}\check@mathfonts
\def\maketag@@@#1{\hbox{\m@th\normalsize\normalfont#1}}%
\begin{equation}\label{globsys1}
D_f\left( U^{i-1/2} + \Delta x \sum_{k=1}^s a_{m,k} K_i^k \right) K_i^m = 
S\left(U^{i-1/2} + \Delta x \sum_{k=1}^s a_{m,k} K_i^k \right) H_x( x_i^m), \quad m= 1, \dots, s.
\end{equation}
\endgroup

    \item Compute: 
    $$
U^{i+1/2}  =  U^{i-1/2} + \Delta x 
\sum_{m = 1}^s b_m K_i^m.
$$

\item Compute:
$$
    U_i^m  =  U^{i-1/2} + \Delta x \sum_{k=1}^s a_{m,k} K_{i}^k, \quad 
    m = 1, \dots, s.
$$

\end{itemize}

\item For $i = i_0, i_0 - 1, i_0 - 2, \dots$

\begin{itemize}
\item Compute $K_{i-1}^m$ by solving the non-linear system:
\begin{adjustwidth}{-4.6cm}{0cm}   
\begin{equation}\label{globsys2}
D_f\left( U^{i-1/2} - \Delta x \sum_{k=1}^s a_{m,k} K_{i-1}^k \right) K_{i-1}^m = 
S\left(U^{i-1/2} - \Delta x \sum_{k=1}^s a_{m,k} K_{i-1}^k \right) H_x( x_{i-1}^m), \quad m= 1, \dots, s.
\end{equation}\end{adjustwidth} 

    \item Compute: 
    $$
U^{i-3/2}  =  U^{i-1/2} - \Delta x 
\sum_{m = 1}^s b_m K_{i-1}^m.
$$

\item Compute:
$$
    U^{m}_{i-1}  =  U^{i-1/2} - \Delta x \sum_{k=1}^s a_{m,k} K_{i-1}^k, \quad m = 1, \dots, s.
$$

\end{itemize}

\end{itemize}

\end{alg}

With the notation of Theorem \ref{th:wb2}, this ODE solver provides the following approximations of a stationary solution $U^*$ at the quadrature points and intercells:
$$
 U_{\Delta x}^{*, i+1/2}  =  U^{i+1/2}, \quad  
 U^{*,m}_{\Delta x, i}    =  U^m_i, \quad m=1, \dots, s, \quad \forall i,
 $$
and the discrete stationary solutions will be then given by
 $$
 \widetilde{U}^*_{\Delta x, i} = \sum_{m = 1}^s b_m U^{*, m}_{\Delta x, i}.
 $$

The non-linear systems \eqref{globsys1} and \eqref{globsys2} are solved by using a fixed-point algorithm. If~the Jacobian of the flux function is never singular, this algorithm is equivalent to apply the collocation RK method to the ODE system in normal form \eqref{sblstnormal}.
When the Jacobian may become singular, special care has to be taken when solving the linear systems involved by the fixed-point algorithm: this difficulty will be discussed in \hl{Section} \ref{ss:resonant}.

\subsection{Local Problem~Solver}

 Let us suppose that the stencils are
$$
\mathcal{S}_i = \{ i-l, \dots, i + r \}.
$$ 
Remember that local problems (LP) consist in, given an index $i$ and a state $\widetilde W \in \Omega$, finding~approximations 
$$
U^{*,m}_{i, j}, \  m=1, \dots, M, \ j = i-l, \dots, i+r; \quad U_{i}^{*, i\pm 1/2},
$$
of the values
$$
U_i^*(x^m_j), \  m=1, \dots, M, \ j \in i-l, \dots, i+r; \quad U_i^*(x_{i\pm 1/2}),
$$
where $U_i^*$ is the stationary solution that satisfies
\begin{equation} 
\sum_{m=1}^M b_m U_i^* (x^m_i)  = \widetilde W.
\end{equation}

 To solve the problem, we look first for an approximation of $U^*_i$ in the cell $I_i$ using the collocation method, i.e.,~we look for a vector valued polynomial $P_i$ of degree $s$ 
 such that:
\begin{equation}\label{nonlinsys1}
 \left\{
 \begin{array}{l}
 \displaystyle \sum_{m=1}^s b_m P_i(x_i^m)  = \widetilde W, \\\\
 \displaystyle D_f(P_i(x_i^m)) P'_i(x_i^m) = S(P_i(x_i^m)) H_x (x_i^m), \quad m = 1,\dots, s. 
 \end{array}
 \right.
 \end{equation}
 Observe that this system has $(s +1)$ vector unknowns (the coefficients of the polynomials) and $(s +1)$ vector~equations. 
 
 A fixed-point algorithm will be used to solve this system. Once it has been solved, the~approximations of $U^*_i$ at the intercells are computed by evaluating the polynomial
 $P_i$:

 $$
 U^{i-1/2} = P_i(x_{i-1/2}), \quad U^{i +1/2} = P_i(x_{i + 1/2}).
 $$
 
From these two values, the~Cauchy solver described in the previous section will be used to compute the approximations of the stationary solution in the other cells of the~stencil.

 Using the standard RK notation, the~local problem solver is as follows:

   \begin{alg} \label{alg:lp} {\hl{Numerical solver} for the local problems (LP) using collocation RK methods.} 
\begin{itemize}

 \item Compute $U^{i-1/2}$, $K_i^m$, $m =1,\dots, s$ by solving the system:
 \begin{adjustwidth}{-4.6cm}{0cm} 
\begin{equation} \label{locsys}
 \left\{
  \begin{array}{l}
\displaystyle D_f\left( U^{i-1/2} + \Delta x \sum_{k=1}^s a_{m,k} K_i^k \right) K_i^m = 
S\left(U^{i-1/2} + \Delta x \sum_{k=1}^s a_{m,k} K_i^k \right) H_x( x_i^m), \quad m= 1, \dots, s, \\
\\
\displaystyle\sum_{m=1}^s b_m \left(U^{i-1/2} + \Delta x \sum_{k=1}^s a_{m,k} K_i^k \right) = \widetilde W.
\end{array}
\right. 
\end{equation}\end{adjustwidth}

\item Compute:
\begin{eqnarray*}
U^{i + 1/2} & = &  U^{i -1/2} + \Delta x \sum_{m=1}^s b_m K_i^m; \\
U_i^m& = & U^{i-1/2} + \Delta x \sum_{k=1}^s a_{m,k} K_i^k, \quad m= 1, \dots, s.
\end{eqnarray*}

\item Apply Algorithm \ref{alg:gp} forward in space from $U^{i+1/2}$ to obtain for $j=i + 1, \dots, i + r$
$$
U^{j + 1/2}; \quad  U_{j}^m, \quad m = 1, \dots, s.
$$

\item Apply Algorithm \ref{alg:gp} backward in space from $U^{i-1/2}$ to obtain for $j=i-l, \dots, i-1$
$$
U^{j-1/2}; \quad  U_{j}^m, \quad m = 1, \dots, s.
$$

\end{itemize}
\end{alg}

The output of the solver is then:
$$
U^{*,m}_{i,j} =  U_j^m, \ m= 1, \dots, s, \ j= i-l, \dots i + r; \quad U_i^{*,i-1/2} = U^{i-1/2}, \quad U_i^{*,i + 1/2}  = U^{i + 1/2}.
$$

A fixed-point algorithm is used again to solve the non-linear system \eqref{locsys}. A~sensible initial guess is given by:
\begin{eqnarray*}
U^{i-1/2, 0} &=& \widetilde W, \\
K_i^{m,0} &=& K, \quad m=1,\dots, s,
\end{eqnarray*}
where $K$ is the solution of the linear system:
$$
D_f(\widetilde W)K = S(\widetilde W)H_x(x_i),
$$
where $x_i$ represents the mid-point of the cell $I_i$.
Let us mention again that special care has to be taken if the Jacobian of the flux function may become singular: this difficulty will be discussed in \hl{Section} \ref{ss:resonant}.

\subsection{Well-Balanced~Property}
To prove that the high-order numerical methods are well-balanced, we have to check that the ODE solver and the local problem solver described in the two last subsections are consistent in the sense of the statement of Theorem \ref{th:wb2}.

Let 
$$U^{*,m}_{\Delta x, i}, \quad m=1, \dots, s; \quad U_{\Delta x}^{*, i+1/2} $$
be the approximations of a stationary solution $U^*$ at the quadrature points and intercells 
$$
U^*(x_i^m), \quad m=1, \dots, s; \quad U^*(x_{i+1/2}),
$$
provided by the ODE solver described in \hl{Section} \ref{ss:cauchy}
in the uniform mesh of nodes $x_{i+1/2} = x_{i-1/2} +   \Delta  x$, $i=i_0,i_0+1,\dots$ .
Let $P^*_i$ be the only polynomial of degree $s$ such that:
\begin{eqnarray*}
 & & P^*_i(x_{i-1/2}) = U_{\Delta x}^{*, i-1/2}, \\
 & & D_f(P^*_i(x_i^m)) \partial_x P^*_i(x_i^m) = S(P^*_i(x_i^m))H_x(x_i^m), \quad i=1, \dots, s.
\end{eqnarray*}
By the double interpretation of the method, we have:
\begin{eqnarray*}
P^*_i(x_i^m) & = &  U^{*,m}_{\Delta x, i}, \quad m=1, \dots, s; \\
P^*_i(x_{i+1/2}) & = & U_{\Delta x}^{*, i+1/2}.
\end{eqnarray*}
To solve now the local problem at the $i$-th cell with
$$
\widetilde W = \widetilde U^*_{\Delta x, i} = \sum_{m=1}^s b_m U^{*,m}_{\Delta x, i},
$$
we have to look for a polynomial $P_i$ satisfying \eqref{nonlinsys1}. Obviously, $P^*_i$ solves this non-linear system; thus,
$P_i = P^*_i$, which implies, in particular,
\begin{eqnarray}
& & U^{*,m}_{i, i} = U^{*,m}_{\Delta x, i}, \quad m=1, \dots, s;\\
& &  U_{i}^{*, i \pm 1/2} = U_{\Delta x}^{*, i \pm 1/2}.
\end{eqnarray}
The equalities
\begin{equation}
U^{*,m}_{i, j} = U^{*,m}_{\Delta x, j}, \quad m=1, \dots, s, \ j =  i + 1, \dots i + r 
\end{equation}
are then trivial, since these values are computed by applying Algorithm \ref{alg:gp} forward in space from the same value $U^*_{\Delta x, i \pm 1/2}$, and the equalities
\begin{equation}
U^{*,m}_{i, j} = U^{*,m}_{\Delta x, j}, \quad m=1, \dots, s, \ j =  i-l, \dots i-1, 
\end{equation}
follow from the symmetry of the collocation method. Therefore, the~hypothesis of \mbox{Theorem \ref{th:wb2}} is satisfied, and the numerical method is then well-balanced \textit{\hl{provided that every stationary solution} $U^*$ can be approximated by the RK collocation method}. 

\subsection{First and Second Order~Methods} \label{ss:fso}

For second order methods, the~MUSCL reconstruction that uses 3-cell centered stencils ($l = r =1$) and
 the second-order 1-stage collocation RK method whose Butcher tableau is
\begin{equation}\label{RKcol2}
\begin{array}{c|c}
1/2 & 1/2 \\\hline
0 & 1\\
\end{array}
 \end{equation}
 will be used. The~solution of local problems is particularly simple, since \eqref{locsys} reduces to~find
 $U^{i-1/2}$, $K$ such that
 $$
  \left\{
  \begin{array}{l}
\displaystyle D_f\left( U^{i-1/2} + \frac{\Delta x}{2}K\right) =  S\left( U^{i-1/2} + \frac{\Delta x}{2}K\right) H_x( x^1_i),  \\
\\
\displaystyle U^{i-1/2} + \frac{\Delta x}{2}K = \widetilde W,
\end{array}
\right. 
$$
where 
$$x^1_i = x_{i-1/2} + \frac{\Delta x}{2} = x_i. $$
Therefore, it is enough to solve
$$
D_f(\widetilde W)K = S(\widetilde W)H_x(x_i)
$$
and then define
$$
U^{i-1/2} = \widetilde W - \frac{\Delta x}{2} K, \quad U^{i + 1/2} = \widetilde W + \frac{\Delta x}{2} K, \quad U^1_0 = \widetilde W.
$$

Then, the~standard method is used forward and bakward in space to compute $U^1_{-1}$ and $U^1_1$. 

Finally, for~first order numerical method, the~trivial pointwise constant reconstruction is considered so that the stencils for the reconstructions consist only of one cell ($l=r=0$). The~second order collocation method \eqref{RKcol2} is also~used.

\begin{Remark}
Notice that, since the 1-stage Gauss-Legendre method is second order, the~schemes considered here~are:
\begin{itemize}
    \item First order numerical methods which are well-balanced with order two.
    \item Second order numerical methods which are well-balanced with order two.
\end{itemize}
\end{Remark}

\subsection{Third-Order~Methods}
The third-order CWENO reconstruction operator will be used here to design third-order methods. This operator uses 3-cell centered stencils, so that $l = r = 1$. 
 The fourth-order 2-stage Gauss-Legendre method corresponding to the choices,
\begin{equation}
    \left[
\begin{array}{cc}
\displaystyle a_{11} & a_{12} \\
\displaystyle a_{21} & a_{22} \\
\end{array}
\right] =  
    \left[
\begin{array}{cc}
\displaystyle \frac{1}{4} & \displaystyle \frac{1}{4} - \frac{\sqrt{3}}{6} \\
\displaystyle \frac{1}{4} + \frac{\sqrt{3}}{6} & \displaystyle \frac{1}{4}  \\
\end{array}
\right],
 \end{equation}
\begin{equation}
  c_1 = \frac{1}{2}-\frac{\sqrt{3}}{6}, \quad c_2 =\frac{1}{2}+\frac{\sqrt{3}}{6},
  \quad b_1 = b_2 = \frac{1}{2},
 \end{equation}
 will be used to solve both the local and global~problems.

\begin{Remark}
Notice that, since the 2-stage Gauss-Legendre method is fourth order, we consider third order numerical methods which are well-balanced with order four.
\end{Remark}

\subsection{Resonant~Problems}\label{ss:resonant}
 When a fixed-point algorithm is used to solve \eqref{globsys1}, \eqref{globsys2}, or~\eqref{locsys}, one or more linear systems of the form
\begin{equation}\label{linearsys}
D_f(W) K = S(W)H_x(\bar x)
\end{equation}
have to be solved at every step, where $W$ is a known state, and $\bar x$ a quadrature point. If~ $W$ is a sonic state, i.e.,~if one of the eigenvalues of $D_f(W)$ vanishes, the~problem is said to be resonant and the system may not have solution or it may have infinitely many~ones:
\begin{enumerate}
\item If $S(W)H_x(x)$ does not belong to the image of the linear application defined by the matrix
$D_f(W)$, the~system has no~solution. 

\item Otherwise, the~system has infinitely many solutions
$$ K^* + \alpha R, \quad \alpha \in \mathbb{R},$$
where $K^*$ is a particular solution, and $R$ is an eigenvector associated to the null eigenvalue.
\end{enumerate}

Since, when applying RK collocation methods, we are looking for solutions of the differential Equation \eqref{sblst}, a~solution $K$ of \eqref{linearsys} is said to be \emph{admissible} if there exists a $C^1$ stationary solution $U^*$ such that
$$
U^*(\bar x) = W, \quad U_x^*(\bar x) = K.
$$

In case 1, there are no admissible solutions. If~this situation arises when solving a local problem, it will be assumed that the problem has no solution and the standard reconstruction will be used in the corresponding~cell. 

In case 2, if~$K$ is an admissible solution and $W$ is an isolated sonic point of the corresponding smooth stationary solution $U^*$, the~following equality is satisfied:
$$
U^*_x(x) = D_f(U^*)^{-1}S(U^*)H_x(x), \quad \forall x \in (\bar x - \varepsilon, \bar x + \varepsilon) - \{ \bar x \},
$$
for some $\varepsilon > 0$.
Therefore, the~limit
\begin{equation}\label{limit}
\lim_{x \to \bar x} D_f(U^*)^{-1}S(U^*)H_x(x)
\end{equation}
has to exist. If~the available information about the PDE system and its stationary solution allows one to compute this limit by any analytic procedure, then the possible values of \eqref{limit} constitute the set of admissible solutions, and the solution of \eqref{linearsys} will then be selected in this set. 
This technique will be illustrated in the particular case of the shallow water equations in \hl{Section} \ref{Subsec:sw_transcritical}.

\begin{Remark}
 Observe that, if~$W$ is sonic and $S(W)$ does not belong to the image of $D_f(W)$, case 2 can only happen if $H_x(\bar x)= 0$. 
In this case, if~there exists a $C^1$ stationary solution $U^*$ such that $U^*(\bar x) = W$, then $U^*(\bar x)$ would be an admissible solution of \eqref{linearsys} which implies that necessarily $H_x(x^*) = 0$, i.e.,~a smooth stationary solution of \eqref{sle} can only reach a sonic state $W$ such that $S(W)$ does not belong to the image of $D_f(W)$ at a critical point of $H$. In~this case,
\eqref{limit} will imply a $0/0$ indeterminate form that can be solved by L'H\^opital rule. This is the case for many systems, including the shallow water equations, for~which a smooth stationary solution can only reach a sonic (critical) state in the minimum points of the depth function $H$. 
\end{Remark}

\section{Numerical~Experiments} \label{numerical_experiments}
Let us describe first the choices that have been made to implement the well-balanced scheme previously~described:
\begin{itemize}
\item The second-order 1-stage collocation RK method is chosen for the first and second order well-balanced numerical schemes and the 2-stage Gauss-Legendre collocation method for third order ones.
\item Concerning the quadrature rules, the~midpoint rule for first and second order schemes are used, and~the 2-point Gauss quadrature rule for third order schemes.
\end{itemize}

In all the numerical tests, Rusanov numerical flux is chosen and MUSCL (see Reference~\cite{MUSCL}) and CWENO reconstructions operators (see Reference~\cite{CWENO,CWENO2}) are used, respectively, for the second and third order schemes. Finally, TVD Runge-Kutta schemes of first, second, and third~\cite{Gottlieb98} are used to integrate the ODE system \eqref{met_num}--\eqref{Siquadbis}.

In those problems where the initial condition is a steady state $U^*$, the~approximation of its cell averages are computed using a quadrature formula if the exact solution is known, or~by approximating the solution of the Cauchy problem \eqref{sblst} with initial condition
$$
U(a) = U^*(a),
$$
using the RK collocation method: see Remark \ref{rem:wpic}.

Note that the procedure that we have been describing in this work is quite general and works for a wide class of $f$, $S$, and $H$.

In this section, we will denote the different numerical methods that we compare as~follows:

\begin{itemize}

\item NWBM$i$, $i=1, 2, 3$: a standard non-well balanced numerical scheme of order $i$. 

\item WBM$i$, $i=1, 2, 3$: well-balanced scheme of order $i$ with reconstruction operator given by Algorithm \ref{alg:ewbnirec}.

\item CWBM$i$, $i=1, 2, 3$: well-balanced scheme of order $i$ with reconstruction operator \ref{alg:wbrec} in which the Cauchy problems are solved following Reference~\cite{GomezCastroPares2020}.

\item CLWBM$i$, $i=1, 2, 3$: well-balanced scheme of order $i$ with reconstruction operator defined by Algorithm \ref{alg:wbrec} in which Cauchy and local problems are integrated by using the Gauss-Legendre collocation method described in the previous~section.

\end{itemize}

The numerical tests have been run on a laptop with Intel(R) Xeon(R) CPU E3-1220 v3 @ 3.10~GHz with 8~Mb cache using one single~core.

\subsection{Problem 1: Burgers Equation with a Non-linear Source Term~I}

We consider the Burgers equation with a non-linear source term:
\begin{equation}\label{burgers}
\begin{cases}
\displaystyle u_t + \left( \frac{u^2}{2} \right)_x= u^2 , \quad x\in \mathbb{R}, \, t>0,\\
u(x,0)=u_0(x).
\end{cases}
\end{equation}
This problem corresponds to the choice
$$U = u, \ f(U)= \displaystyle \frac{u^2}{2}, \ S(U)=u^2, \quad H(x)=x.$$
The stationary solutions satisfy the ODE
\begin{equation}\label{ODEBurgers1}
u' = u,
\end{equation}
whose general solution is 
$$u(x)=C e^x, \quad C \in \mathbb{R}.$$

As the expression of the steady states is known, the~stationary solution $u^*_i$ satisfying \eqref{step1qf} could be easily computed for any quadrature formula
\begin{equation}\label{expresionpaso1}
u^*_i(x)= \frac{u_i}{\sum_{m=0}^M \alpha_m^i e^{x_m^i}} e^x.
\end{equation}

Although in this test WBM$i$ can be easily computed, we also consider CWBM$i$ and CLWBM$i$ to compare the efficiency of the~methods. 

\subsubsection*{Test 1.1}\label{test11}
Let us consider the space interval $[-1, 1]$ and the time interval $[0,5]$. $u_0=e^x$ is set as initial condition and the CFL parameter is set to 0.9. The~boundary condition is imposed at the left boundary and open boundary conditions at the right~boundary. 

 The convergence criterion chosen in Algorithm \ref{alg:lp} is the following: 
we stop it if the maximum absolute value between the approximations of the values of the stationary solution at the extremes of the cell and the quadrature points of the cells belonging to the stencil for two consecutive iterations is less than a~threshold.

In what follows, the~tolerance considered to stop Algorithm \ref{alg:lp} is $\varepsilon=10^{-15}$.

Note that, at the discrete level, the~cell-averages of $u_0(x)=e^x$ have to be computed. We consider $4$ different ways of computing~them:
\begin{itemize}
    \item[(a)] Exactly computed:
    $$u^0_i = \frac{e^{\Delta x} - 1}{\Delta x } e^{x_{i-1/2}}, \quad \forall i.
    $$
    \item[(b)] Using a quadrature formula:
    $$
    u^0_i = \sum_{m=0}^M \alpha_m^i e^{x_m^i}, \quad \forall i.
    $$
    \item[(c)] The cell averages are approximated by
\begin{equation} \label{test1:ic}
u^0_i = \sum_{m = 0}^M \alpha^i_m u^{*,i}_{h,m},\quad \forall i,
\end{equation}
where $u^{*,i}_{h,m}$ are the approximation the Cauchy problem
\begin{equation}\label{test1:cauchy}
    \begin{cases}
    u'=u,\\
    u(-1)=e^{-1},
    \end{cases}
    \end{equation}
computed with the standard RK4 method on the mesh selected to solve the control problems described in Reference~\cite{GomezCastroPares2020}.

\item[(d)] The cell averages are approximated by \eqref{test1:ic}, where now 
 $u^{*,i}_{h,m}$ are the approximation of the exponential at the quadrature points obtained by solving \eqref{test1:cauchy} with the RK collocation~method.

\end{itemize}

\begin{figure}[H]
{\captionsetup{position=bottom,justification=centering}
\hspace{-0.5cm}
  \subfloat[ NWBM$i$, $i =1,2,3$ ]{
   \includegraphics[width=0.38\textwidth]{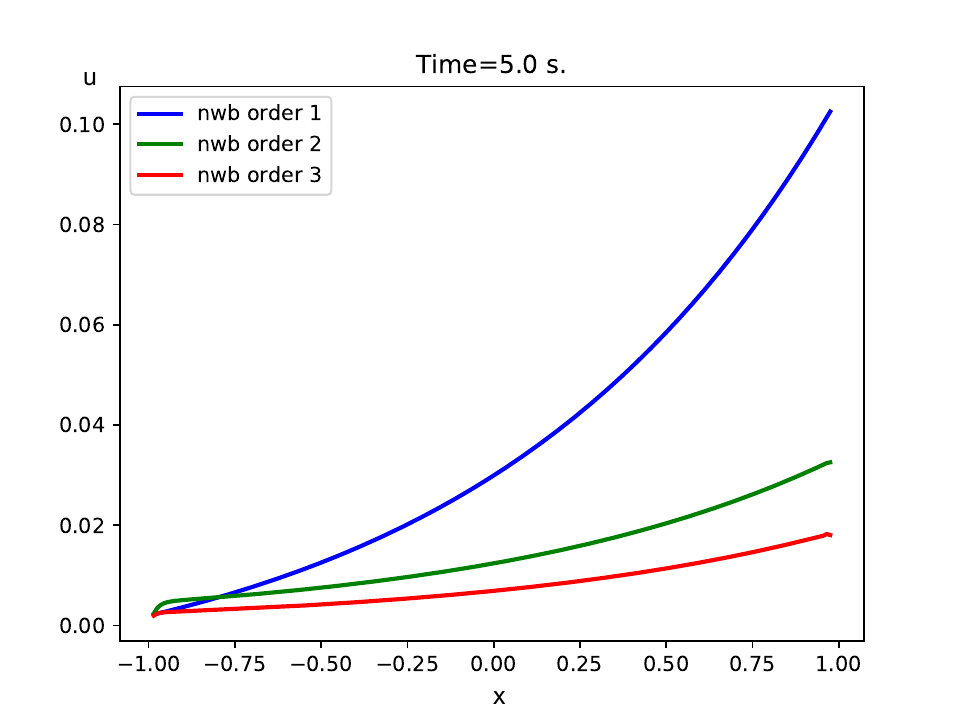}}
 \subfloat[CLWBM$i$, $i = 1,2,3$]{
   \includegraphics[width=0.38\textwidth]{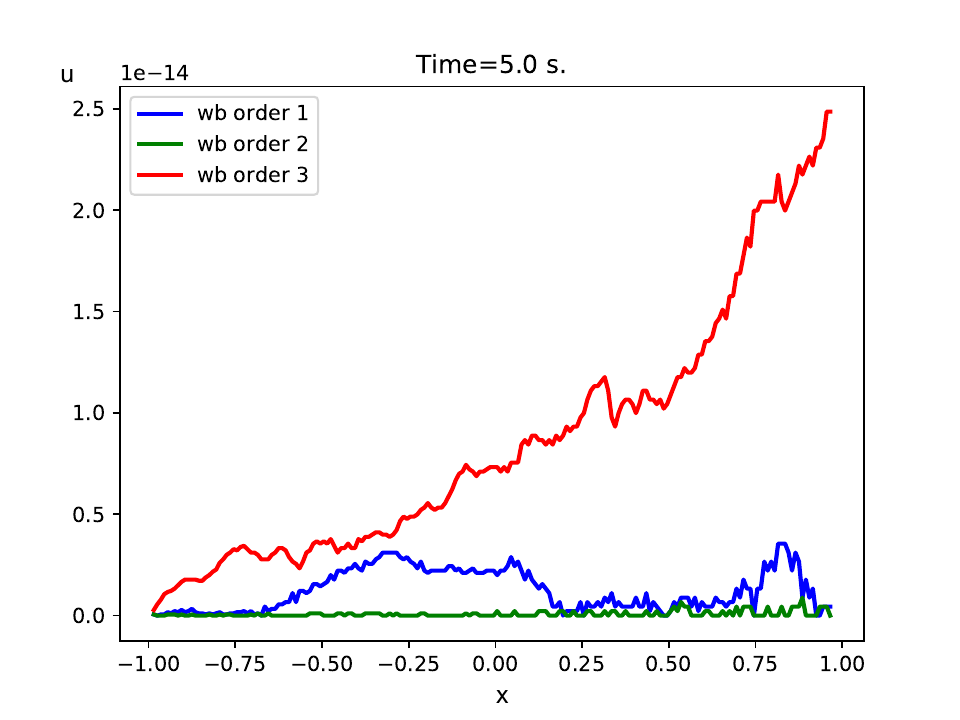}}}
    \caption{\hl{Test 1.1. Differences} between the stationary and the numerical solutions at time $t = 5$~s for a $200-$cell~mesh.} \label{fluc_test1_1}
\end{figure}

The numerical methods are integrated until the final time $t = 5$, which is more than twice the travel time of a small perturbation of the steady state along the domain. $L^1$ errors between the initial and final cell-averages will be computed. According to Theorem \ref{th:wb1}, errors should be of the order of the machine precision (MP) if option (d) is considered to compute the initial condition for CLWBM$i$, $i=1,2,3$. Figure~\ref{fluc_test1_1} (right) shows the errors corresponding to CLWBM$i$, $i = 1,2,3$ for a 200-cell mesh. As~expected, they are of the order of~MP.

Notice that, in the other cases, if~(a), (b), or (c) are used for CLWBM$i$, $i = 1,2,3$, the~errors decrease with $\Delta x$. Thus, if~option (b) is chosen and a third order reconstruction operator is used, then the convergence rate is expected to be $4$ as errors are generated by the use of the 2-stage collocation RK method for solving the local problems. If~option (c) is used, the~same convergence rate is expected, and the errors are due to the use of RK4 for computing the initial condition and the use of the 2-stage collocation RK method.
Finally, in~the case of option (a), the errors are due to the use of both numerical integration and the 2-stage collocation method in the solution of local problems: again, the errors are expected to be of order~4. 

Similar results are obtained for first and second-order methods: the order of convergence is two if option (a), (b), and (c) are used since the mid-point method is second-order accurate, and we consider the 1-stage collocation RK~method.

\begin{specialtable}[H]
 \caption{\hl{Test 1.1.} $L^1$ errors at $t = 5$ and convergence rates for CLWBM3 when option (a), (b), (c), and (d) are chosen to compute the initial cell~averages.}
 \label{test11_orden3_condicionesiniciales}
 \begin{adjustbox}{width=\linewidth,left} 
\begin{tabular}{ccccccccc} \toprule
\multirow{2}{*}{\textbf{Cells}}& \multicolumn{2}{c}{\textbf{Option (a)}} &\multicolumn{2}{c}{\textbf{Option (b)}} & \multicolumn{2}{c}{\textbf{Option (c)}} & \multicolumn{2}{c}{\textbf{Option (d)}} \\
 &  \textbf{Error} & \textbf{Order}&  \textbf{Error} & \textbf{Order}&  \textbf{Error} & \textbf{Order} &  \textbf{Error} & \textbf{Order} \\ \midrule
5& 1.96 $\times$ $10^{-4}$ & - & 1.82 $\times$ $10^{-4}$  & - & 1.46 $\times$ $10^{-4}$ & - & 4.88 $\times$ $10^{-16}$ & -\\
10 & 1.21 $\times$ $10^{-5}$ & 4.015 & 1.12 $\times$ $10^{-5}$ & 4.016 & 8.81 $\times$ $10^{-6}$ & 4.054 & 7.44 $\times$ $10^{-16}$ & -\\
20 & 7.56 $\times$ $10^{-7}$ & 4.003 & 7.01 $\times$ $10^{-7}$ & 4.003 & 5.42 $\times$ $10^{-7}$ & 4.023 & 3.00 $\times$ $10^{-16}$ & -\\
40 & 4.72 $\times$ $10^{-8}$ & 4.001 & 4.38 $\times$ $10^{-8}$ & 4.001 & 3.36 $\times$ $10^{-8}$ & 4.011 & 7.94 $\times$ $10^{-15}$ & -\\
80 & 2.95 $\times$ $10^{-9}$ & 4.000 & 2.74 $\times$ $10^{-9}$ & 4.000 & 2.09 $\times$ $10^{-9}$ & 4.006 & 9.00 $\times$ $10^{-15}$ & -\\
160 & 1.84 $\times$ $10^{-10}$ & 4.000 & 1.71 $\times$ $10^{-10}$ & 4.000 & 1.30 $\times$ $10^{-10}$ & 4.003 & 1.23 $\times$ $10^{-14}$ & -\\\bottomrule
\end{tabular}
\end{adjustbox} 
\end{specialtable}

Now, we compare the results obtained with NWBM$i$, WBM$i$, CWBM$i$, CLWBM$i$, $i=1,2,3$. The~cell-averages are computed using option (b) for NWBM$i$ and WBM$i$, option (c) for CWBM$i$, and option (d) for CLWBM$i$. $L^1$ errors with respect to the initial condition are~computed.
\vspace{-20pt}

Figure~\ref{fluc_test1_1} shows the differences between the stationary and the numerical solutions obtained with NWBM$i$ and CLWBM$i$, $i=1,2,3$. Similar results are obtanied for WBM$i$ and CWBM$i$, $i=1,2,3$.
\hl{Tables}~\ref{ex1_error_nwb} and \ref{ex1_error_wb} show the errors corresponding to NWBM$i$, WBM$i$, CWBM$i$, CLWBM$i$, $i = 1,2,3$. 
\begin{specialtable}[H]
\caption{Test 1.1. Differences in $L^1$-norm with respect to the stationary solution and convergence rates for NWBM$i$, $i= 1,2,3$.} \label{ex1_error_nwb}
\begin{adjustbox}{width=\linewidth,left} 
\begin{tabular}{ccccccc} \toprule
\textbf{Cells} & \textbf{NWBM1: Error} &\textbf{Order}& \textbf{NWBM2: Error} & \textbf{Order}&\textbf{NWBM3: Error}&\textbf{Order}\\\midrule
100& 7.53 $\times$ $10^{-2}$ & - & 2.44 $\times$ $10^{-3}$  & - & 7.66 $\times$ $10^{-6}$ & -\\
200 & 3.78 $\times$ $10^{-2}$ & 0.995 & 8.09 $\times$ $10^{-4}$ & 1.591 & 9.62 $\times$ $10^{-7}$ & 2.993\\
400 & 1.89 $\times$ $10^{-2}$ & 1.002 & 2.16 $\times$ $10^{-4}$ & 1.905 & 1.21 $\times$ $10^{-7}$ & 2.995\\
800 & 9.43 $\times$ $10^{-3}$ & 1.000 & 5.54 $\times$ $10^{-5}$ & 1.963 & 1.51 $\times$ $10^{-8}$ & 2.998\\ \bottomrule
\end{tabular}
\end{adjustbox} 
\end{specialtable}
\unskip

\end{paracol}
\nointerlineskip
\begin{specialtable}[H]
\widetable
\caption{Test 1.1. Differences in $L^1$-norm with respect to the stationary solution for WBM$i$, CWBM$i$, CLWBM$i$, $i= 1,2,3$.} \label{ex1_error_wb}
\begin{adjustbox}{width=\linewidth,left} 
\begin{tabular}{cccccccccc} \toprule
\multirow{2}{*}{\textbf{Cells}} & \multicolumn{3}{c}{\textbf{Order 1: Error}} & \multicolumn{3}{c}{\textbf{Order 2: Error}} & \multicolumn{3}{c}{\textbf{Order 3: Error}} \\

 & \textbf{WBM} & \textbf{CWBM} &\textbf{CLWBM} & \textbf{WBM} & \textbf{CWBM} &\textbf{CLWBM} & \textbf{WBM} & \textbf{CWBM} &\textbf{CLWBM} \\\midrule

100& 4.21 $\times$ $10^{-15}$ & 3.55 $\times$ $10^{-15}$ & 2.50 $\times$ $10^{-15}$ & 8.87 $\times$ $10^{-16}$ & 3.63 $\times$ $10^{-16}$ & 2.03 $\times$ $10^{-16}$ & 3.20 $\times$ $10^{-16}$ & 1.43 $\times$ $10^{-14}$ & 8.17 $\times$ $10^{-15}$ \\
200 & 2.90 $\times$ $10^{-15}$ & 5.54 $\times$ $10^{-13}$ & 2.51 $\times$ $10^{-15}$ & 4.42 $\times$ $10^{-16}$ & 1.23 $\times$ $10^{-15}$ & 1.66 $\times$ $10^{-16}$ & 2.54 $\times$ $10^{-16}$ & 2.43 $\times$ $10^{-14}$ & 1.76 $\times$ $10^{-14}$ \\
400 & 1.84 $\times$ $10^{-14}$ & 2.05 $\times$ $10^{-14}$ & 1.12 $\times$ $10^{-15}$ & 1.82 $\times$ $10^{-15}$ & 3.64 $\times$ $10^{-16}$ & 2.89 $\times$ $10^{-16}$  & 7.40 $\times$ $10^{-14}$ & 4.47 $\times$ $10^{-14}$ & 4.45 $\times$ $10^{-14}$\\
800 & 4.45 $\times$ $10^{-16}$ & 2.67 $\times$ $10^{-15}$ &  2.77 $\times$ $10^{-15}$ & 1.83 $\times$ $10^{-16}$ & 2.03 $\times$ $10^{-16}$ & 
2.05 $\times$ $10^{-16}$ & 2.61 $\times$ $10^{-15}$ &  9.48 $\times$ $10^{-14}$ & 
  7.88 $\times$ $10^{-14}$ \\  \bottomrule
\end{tabular}
\end{adjustbox} 
\end{specialtable}
\begin{paracol}{2}
\switchcolumn

Notice that the errors for NWBM decrease with $\Delta x$ at the expected rate. WBM preserve the exact solution up to MP, as well as CWBM and CLWBM, provided that the tolerance used in the Newton's method for the former and Algorithm 
\ref{alg:lp} for the latter is small enough. The~computational costs are shown in Table~\ref{ex1_times}. Nonetheless, since the goal of carrying out these experiments is to prove that the introduced well-balanced technique works, the~implementation is not optimized. It can be seen that the well-balanced reconstruction given by Algorithm \ref{alg:ewbnirec} increases the computational cost by a factor ranging from 1.5 to 7.5 (see Table~\ref{ex1_times}). On~the other hand, the~use of Algorithm \ref{alg:wbrec} increases the computational cost by a factor ranging from 1 to 1.8
if collocation RK method are used for solving the local problems and a factor ranging from 1.1 to 3.8 if control techniques with $N_p = 1$ (see Reference~\cite{GomezCastroPares2020}) are used. Therefore, the~CLWBM schemes introduced in this work are more efficient than the CWBM schemes presented in Reference~\cite{GomezCastroPares2020}.

\begin{specialtable}[H]
\caption{Test 1.1. Computational times (milliseconds).}\label{ex1_times}

\setlength{\cellWidtha}{\columnwidth/6-2\tabcolsep-0in}
\setlength{\cellWidthb}{\columnwidth/6-2\tabcolsep+0in}
\setlength{\cellWidthc}{\columnwidth/6-2\tabcolsep-0.0in}
\setlength{\cellWidthd}{\columnwidth/6-2\tabcolsep-0.0in}
\setlength{\cellWidthe}{\columnwidth/6-2\tabcolsep-0.0in}
\setlength{\cellWidthf}{\columnwidth/6-2\tabcolsep-0.0in}

\scalebox{1}[1]{\begin{tabularx}{\columnwidth}{>{\PreserveBackslash\centering}m{\cellWidtha}>{\PreserveBackslash\centering}m{\cellWidthb}>{\PreserveBackslash\centering}m{\cellWidthc}>{\PreserveBackslash\centering}m{\cellWidthd}>{\PreserveBackslash\centering}m{\cellWidthe}>{\PreserveBackslash\centering}m{\cellWidthf}}
\toprule 
     \textbf{Cells }             & \boldmath{\textbf{Order ($i$)} }& \boldmath{\textbf{NWBM$i$}} & \boldmath{\textbf{WBM$i$} }& \boldmath{\textbf{CWBM$i$}  }&    \boldmath{ \textbf{CLWBM$i$}}     \\  \midrule
\multirow{3}{*}{100} & 1 & 20 & 30 & 70& 30\\ 
                  & 2 & 30 & 60 & 140 & 110 \\ 
                  & 3 & 40 & 190 & 240 & 230 \\ \midrule
\multirow{3}{*}{200} & 1 & 20 & 60 & 230& 100 \\ 
                  & 2 & 40 & 190 & 330 & 280 \\ 
                  & 3 & 110 & 480 & 530& 520 \\ \midrule         
\multirow{3}{*}{400} & 1& 50 & 180 & 520 & 220 \\ 
                  & 2 & 100 & 530 &1150&  720 \\ 
                  & 3 & 350 & 1680 & 1980&  1870 \\ \midrule
\multirow{3}{*}{800} & 1 & 140 & 570 & 2020& 870 \\ 
                  & 2 & 270 & 2040 & 3580& 2820 \\ 
                  & 3& 1080 & 5540 & 6600& 5960 \\ \bottomrule                  
\end{tabularx}}

\end{specialtable}

\subsection{Problem 2: Burgers Equation with a Non-linear Source Term~II}
We now consider Burgers equation with another non-linear source term:
\begin{equation} \label{burgers2}
\begin{cases}
u_t +  \left( \displaystyle \frac{u^2}{2} \right)_x= \sin(u), \quad x\in \mathbb{R}, \, t>0,\\
u(x,0)=u_0(x).
\end{cases}
\end{equation}
The stationary solutions satisfy the ODE
\begin{equation}\label{ODEBurgers2}
u' = \frac{ \sin(u)}{u}.
\end{equation}

Note that no simple expression of the stationary solutions is known so \eqref{step1} has to be computed~numerically.

\subsubsection{Test~2.1}
We consider again the same computational domain, and~the system is integrated up to $T=5$. The~CFL parameter is imposed to be $0.9$. As~initial condition, we consider the stationary solution that corresponds to the solution of the Cauchy problem \eqref{ODEBurgers2} with 
$$
u(-1) = 2,
$$
that is approximated by the Gauss-Legendre collocation method. The~left boundary condition is set to $u(-1,t)=2$, and open boundary conditions are considered~upstream. 

Figure~\ref{fluc_test21} shows the differences between the stationary and the numerical solutions computed with NWBM$i$, $i=1,2,3$ (left) and CLWBM$i$, $i=1,2,3$ (right), (the figures corresponding to CWBM$i$, $i=1,2,3$ are similar).
\vspace{-30pt}

\begin{figure}[H]
{\captionsetup{position=bottom,justification=centering}
  \subfloat[NWBM$i$, $i=1,2,3$]{
   \includegraphics[width=0.38\textwidth]{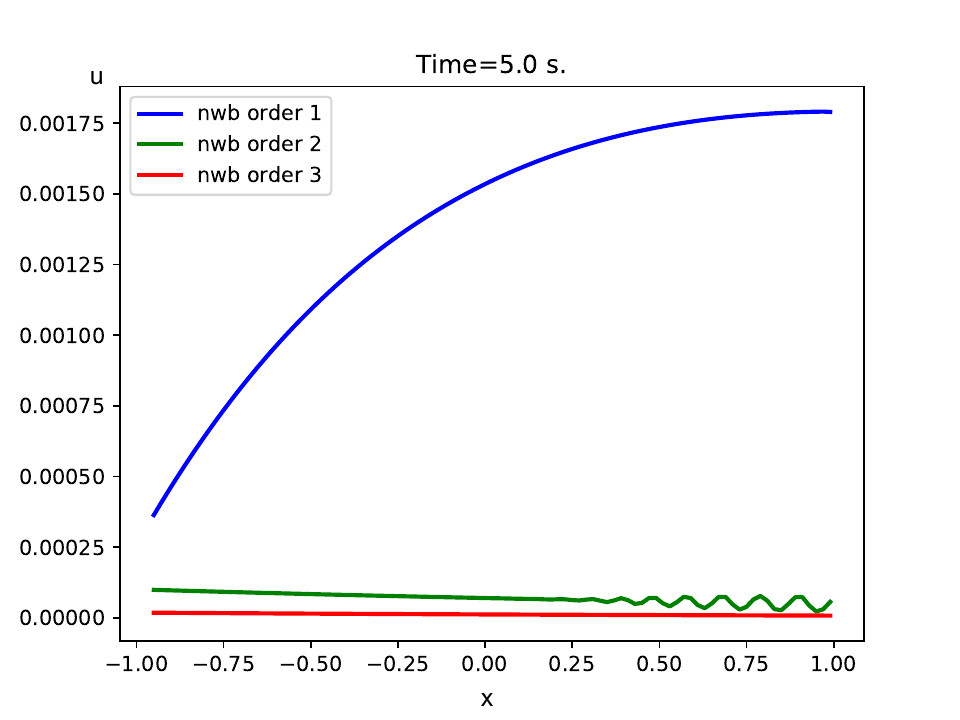}}
 \subfloat[CLWBM$i$, $i=1,2,3$ ]{
   \includegraphics[width=0.38\textwidth]{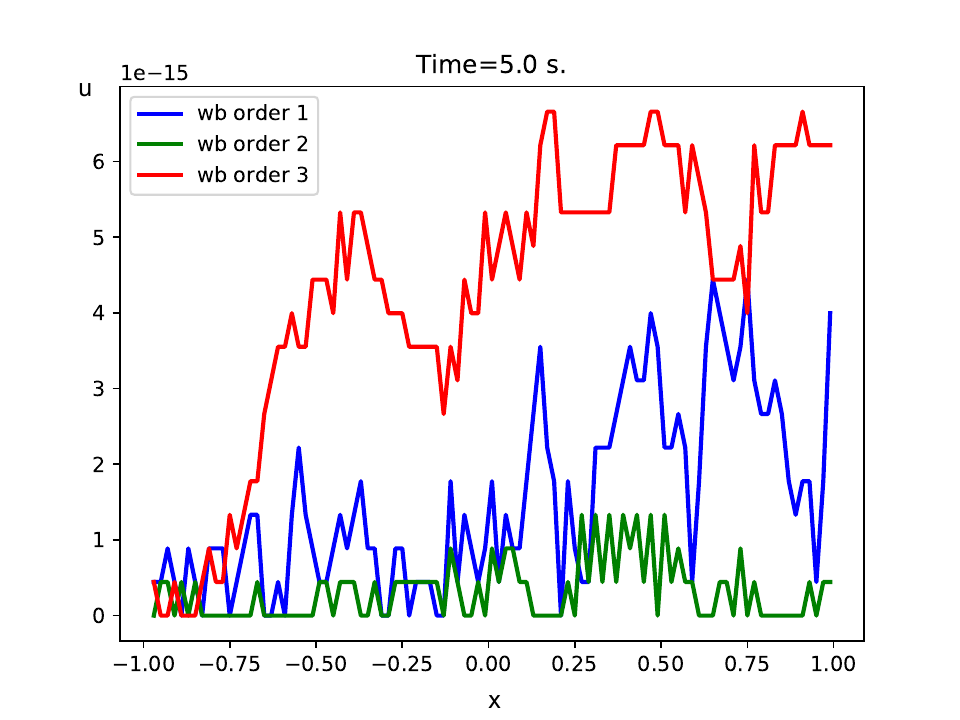}}}
    \caption{\hl{Test 2.1. Differences} between the stationary and the numerical solutions at $t = 5$~s using a mesh with $100$ cells.} \label{fluc_test21}
\end{figure}

The errors corresponding to NWBM$i$, CWBM$i$, and CLWBM$i$, $i=1,2,3$, respectively, are shown in Tables~\ref{test21_error_nwb} and \ref{test21_error_wb_col}. Computational costs are shown in Table~\ref{test21_times}. Note that the collocation strategy increases the computational cost from 7.3 to 14, whereas the method proposed in Reference~\cite{GomezCastroPares2020} increases the computational cost by a factor ranging from 27 to 40. After~the study of these two experiments, one can conclude that the results obtained with the collocation strategy introduced in this paper are better than the ones got by CWBM$i$, $i=1,2,3$ as described in Reference~\cite{GomezCastroPares2020}, considering that the errors are similar for both procedures, but the first strategy is much more efficient. Consequently, we will only consider the collocation strategy in the following~experiments.

\begin{specialtable}[H]
\caption{Test 2.1. Differences in $L^1$-norm with respect to the stationary solution and convergence rates for NWBM$i$, $i= 1,2,3$.} \label{test21_error_nwb}
\setlength{\cellWidtha}{\columnwidth/7-2\tabcolsep-0.2in}
\setlength{\cellWidthb}{\columnwidth/7-2\tabcolsep+0.2in}
\setlength{\cellWidthc}{\columnwidth/7-2\tabcolsep-0.1in}
\setlength{\cellWidthd}{\columnwidth/7-2\tabcolsep+0.2in}
\setlength{\cellWidthe}{\columnwidth/7-2\tabcolsep-0.1in}
\setlength{\cellWidthf}{\columnwidth/7-2\tabcolsep+0.2in}
\setlength{\cellWidthg}{\columnwidth/7-2\tabcolsep-0.2in}

\scalebox{1}[1]{\begin{tabularx}{\columnwidth}{>{\PreserveBackslash\centering}m{\cellWidtha}>{\PreserveBackslash\centering}m{\cellWidthb}>{\PreserveBackslash\centering}m{\cellWidthc}>{\PreserveBackslash\centering}m{\cellWidthd}>{\PreserveBackslash\centering}m{\cellWidthe}>{\PreserveBackslash\centering}m{\cellWidthf}>{\PreserveBackslash\centering}m{\cellWidthg}}
\toprule

\textbf{Cells} & \boldmath{\textbf{Error ($i=1$)}}&\textbf{Order}&\boldmath{\textbf{Error ($i=2$)}}& \textbf{Order}&\boldmath{\textbf{Error ($i=3$)}}&\textbf{Order}\\\midrule
100& 2.72 $\times$ $10^{-3}$ & - & 1.43 $\times$ $10^{-4}$ & - & 2.53 $\times$ $10^{-5}$ & -\\
200 & 1.34 $\times$ $10^{-3}$ & 1.021 & 2.43 $\times$ $10^{-6}$ & 5.879 & 1.74 $\times$ $10^{-8}$ & 10.503\\
400 & 6.58 $\times$ $10^{-4}$ & 1.026 & 8.19 $\times$ $10^{-7}$ & 1.569 & 1.14 $\times$ $10^{-10}$ & 7.250\\
800 & 3.24 $\times$ $10^{-4}$ & 1.022 & 2.34 $\times$ $10^{-7}$ & 1.806 & 1.41 $\times$ $10^{-11}$ & 3.016\\ \bottomrule
\end{tabularx}}
\end{specialtable}
\unskip

\end{paracol}
\nointerlineskip
\begin{specialtable}[H]
\widetable
\caption{Test 2.1. Differences in $L^1$-norm with respect to the stationary solution for CWBM$i$, CLWBM$i$, $i= 1,2,3$.} \label{test21_error_wb_col}
\setlength{\cellWidtha}{\columnwidth/7-2\tabcolsep-0.2in}
\setlength{\cellWidthb}{\columnwidth/7-2\tabcolsep+0.2in}
\setlength{\cellWidthc}{\columnwidth/7-2\tabcolsep-0.1in}
\setlength{\cellWidthd}{\columnwidth/7-2\tabcolsep+0.2in}
\setlength{\cellWidthe}{\columnwidth/7-2\tabcolsep-0.1in}
\setlength{\cellWidthf}{\columnwidth/7-2\tabcolsep+0.2in}
\setlength{\cellWidthg}{\columnwidth/7-2\tabcolsep-0.2in}

\scalebox{1}[1]{\begin{tabularx}{\columnwidth}{>{\PreserveBackslash\centering}m{\cellWidtha}>{\PreserveBackslash\centering}m{\cellWidthb}>{\PreserveBackslash\centering}m{\cellWidthc}>{\PreserveBackslash\centering}m{\cellWidthd}>{\PreserveBackslash\centering}m{\cellWidthe}>{\PreserveBackslash\centering}m{\cellWidthf}>{\PreserveBackslash\centering}m{\cellWidthg}}
\toprule

\multirow{2}{*}{\textbf{Cells}} & \multicolumn{2}{c}{\textbf{Order 1: Error}} & \multicolumn{2}{c}{\textbf{Order 2: Error}} & \multicolumn{2}{c}{\textbf{Order 3: Error}} \\

  & \textbf{CWBM} &\textbf{CLWBM} &  \textbf{CWBM} &\textbf{CLWBM}  & \textbf{CWBM} &\textbf{CLWBM} \\\midrule
100& 9.71 $\times$ $10^{-14}$ & 3.00 $\times$ $10^{-15}$ & 1.76 $\times$ $10^{-13}$ & 6.39 $\times$ $10^{-16}$  & 1.99 $\times$ $10^{-13}$ & 8.50 $\times$ $10^{-15}$\\
200 & 7.56 $\times$ $10^{-15}$ & 5.37 $\times$ $10^{-15}$ & 3.46 $\times$ $10^{-15}$ & 5.15 $\times$ $10^{-16}$ &  2.97 $\times$ $10^{-14}$ & 2.51 $\times$ $10^{-14}$\\
400 & 4.00 $\times$ $10^{-15}$ & 5.68 $\times$ $10^{-15}$ &  7.53 $\times$ $10^{-16}$ & 5.73 $\times$ $10^{-16}$ &  3.31 $\times$ $10^{-14}$ & 4.85 $\times$ $10^{-14}$\\ 
800 & 5.97 $\times$ $10^{-15}$ & 4.63 $\times$ $10^{-15}$ &  8.54 $\times$ $10^{-16}$ & 5.31 $\times$ $10^{-16}$ & 6.63 $\times$ $10^{-14}$ &  9.61 $\times$ $10^{-14}$\\ \bottomrule
\end{tabularx}}
\end{specialtable}
\begin{paracol}{2}
\switchcolumn
\vspace{-10pt}

\begin{specialtable}[H]
\caption{Test 2.1. Computational cost (milliseconds). $t=5$~s.}\label{test21_times}
\setlength{\cellWidtha}{\columnwidth/5-2\tabcolsep-0in}
\setlength{\cellWidthb}{\columnwidth/5-2\tabcolsep+0in}
\setlength{\cellWidthc}{\columnwidth/5-2\tabcolsep-0.0in}
\setlength{\cellWidthd}{\columnwidth/5-2\tabcolsep-0.0in}
\setlength{\cellWidthe}{\columnwidth/5-2\tabcolsep-0.0in}

\scalebox{1}[1]{\begin{tabularx}{\columnwidth}{>{\PreserveBackslash\centering}m{\cellWidtha}>{\PreserveBackslash\centering}m{\cellWidthb}>{\PreserveBackslash\centering}m{\cellWidthc}>{\PreserveBackslash\centering}m{\cellWidthd}>{\PreserveBackslash\centering}m{\cellWidthe}>{\PreserveBackslash\centering}m{\cellWidthf}}
\toprule 
          \textbf{Cells}        & \boldmath{\textbf{Order($i$)}} & \boldmath{\textbf{NWBM$i$} }& \boldmath{\textbf{CWBM$i$} }& \boldmath{\textbf{CLWBM$i$} }\\
 \midrule
\multirow{3}{*}{100} & 1& 10 & 340 & 130 \\ 
                  & 2 & 20 & 690 & 280 \\ 
                  & 3& 40 & 1390 & 490 \\ \midrule
\multirow{3}{*}{200} & 1& 30 & 1280 & 230  \\ 
                  & 2& 60 & 2350 & 440 \\ 
                  & 3& 180 & 5190 & 1230\\ \bottomrule

\end{tabularx}}
\end{specialtable}

\subsubsection{Test~2.2}
Now, we consider the evolution of a small perturbation of the previous stationary solution. Thus, the~initial solution is set to 
$$u_0(x)=u^*(x)+0.3e^{-200(x+0.5)^2},$$
where $u^*(x)$ is the previous stationary~solution.

Figure~\ref{fluc_test22} shows the difference between the stationary and the numerical solutions computed with NWBM$i$, $i=1,2,3$ and CLWBM$i$, $i=1,2,3$ at times $t=0.3, 5$ (the results are similar for CWBM$i$, $i=1,2,3$). In~black, we show a reference solution computed on a fine mesh with 12,800 cells using a first order well-balanced~method.
\vspace{-20pt}

\begin{figure}[H]

{\captionsetup{position=bottom,justification=centering}
\hspace{-0.3cm}
  \subfloat[NWBM$i$, $i =1,2,3$. $t = 0.3$~s. ]{
   \includegraphics[width=0.38\textwidth]{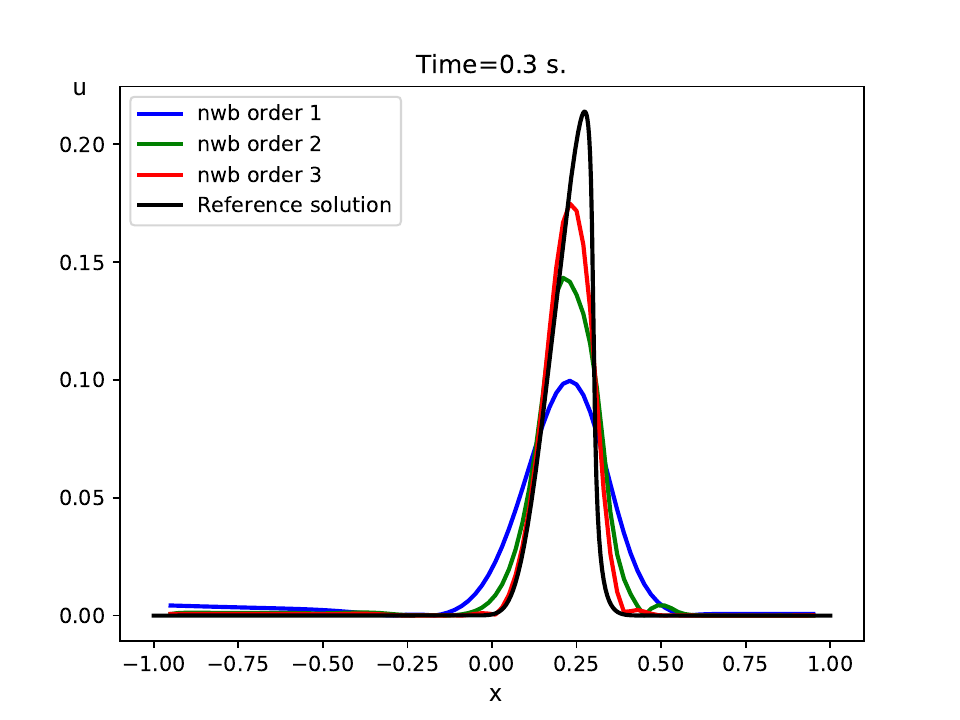}}
 \subfloat[CLWBM$i$, $i=1,2,3$. $t = 0.3$~s.]{
   \includegraphics[width=0.38\textwidth]{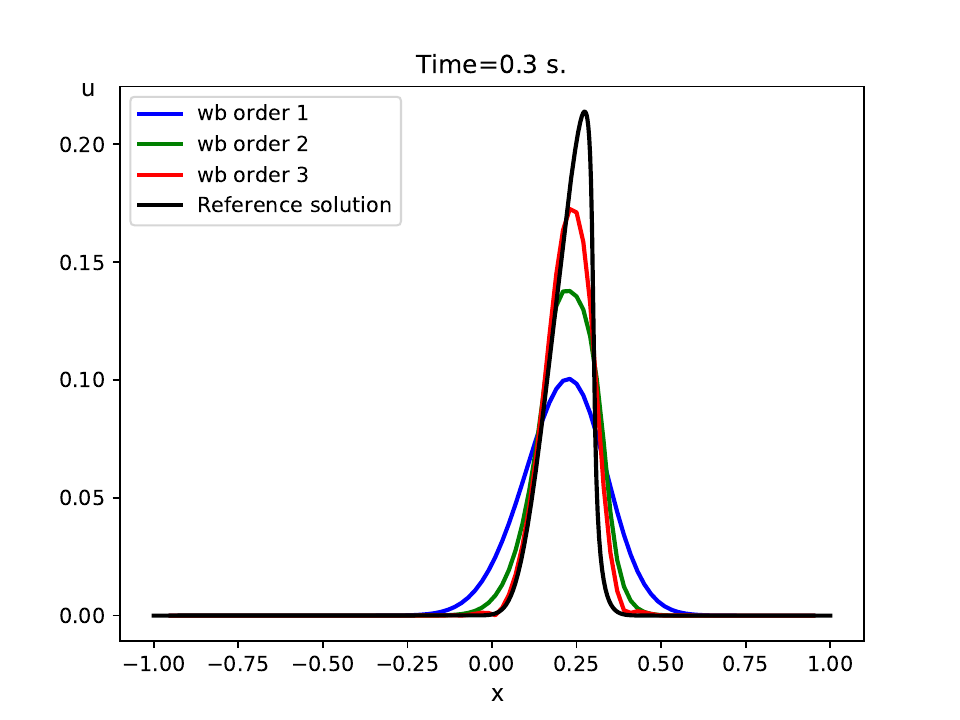}}   \vspace{-10pt}
  \subfloat[NWBM$i$, $i =1,2,3$. $t = 5$~s.]{
   \hspace{-0.3cm}\includegraphics[width=0.38\textwidth]{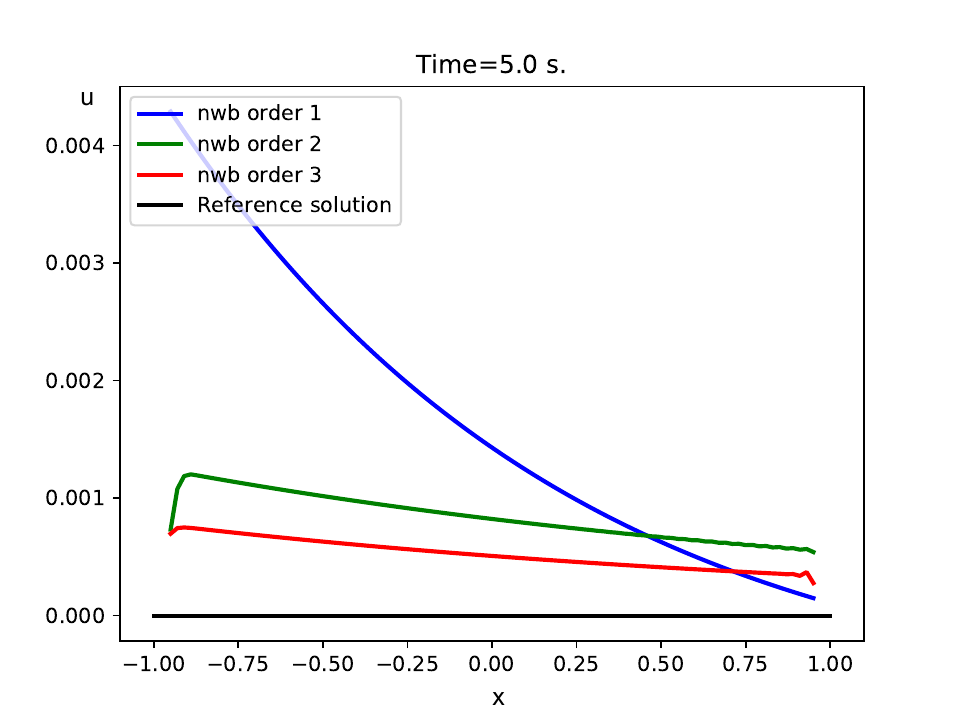}} 
   \subfloat[CLWBM$i$, $i=1,2,3$. $t = 5$~s.]{
   \includegraphics[width=0.38\textwidth]{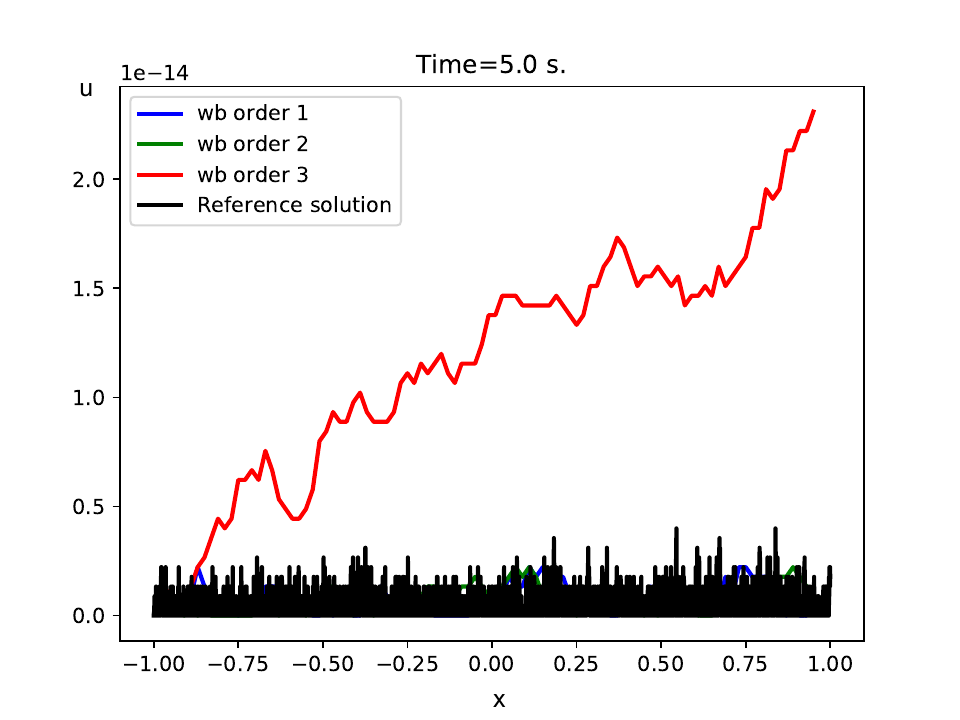}}}
    \caption{\hl{Test 2.2. Reference} and numerical solutions: differences with the stationary solution at times $t = 0.3, 5$~s for a $100$-cell~mesh.} \label{fluc_test22}

\end{figure}
Observe that, as~expected, the~well-balanced methods are able to recover the stationary solution after the perturbation left the domain, while non-well balanced methods introduce an error that remains in the domain. This could be also observed in Table~\ref{test22_errores}, where the errors with respect to the stationary solution at time $t = 5$~s are shown for the $100$-cell~mesh. 

\begin{specialtable}[H]
\caption{Test 2.2. Differences in $L^1$-norm with respect to the stationary solution for NWBM$i$, CWBM$i$, and CLWBM$i$ ($i=1,2,3$) at time $t=5$ s.} \label{test22_errores}
\setlength{\cellWidtha}{\columnwidth/4-2\tabcolsep-0in}
\setlength{\cellWidthb}{\columnwidth/4-2\tabcolsep+0in}
\setlength{\cellWidthc}{\columnwidth/4-2\tabcolsep-0.0in}
\setlength{\cellWidthd}{\columnwidth/4-2\tabcolsep-0.0in}

\scalebox{1}[1]{\begin{tabularx}{\columnwidth}{>{\PreserveBackslash\centering}m{\cellWidtha}>{\PreserveBackslash\centering}m{\cellWidthb}>{\PreserveBackslash\centering}m{\cellWidthc}>{\PreserveBackslash\centering}m{\cellWidthd}}
\toprule
\textbf{Method} &\boldmath{\textbf{Error ($i=1$)}}&\boldmath{\textbf{Error ($i = 2$)}}& \boldmath{\textbf{Error ($i = 3$)}}\\\midrule
NWBM$i$& 2.43 $\times$ $10^{-3}$ &  1.71 $\times$ $10^{-3}$  &  1.05 $\times$ $10^{-3}$\\
CWBM$i$& 1.09 $\times$ $10^{-13}$&  4.47 $\times$ $10^{-15}$ &  3.34 $\times$ $10^{-14}$ \\
CLWBM$i$& 2.52 $\times$ $10^{-15}$ &  1.19 $\times$ $10^{-15}$  &  1.24 $\times$ $10^{-14}$\\\bottomrule
\end{tabularx}}
\end{specialtable}

\subsection{Problem 3: Shallow Water~Equations}\label{Subsec:sw_transcritical}

Let us consider now the 1D shallow water system which corresponds to the choices $N=2$, with~$$ U =\begin{pmatrix}
h \\
q \\
\end{pmatrix} , \quad f(U) =\begin{pmatrix}
q \\
\displaystyle \frac{q^2}{h}+\displaystyle \frac{g}{2}h^2\\
\end{pmatrix}, \quad S(U) =\begin{pmatrix}
0\\
gh\\
\end{pmatrix},$$
where $x$ corresponds to to the axis of the channel, and $t$ is the time. The~unknown $h(x,t)$ and $q(x,t)$ are the thickness and discharge, respectively; the function $H(x)$ is the depth function, and $g$ is the gravity; $u=q/h$ is the depth-averaged velocity, and $c=\sqrt{gh}$. 

The eigenvalues of the system are $\lambda^\pm=u\pm c$, and the Froude number defined by
\begin{equation}
Fr(U)= \displaystyle \frac{|u|}{c},
\end{equation}
characterizes the flow regime: subcritical ($Fr<1$), critical ($Fr=1$), or supercritical ($Fr>1$). 

The stationary solutions $U^{*}$ satisfy the ODE system :
\begin{equation}
\begin{cases}
q_x=0,\\
\left( \displaystyle \frac{q^2}{h} + \frac{1}{2} g h^2 \right)_x= ghH_x,\\
\end{cases}
\end{equation}
which can be written as follows:
\begin{equation}\label{u'=Gsw}
\begin{cases}
\left( - {u^2}+ gh \right) h_x= ghH_x ,\\
q_x=0.
\end{cases}
\end{equation}

The stationary solutions are given in implicit form by:
\begin{equation}\label{incondimpl}
q = C_1, \quad \frac{q^2}{2 h^2} + gh -g H = C_2, \quad C_1, C_2 \in \mathbb{R}.
\end{equation}
In Reference~\cite{lopez2013}, a family of high-order well-balanced numerical methods was presented,
based in the Algorithm \ref{alg:ewbnirec}, in which local problems were solved on the basis of
this implicit form. Here, Algorithm \ref{alg:wbrec} with collocation RK methods to solve the local problems is~used.

The numerical treatment of resonant situations discussed in \hl{Section} \ref{ss:resonant} will be illustrated here in the particular case of the shallow water system. It can be easily checked that system 
\eqref{linearsys}, which, in our case, becomes
\begin{equation} \label{sw_sistema_ucritico}
\left[ 
\begin{array}{cc}
0 & 1 \\
-u^2+gh & 2u 
\end{array}
\right] 
\left[
\begin{array}{c}
K_1 \\
K_2
\end{array}
\right] 
= 
\left[
\begin{array}{c}
0 \\
g h H_x(x)
\end{array}
\right] ,
\end{equation}
when $U^*$ is a critical state, i.e., when $u^2=gh$, reduces to:
\begin{equation} \label{sw_sistema_ucritico2}
\left[ 
\begin{array}{cc}
0 & 1 \\
0 & 2u 
\end{array}
\right] 
\left[
\begin{array}{c}
K_1 \\
K_2
\end{array}
\right] 
= 
\left[
\begin{array}{c}
0 \\
g h H_x(x)
\end{array}
\right] .
\end{equation}
Therefore, the~system has solutions only if $H_x(x)=0$: in this case, the solutions are
$$
K = \alpha [1, 0]^T, \quad \alpha \in \mathbb{R}.
$$
In fact, it can be shown (see Reference~\cite{lopez2013}, for instance) that a smooth stationary solution can only reach a critical state at a minimum point $x^c$ of the depth function $H$. 
Differentiating the second relation of \eqref{incondimpl} with respect to $x$, and~using the relation $H_x(x_c)=0$, one obtains, at~$x=x_c$, 
\[
    \frac{qq_x}{h^2} + h_x\left(g-\frac{q^2}{h^3}\right) = 0,
\]
which, assuming $h_x(x_c)\neq0$, implies
\begin{equation}
    h^c(q) = \displaystyle \frac{q^{2/3}}{g^{1/3}}.
    \end{equation}

With this information in mind, let us compute the limit
$$
\lim_{x \to x^c} \frac{g h H_x}{- {u^2}+ gh},
$$
to determine the value of $h_x$ at $x^c$. This is a $0/0$ indeterminate limit, and L'H\text{\^o}pital's rule
can be applied. To~do this, first, the limit is rewritten as follows:
\begin{equation} \label{sw_critico_edo}
   \lim_{x \to x^c}      \frac{g H_x}{g-\displaystyle \frac{q^2}{h^3}} . 
    \end{equation}
Some easy computation leads to
\begin{equation}
h_x(x^c)  = \frac{h^c(q)^4}{3 h_x(x^c) q^2} g  H_{xx}(x^c);
\end{equation}
therefore:
\begin{equation}
h_x(x^c) =  \pm   \frac{h^c(q)^2}{\sqrt{3} q} g^{1/2} \sqrt{ H_{xx}(x^c)} = \pm \sqrt{ \frac{h^c(q) H_{xx}(x^c)}{3}} =  \pm \sqrt{ \frac{q^{2/3} H_{xx}(x^c)}{3 g^{1/3}}}.
\end{equation}
The above expression shows that, if $q\neq 0$ and $H_{xx}(x_c)\neq 0$, it is $h_x(x_c)\neq 0$, thus justifying the assumption used before.
Then, the~chosen solution of \eqref{linearsys} will be 
$$
K = \left[ \pm \sqrt{ \frac{q^{2/3} H_{xx}(x^c)}{3 g^{1/3}}},0 \right]^T.
$$

More precisely, when applying Algorithms \ref{alg:lp} or \ref{alg:gp}, a~threshold $\epsilon > 0$ is selected to detect critical states. If~a system \eqref{linearsys} has to be solved in which
$$
|Fr(U^*) -1 |< \epsilon,
$$
at a point $\bar x$, then
\begin{itemize}
    \item If $\bar x$ is not close to a minimum point of $H$, it is assumed that is not possible to find a smooth stationary solution that solves the problem and the algorithm is stopped. 
    \item Otherwise, the~selected solution of \eqref{linearsys} is 
 \begin{itemize}
    \item $K= \left[\displaystyle \sqrt{ \frac{q^{2/3} H_{xx}(\bar{x})}{3 g^{1/3}}}, 0\right]^T$, if $h$ is increasing close to $\bar x$, 
    \item $K=  \left[ -\displaystyle \sqrt{ \frac{q^{2/3} H_{xx}(\bar{x})}{3 g^{1/3}}}, 0 \right]^T$, if $h$ is decreasing close to $\bar x$,
\end{itemize}
and the algorithm goes~on.

\end{itemize}

\subsubsection{Test 3.1.: Transcritical Smooth Stationary~Solution}
The objective with this test is to check the correct behavior of the numerical methods of this work in presence of critical states. Following Reference~\cite{lopez2013}, we consider here a smooth transcrital stationay solution for the 1D shallow water system on the interval $[0,3]$. The~integration time is $T=1$, and the bathymetry is given by
\begin{equation} \label{sw_fondo_estacionaria_sub}
H(x)= \left\{ \begin{array}{lcc}
             -0.25(1+\cos(5 \pi (x+0.5))) &   \text{if}  & 1.3 \leq x \leq 1.7, \\
             \\ \hspace{20mm} 0 &  &\text{otherwise}. 
             
             \end{array}
   \right.
\end{equation}
The initial condition corresponds to the solution of the Cauchy problem \eqref{u'=Gsw} with initial condition 
$h(0) = 1.67750727$ and $q(0)=2.5$. It can be checked that it has a critical state at $x^c=1.5$. The~shallow water system is integrated up to time $T=1$, setting $q$ downstream, and the water height upstream. The~difference between the stationary and the numerical solutions computed with NWBM$i$, $i=1,2,3$ and CLWBM$i$, $i=1,2,3$ are shown in Figure~\ref{test41}. 
$L^1$ Errors are shown in Tables~\ref{test41_error_nwb} and \ref{test41_error_wb_col}. Notice that Reference~\cite{lopez2013} obtained similar results. The~computational costs are shown in Table~\ref{swtrans_times}.

\begin{specialtable}[H]
\caption{Test 3.1. Differences in $L^1$-norm with respect to the stationary solution and convergence rates for NWBM$i$, $i= 1,2,3$.} \label{test41_error_nwb}

\setlength{\cellWidtha}{\columnwidth/7-2\tabcolsep-0.2in}
\setlength{\cellWidthb}{\columnwidth/7-2\tabcolsep+0.2in}
\setlength{\cellWidthc}{\columnwidth/7-2\tabcolsep-0.1in}
\setlength{\cellWidthd}{\columnwidth/7-2\tabcolsep+0.2in}
\setlength{\cellWidthe}{\columnwidth/7-2\tabcolsep-0.1in}
\setlength{\cellWidthf}{\columnwidth/7-2\tabcolsep+0.2in}
\setlength{\cellWidthg}{\columnwidth/7-2\tabcolsep-0.2in}

\scalebox{1}[1]{\begin{tabularx}{\columnwidth}{>{\PreserveBackslash\centering}m{\cellWidtha}>{\PreserveBackslash\centering}m{\cellWidthb}>{\PreserveBackslash\centering}m{\cellWidthc}>{\PreserveBackslash\centering}m{\cellWidthd}>{\PreserveBackslash\centering}m{\cellWidthe}>{\PreserveBackslash\centering}m{\cellWidthf}>{\PreserveBackslash\centering}m{\cellWidthg}}
\toprule
\textbf{Cells} &\boldmath{\textbf{ Error ($i=1$)}  }&\textbf{Order}&\boldmath{\textbf{Error ($i=2$)}}  & \textbf{Order}&\boldmath{\textbf{Error ($i=3$) }}&\textbf{Order}\\
 & \multicolumn{2}{c}{\boldmath{$h$} }  & \multicolumn{2}{c}{\boldmath{$h$} } & \multicolumn{2}{c}{\boldmath{$h$}}  \\\midrule
100& 4.99 $\times$ $10^{-2}$ & - & 7.63 $\times$ $10^{-3}$  & - & 5.99 $\times$ $10^{-3}$ & -\\
200 & 1.31 $\times$ $10^{-2}$ & 1.923 & 1.27 $\times$ $10^{-3}$ & 2.583 & 8.86 $\times$ $10^{-4}$ & 2.757\\
400 & 3.87 $\times$ $10^{-3}$ & 1.766 & 1.84 $\times$ $10^{-4}$ & 2.790 & 6.20 $\times$ $10^{-5}$ & 3.838\\
800 & 1.56 $\times$ $10^{-3}$ & 1.314 & 5.31 $\times$ $10^{-5}$ & 1.794 & 6.88 $\times$ $10^{-6}$ & 3.172\\  \midrule
\textbf{Cells} & \boldmath{\textbf{Error ($i=1$)} } &\textbf{Order}&\boldmath{\textbf{Error ($i=2$) }} & \textbf{Order}&\boldmath{\textbf{Error ($i=3$)}} &\textbf{Order}\\
 & \multicolumn{2}{c}{\boldmath{$q$}}  & \multicolumn{2}{c}{\boldmath{$q$} }& \multicolumn{2}{c}{\boldmath{$q$}} \\\midrule
100& 1.28 $\times$ $10^{-1}$ & - & 1.89 $\times$ $10^{-2}$  & - & 1.72 $\times$ $10^{-2}$ & -\\
200 & 2.81 $\times$ $10^{-2}$ & 2.188 & 3.18 $\times$ $10^{-3}$ & 2.575 & 2.35 $\times$ $10^{-3}$ & 2.866\\
400 & 9.29 $\times$ $10^{-3}$ & 1.161 & 4.81 $\times$ $10^{-4}$ & 2.724 & 2.09 $\times$ $10^{-4}$ & 3.489\\
800 & 4.09 $\times$ $10^{-3}$ & 1.168 & 1.37 $\times$ $10^{-4}$ & 1.812 & 2.29 $\times$ $10^{-5}$ & 3.190\\  \bottomrule
\end{tabularx}}
\end{specialtable}
\unskip

\begin{specialtable}[H]
\caption{Test 3.1. Differences in $L^1$-norm with respect to the stationary solution for CLWBM$i$, $i=1,2,3$. } \label{test41_error_wb_col}
\begin{adjustbox}{width=\linewidth,left}  
\begin{tabular}{ccccccc} \toprule
\textbf{Cells} & \multicolumn{2}{c}{\textbf{Error (\emph{i} = 1)}} &\multicolumn{2}{c}{\textbf{Error (\emph{i} = 2)}}&\multicolumn{2}{c}{\textbf{Error (\emph{i} = 3)}}\\

  & \boldmath{$h$}&\boldmath{$q$ }&\boldmath{$h$}& \boldmath{$q$ }&\boldmath{$h$}&\boldmath{$q$}\\\midrule
100& 1.46 $\times$ $10^{-15}$ & 2.13 $\times$ $10^{-15}$ & 2.80 $\times$ $10^{-16}$  & 1.44 $\times$ $10^{-15}$ & 2.88 $\times$ $10^{-16}$ & 3.63 $\times$ $10^{-15}$\\
200 & 4.95 $\times$ $10^{-16}$ & 3.00 $\times$ $10^{-16}$ & 3.03 $\times$ $10^{-15}$ & 1.44 $\times$ $10^{-14}$ & 3.94 $\times$ $10^{-14}$ & 5.53 $\times$ $10^{-14}$\\
400 & 2.94 $\times$ $10^{-16}$ & 1.74 $\times$ $10^{-15}$ & 4.75 $\times$ $10^{-16}$ & 1.20 $\times$ $10^{-15}$ & 8.89 $\times$ $10^{-14}$ & 1.45 $\times$ $10^{-13}$\\
800 & 1.50 $\times$ $10^{-15}$ & 6.92 $\times$ $10^{-15}$ & 3.25 $\times$ $10^{-16}$ & 1.21 $\times$ $10^{-15}$ & 1.04 $\times$ $10^{-13}$ & 1.55 $\times$ $10^{-13}$\\ \bottomrule
\end{tabular}
\end{adjustbox} 
\end{specialtable}
\unskip

\begin{specialtable}[H]
\caption{Test 3.1. Computational cost (milliseconds). $t=1$~s.}\label{swtrans_times}
\setlength{\cellWidtha}{\columnwidth/4-2\tabcolsep-0in}
\setlength{\cellWidthb}{\columnwidth/4-2\tabcolsep+0in}
\setlength{\cellWidthc}{\columnwidth/4-2\tabcolsep-0.0in}
\setlength{\cellWidthd}{\columnwidth/4-2\tabcolsep-0.0in}

\scalebox{1}[1]{\begin{tabularx}{\columnwidth}{>{\PreserveBackslash\centering}m{\cellWidtha}>{\PreserveBackslash\centering}m{\cellWidthb}>{\PreserveBackslash\centering}m{\cellWidthc}>{\PreserveBackslash\centering}m{\cellWidthd}}
\toprule

          \textbf{Cells}        & \boldmath{\textbf{Order ($i$)} }& \boldmath{\textbf{NWBM$i$} }&  \boldmath{\textbf{CLWBM$i$}} \\
 \midrule
\multirow{3}{*}{100} & 1& 20 & 30  \\ 
                  & 2 & 30 & 130  \\ 
                  & 3& 70 & 320  \\ \midrule
\multirow{3}{*}{200} & 1& 30 & 60  \\ 
                  & 2& 80 & 370 \\ 
                  & 3& 200 &1100\\ \bottomrule

\end{tabularx}}
\end{specialtable}
\unskip

\vspace{-20pt}

\begin{figure}[H]
\hspace{-0.5cm}
{\captionsetup{position=bottom,justification=centering}
  \subfloat[NWBM$i$, $i=1,2,3$. $h$.]{
   \includegraphics[width=0.38\textwidth]{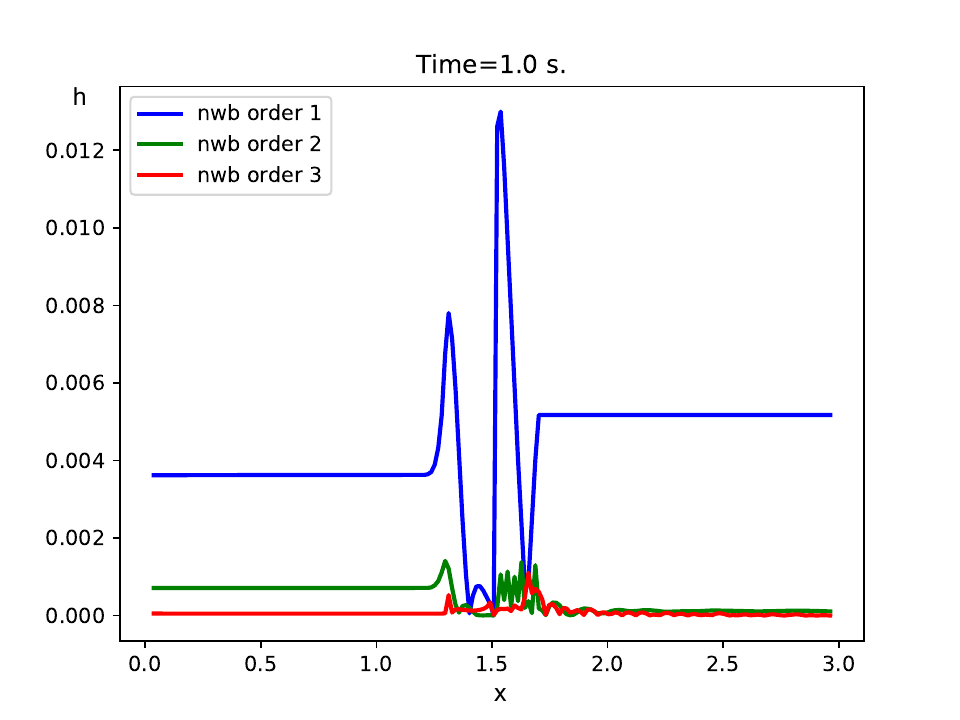}}
   \subfloat[CLWBM$i$, $i=1,2,3$. $h$.]{
   \includegraphics[width=0.38\textwidth]{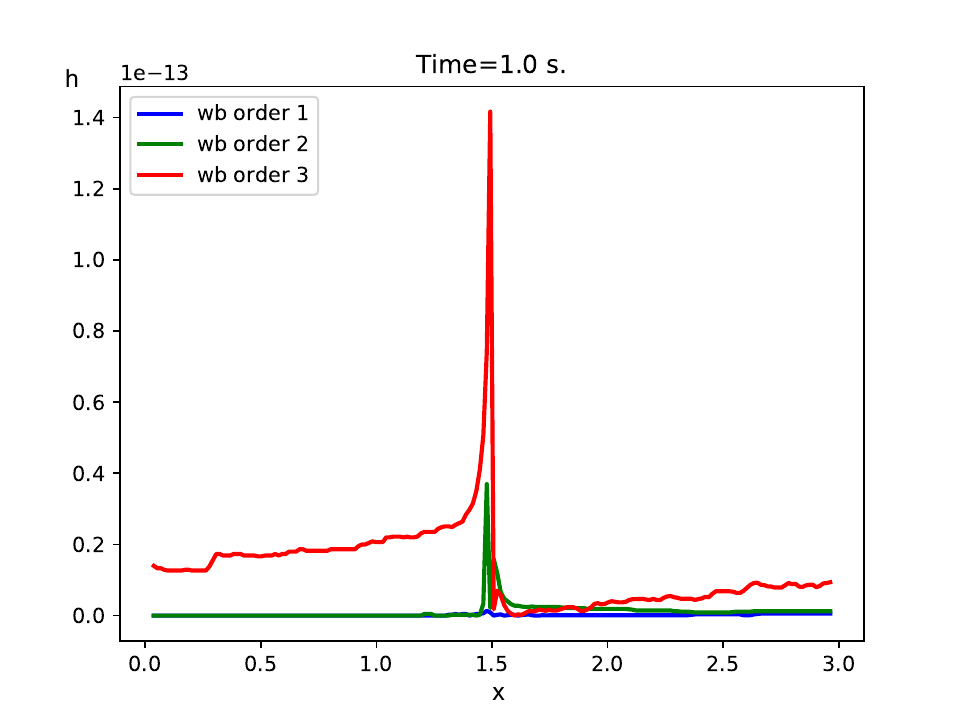}}\vspace{-13pt}
 \subfloat[NWBM$i$, $i=1,2,3$. $q$.]{
  \hspace{-0.5cm} \includegraphics[width=0.38\textwidth]{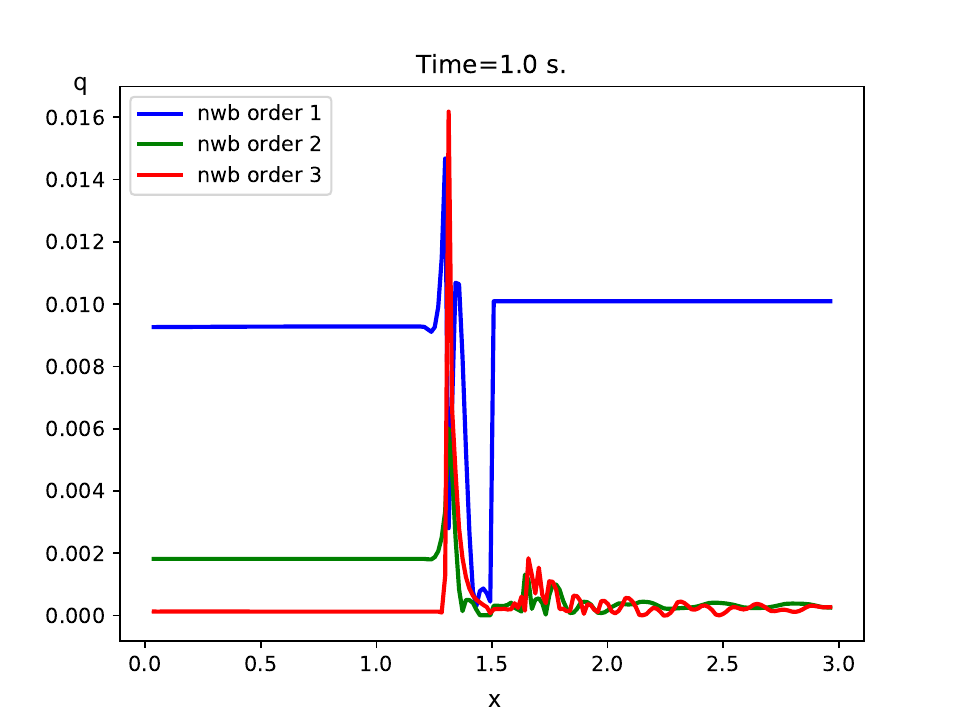}}
   \subfloat[CLWBM$i$, $i=1,2,3$. $q$.]{
   \includegraphics[width=0.38\textwidth]{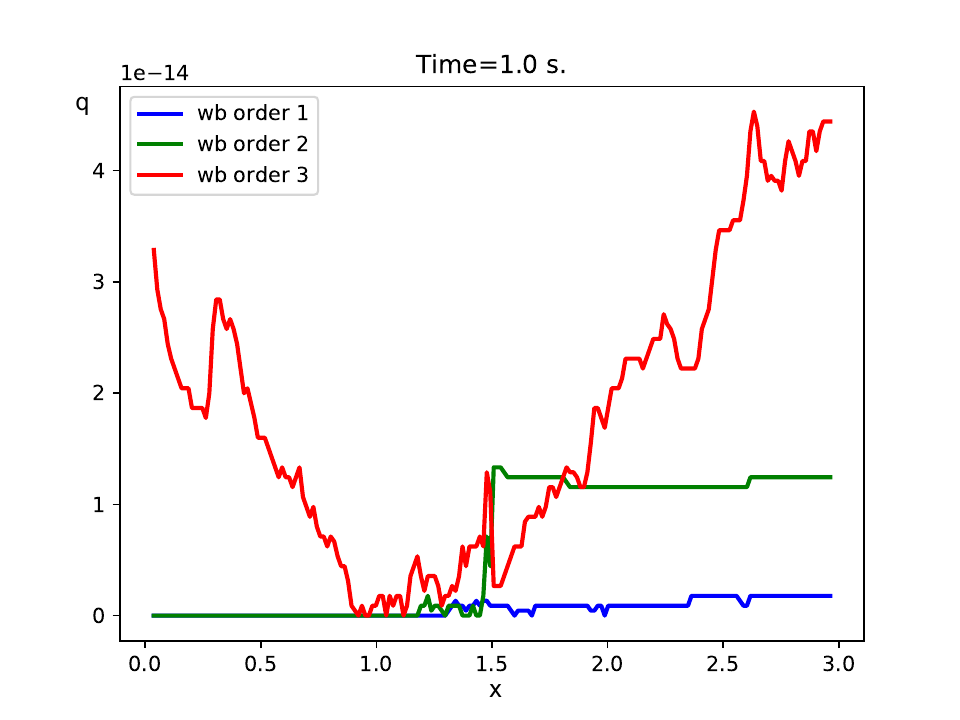}}}
    \caption{\hl{Test 3.1. Reference} and numerical solutions: differences with the stationary solution at time $t = 1$~s for the $200$-cell mesh.} \label{test41}

\end{figure}

\subsubsection{Test 3.2.: Perturbation of a Smooth Transcritical Stationary~Solution}
In this test, we consider the evolution of a small perturbation of the previous transcritical stationary solution. More precisely, a~small perturbation of the water thickness of amplitude $\Delta h= 0.02$ is set at the interval $[1.1,1.2]$, that is in the subcritical region. Therefore, this perturbation splits in two wawes, one moving downstream and the other upstream, as~shown in Figures~\ref{fluc_test42_h} and \ref{fluc_test42_q} at time
$t=0.1$~s. Figures~\ref{fluc_test42_h} and \ref{fluc_test42_q} show the difference between the stationary and the numerical solutions computed with NWBM$i$, $i=1,2,3$ and CLWBM$i$, $i=1,2,3$ at times $t=0.1, 5$~s. A~reference solution is also computed using a first order well-balanced method with fine mesh (2000 cells).
\vspace{-20pt}

\begin{figure}[H]
{\captionsetup{position=bottom,justification=centering}
  \subfloat[NWBM$i$, $i =1,2,3$. $t = 0.1$~s. ]{
   \includegraphics[width=0.38\textwidth]{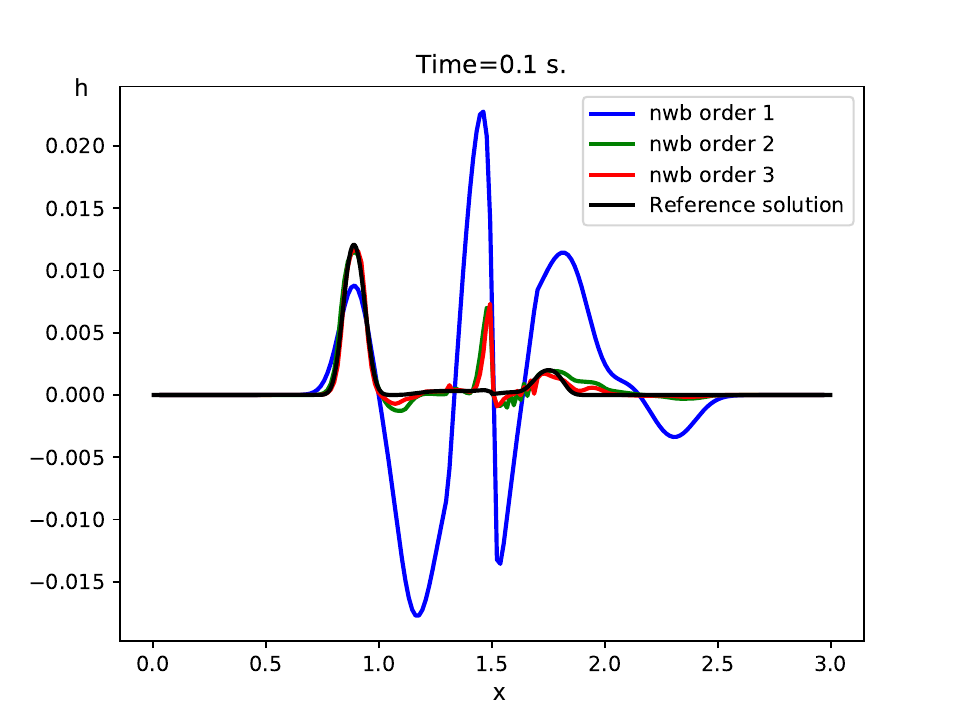}}
 \subfloat[CLWBM$i$, $i=1,2,3$. $t = 0.1$~s.]{
   \includegraphics[width=0.38\textwidth]{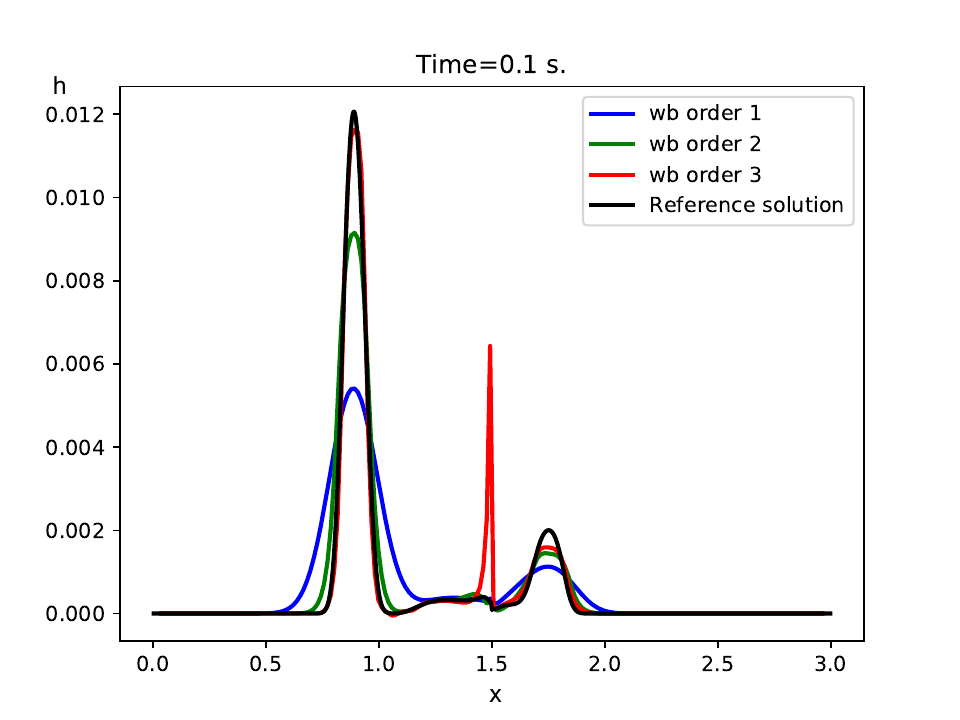}}
   \vspace{-10pt}

   \subfloat[NWBM$i$, $i =1,2,3$. $t = 5$~s.]{
   \includegraphics[width=0.38\textwidth]{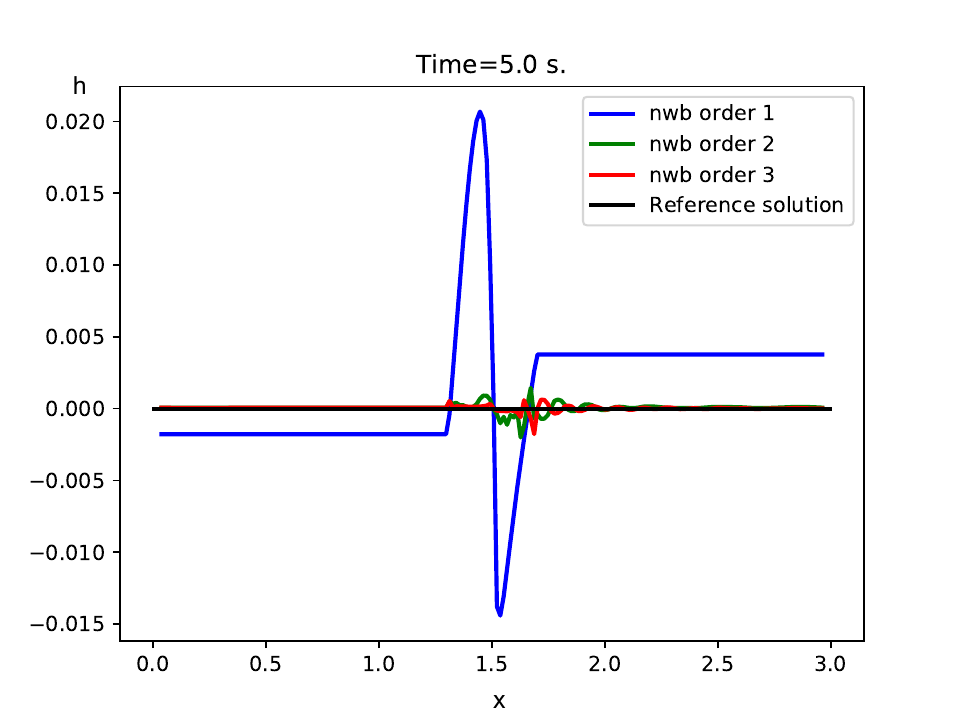}}
   \subfloat[CLWBM$i$, $i=1,2,3$. $t = 5$~s.]{
   \includegraphics[width=0.38\textwidth]{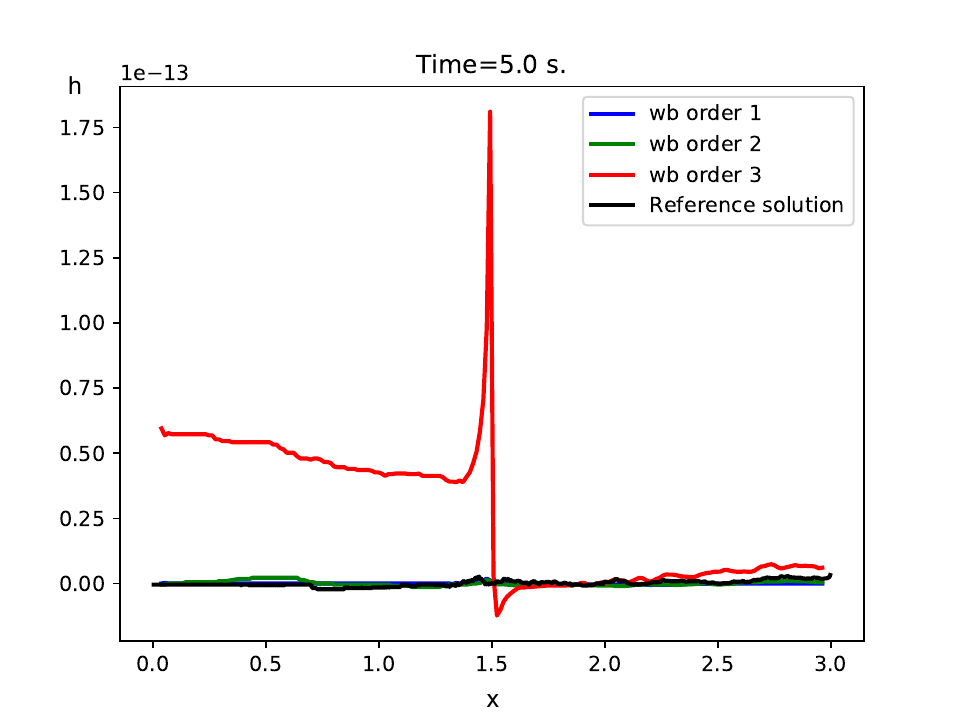}}}
    \caption{\hl{Test 3.2. Reference} and numerical solutions: differences with the stationary solution at times $t = 0.1, 5$~s for $h$ with a $200$-cell~mesh.} \label{fluc_test42_h}
\end{figure}
\vspace{-30pt}

\begin{figure}[H]\setcounter{subfigure}{0}
{\captionsetup{position=bottom,justification=centering}
 \subfloat[NWBM$i$, $i =1,2,3$. $t = 0.1$~s. ]{
   \includegraphics[width=0.38\textwidth]{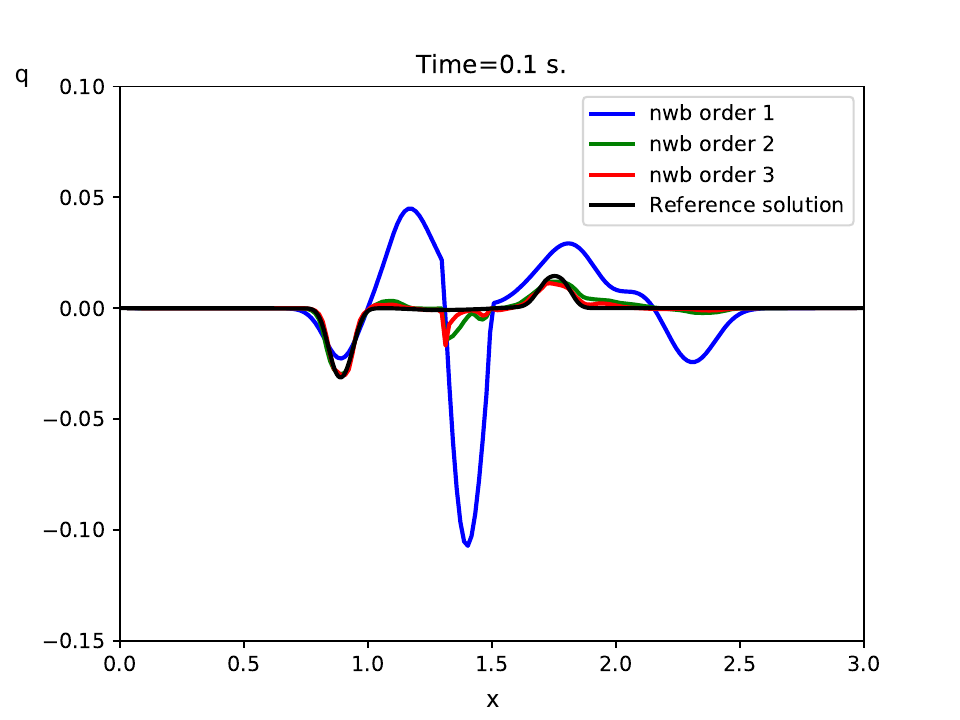}}
 \subfloat[CLWBM$i$, $i=1,2,3$. $t = 0.1$~s.]{
   \includegraphics[width=0.38\textwidth]{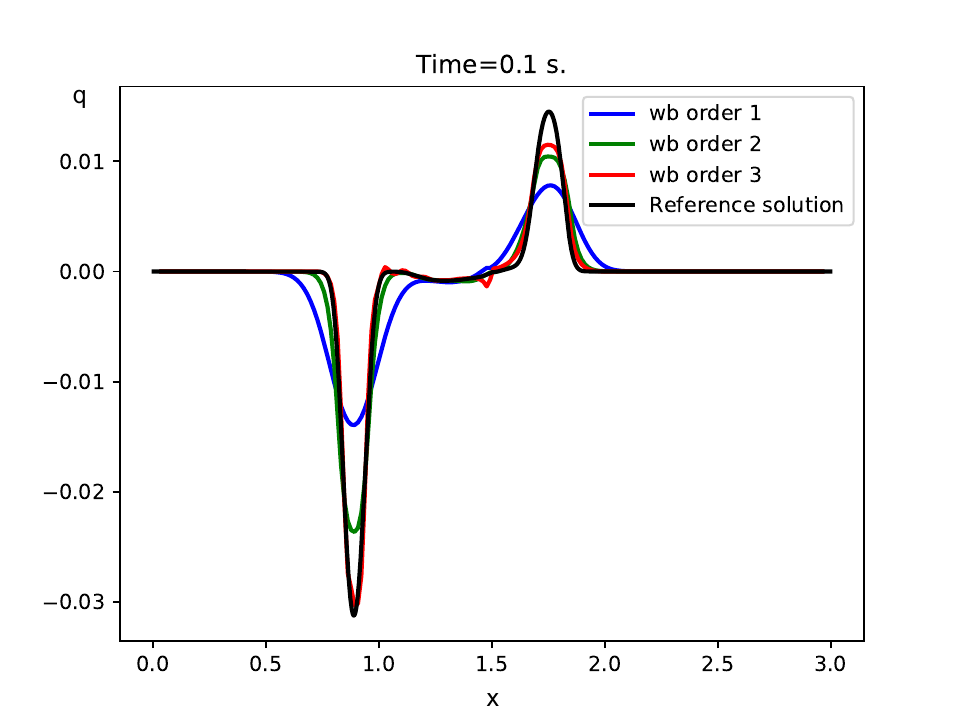}}}
   \caption{\textit{Cont}.}
\end{figure}

\begin{figure}[H]\ContinuedFloat
\centering
\setcounter{subfigure}{2}
{\captionsetup{position=bottom,justification=centering}

   \subfloat[NWBM$i$, $i =1,2,3$. $t = 5$~s.]{
 \includegraphics[width=0.38\textwidth]{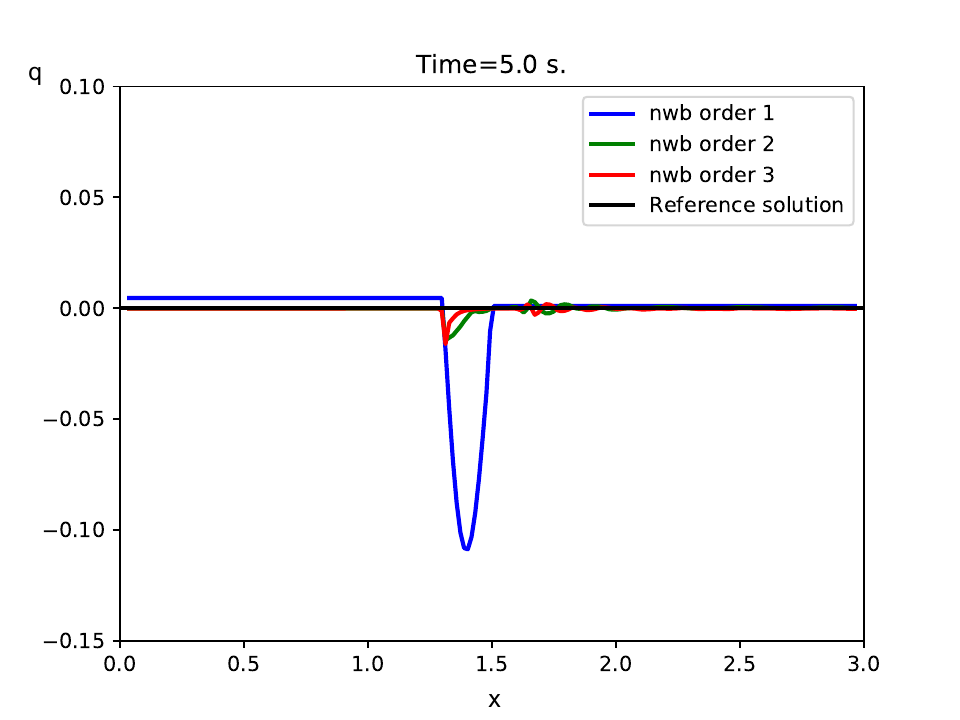}}
   \subfloat[CLWBM$i$, $i=1,2,3$. $t = 5$~s.]{
   \includegraphics[width=0.38\textwidth]{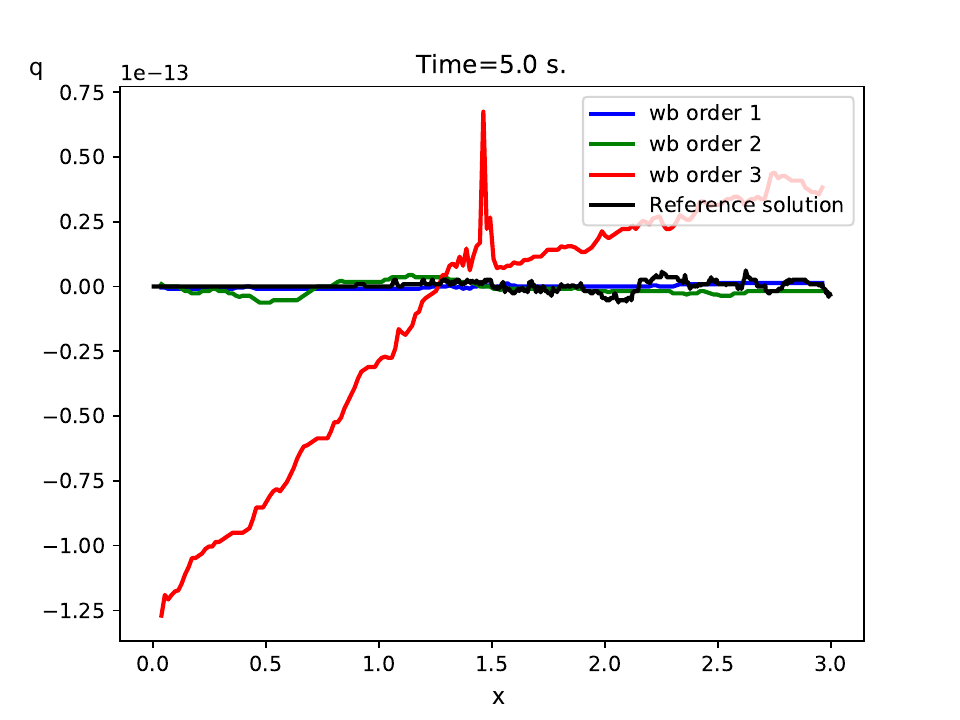}}}
    \caption{\hl{Test 3.2. Reference }and numerical solutions: differences with the stationary solution at 
    times $t = 0.1, 5$~s for $q$ for a $200$-cell mesh.} \label{fluc_test42_q}
\end{figure}
Observe that, as in the previous examples, only the well-balanced methods are able to preserve the stationary solution. The~non-well-balanced method generates strong perturbations that are bigger than the initial perturbation. Differences between the stationary and the numerical solutions in $L^1$-norm at time $t=5$~s for a $200$-cell mesh are shown in \mbox{Table~\ref{test42_errores}}. Finally, let us point out that the third order well-balanced method presents a small spurious oscillation in the location of the critical point, that disappears with time. This is not the case for the first and second order well-balanced~methods. 

\end{paracol}
\nointerlineskip
\begin{specialtable}[H]
\widetable
\caption{Test 3.2. Differences in $L^1$-norm for NWBM$i$ and CLWBM$i$ ($i=1,2,3$) with respect to the stationary solution using a $200$-cell mesh at time $t=5$~s.} \label{test42_errores}
\setlength{\cellWidtha}{\columnwidth/7-2\tabcolsep-0.2in}
\setlength{\cellWidthb}{\columnwidth/7-2\tabcolsep+0.2in}
\setlength{\cellWidthc}{\columnwidth/7-2\tabcolsep-0.1in}
\setlength{\cellWidthd}{\columnwidth/7-2\tabcolsep+0.2in}
\setlength{\cellWidthe}{\columnwidth/7-2\tabcolsep-0.1in}
\setlength{\cellWidthf}{\columnwidth/7-2\tabcolsep+0in}
\setlength{\cellWidthg}{\columnwidth/7-2\tabcolsep-0in}

\scalebox{1}[1]{\begin{tabularx}{\columnwidth}{>{\PreserveBackslash\centering}m{\cellWidtha}>{\PreserveBackslash\centering}m{\cellWidthb}>{\PreserveBackslash\centering}m{\cellWidthc}>{\PreserveBackslash\centering}m{\cellWidthd}>{\PreserveBackslash\centering}m{\cellWidthe}>{\PreserveBackslash\centering}m{\cellWidthf}>{\PreserveBackslash\centering}m{\cellWidthg}}
\toprule
\textbf{Method} & \multicolumn{2}{c}{\boldmath{\textbf{Error ($i=1$)}}} &\multicolumn{2}{c}{\boldmath{\textbf{Error ($i=2$)}}}&\multicolumn{2}{c}{\boldmath{\textbf{Error ($i=3$)}}}\\

  & \boldmath{$h$}&\boldmath{$q$ }&\boldmath{$h$}&\boldmath{ $q$} &\boldmath{$h$}&\boldmath{$q$}\\\midrule
NWBM$i$& 1.01 $\times$ $10^{-2}$  & 2.13 $\times$ $10^{-2}$ & 4.83 $\times$ $10^{-4}$ & 2.00 $\times$ $10^{-3}$ & 2.70 $\times$ $10^{-4}$ & 1.30 $\times$ $10^{-3}$\\
CLWBM$i$& 5.11 $\times$ $10^{-14}$ & 4.32 $\times$ $10^{-16}$ & 1.95 $\times$ $10^{-14}$ & 6.47 $\times$ $10^{-15}$ & 7.89 $\times$ $10^{-14}$ & 1.38 $\times$ $10^{-13}$\\\bottomrule

\end{tabularx}} 
\end{specialtable}
\begin{paracol}{2}
\switchcolumn

\vspace{-10pt}

\subsubsection{Test 3.3.: Perturbation of a Smooth Subcritical Stationary~Solution}
Let us consider evolution of a small perturbation of a smooth subcritical stationary solution. The~initial condition $U_0(x)=(h_0(x),q_0(x))^T$ is given by: 
$$
h_0(x)=\begin{cases}
h^*(x)+0.02,
& \mbox{if $0.7 \leq x \leq 1.0$,}\\
h^*(x), & \mbox{otherwise,}
\end{cases}
$$
$$
q_0(x)=q^*(x),$$
where $U^*(x)=(h^*(x),q^*(x))^T$ is the solution of the Cauchy problem
\begin{equation} \label{sw_estacionaria_sub}
\begin{cases}
q_x=0,\\
\left( - {u^2}+ gh \right) h_x= ghH_x,\\
h(0)=2,\, q(0)=3.5.
\end{cases}
\end{equation}
The depth function is given again by \eqref{sw_fondo_estacionaria_sub}.

Figure~\ref{fluc_test33} shows the differences between the stationary solution and the numerical solutions obtained with NWBM$i$, $i=1,2,3$ and CLWBM$i$, $i=1,2,3$ at times $t=0.1, 5$ for $h$ (the graphs are similar for $q$). A~reference solution has been computed with a first order well-balanced scheme on a fine mesh (3200 cells). 

Notice that, since the perturbation amplitude is small enough, the~well-balanced methods perform much better than the non-well-balanced methods, that generates spurious waves of amplitudes bigger than the original perturbation. In~Table~\ref{test33_errores}, the differences with respect to the stationary solution at time $t = 5$~s are shown for the $100$-cell mesh. {In order to check the efficiency of the methods, Figure~\ref{test33_eficiencia} shows the errors in $L^1$-norm with respect to the reference solution versus the CPU times in milliseconds for the NWBM$2$ and CLWBM$2$ at time $t=1$~s, when the perturbation is still in the domain. Similar results have been obtained for the first and third order methods. Notice that the CLWBM are computationally more efficient. See, for~instance, that, in order to obtain an error of approximately 2 $\times$ $10^{-4}$, the~computational cost for the NWBM$2$ is increased by a factor of about 2.2 with respect to the CLWBM$2$.}

\end{paracol}
\nointerlineskip
\begin{specialtable}[H]
\widetable
\caption{Test 3.3. Differences in $L^1$-norm for NWBM$i$ and CLWBM$i$ ($i=1,2,3$) with respect to the stationary solution for the $100$-cell mesh at time $t=5$~s.} \label{test33_errores}
\setlength{\cellWidtha}{\columnwidth/7-2\tabcolsep-0.2in}
\setlength{\cellWidthb}{\columnwidth/7-2\tabcolsep+0.2in}
\setlength{\cellWidthc}{\columnwidth/7-2\tabcolsep-0.1in}
\setlength{\cellWidthd}{\columnwidth/7-2\tabcolsep+0.2in}
\setlength{\cellWidthe}{\columnwidth/7-2\tabcolsep-0.1in}
\setlength{\cellWidthf}{\columnwidth/7-2\tabcolsep+0in}
\setlength{\cellWidthg}{\columnwidth/7-2\tabcolsep-0in}

\scalebox{1}[1]{\begin{tabularx}{\columnwidth}{>{\PreserveBackslash\centering}m{\cellWidtha}>{\PreserveBackslash\centering}m{\cellWidthb}>{\PreserveBackslash\centering}m{\cellWidthc}>{\PreserveBackslash\centering}m{\cellWidthd}>{\PreserveBackslash\centering}m{\cellWidthe}>{\PreserveBackslash\centering}m{\cellWidthf}>{\PreserveBackslash\centering}m{\cellWidthg}}
\toprule
\textbf{Method} & \multicolumn{2}{c}{\boldmath{\textbf{Error ($i=1$)}}} &\multicolumn{2}{c}{\boldmath{\textbf{Error ($i=2$)}}}&\multicolumn{2}{c}{\boldmath{\textbf{Error ($i=3$)}}}\\

  & \boldmath{$h$}&\boldmath{$q$ }&\boldmath{$h$}& \boldmath{$q$} &\boldmath{$h$}&\boldmath{$q$}\\\midrule
NWBM$i$& 5.08 $\times$ $10^{-2}$  & 1.94 $\times$ $10^{-1}$ & 9.33 $\times$ $10^{-3}$ & 3.51 $\times$ $10^{-2}$ & 6.02 $\times$ $10^{-3}$ & 2.12 $\times$ $10^{-2}$\\
CLWBM$i$& 1.81 $\times$ $10^{-15}$ & 5.54 $\times$ $10^{-15}$ & 1.95 $\times$ $10^{-15}$ & 4.45 $\times$ $10^{-15}$ & 2.46 $\times$ $10^{-14}$ & 5.20 $\times$ $10^{-14}$\\\bottomrule

\end{tabularx}}
\end{specialtable}

\begin{paracol}{2}
\switchcolumn

\vspace{-30pt}
\begin{figure}[H]
\setcounter{subfigure}{0}
\hspace{-0.3cm}
{\captionsetup{position=bottom,justification=centering}
  \subfloat[NWBM$i$, $i =1,2,3$. $t = 0.1$~s.]{
   \includegraphics[width=0.38\textwidth]{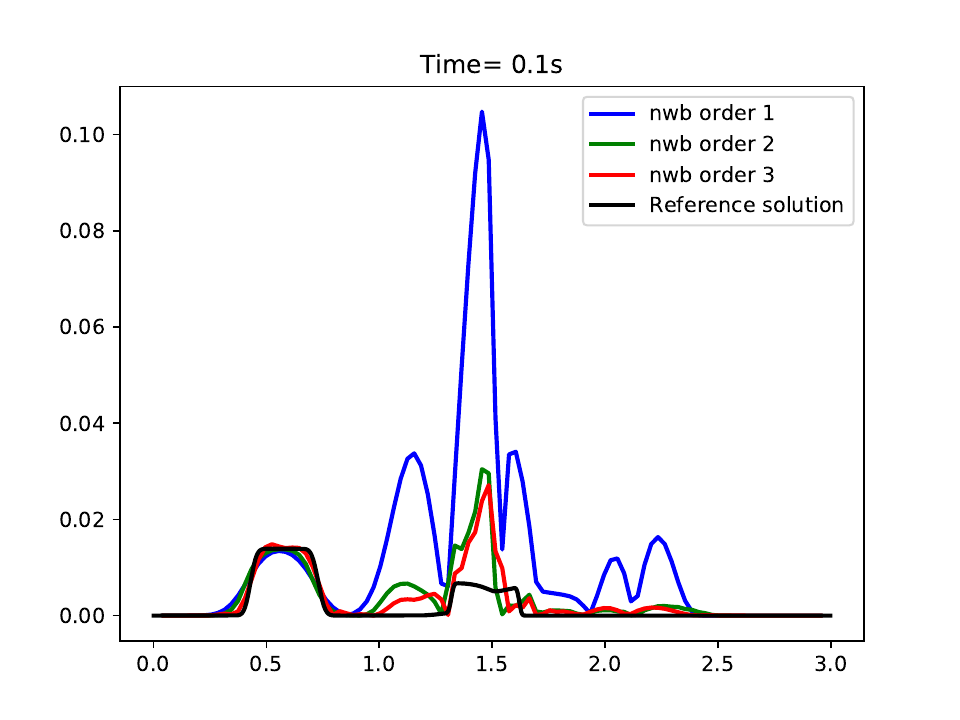}}
   \subfloat[CLWBM$i$, $i =1,2,3$. $t = 0.1$~s. ]{
   \includegraphics[width=0.38\textwidth]{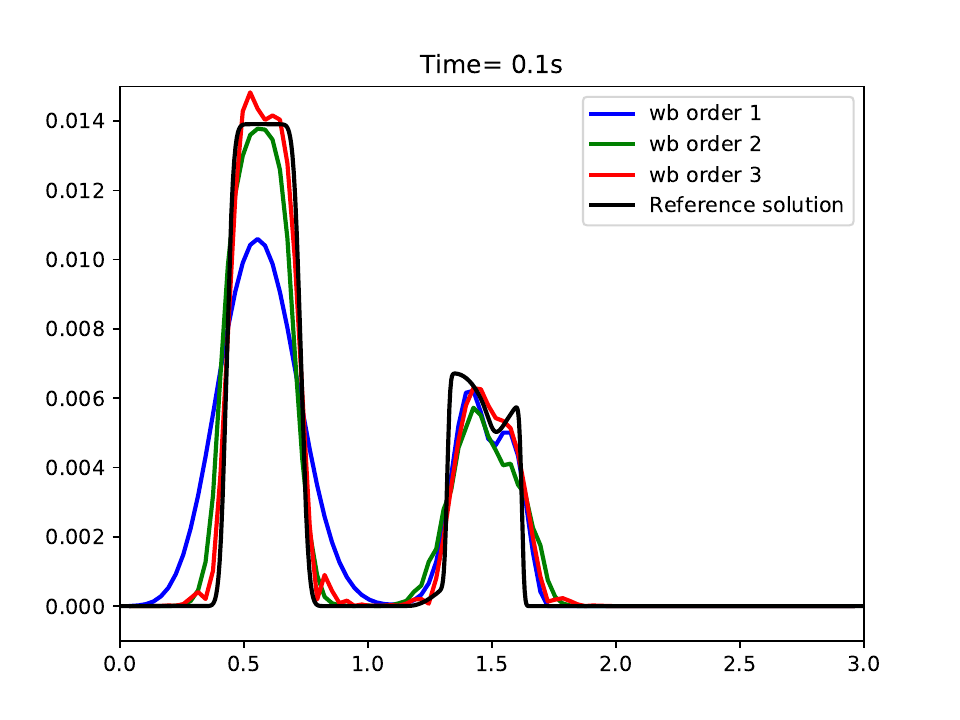}}\vspace{0.0001mm}
   \vspace{-10pt}

   \subfloat[NWBM$i$, $i =1,2,3$. $t = 5$~s.]{
   \hspace{-0.3cm}\includegraphics[width=0.38\textwidth]{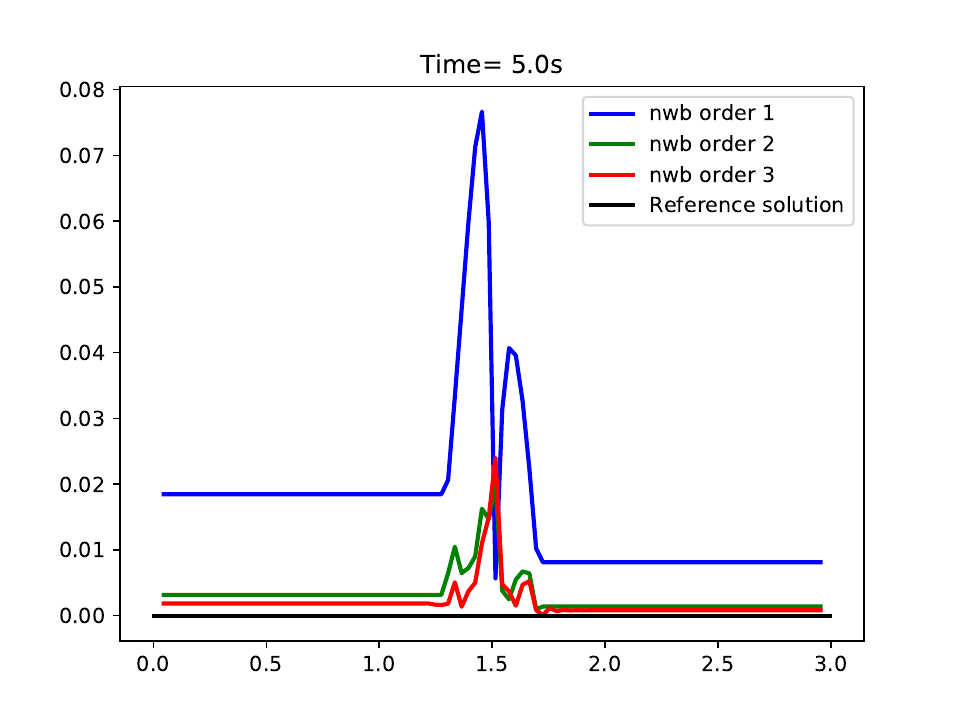}}
   \subfloat[CLWBM$i$, $i=1,2,3$. $t = 5$~s.]{
   \includegraphics[width=0.38\textwidth]{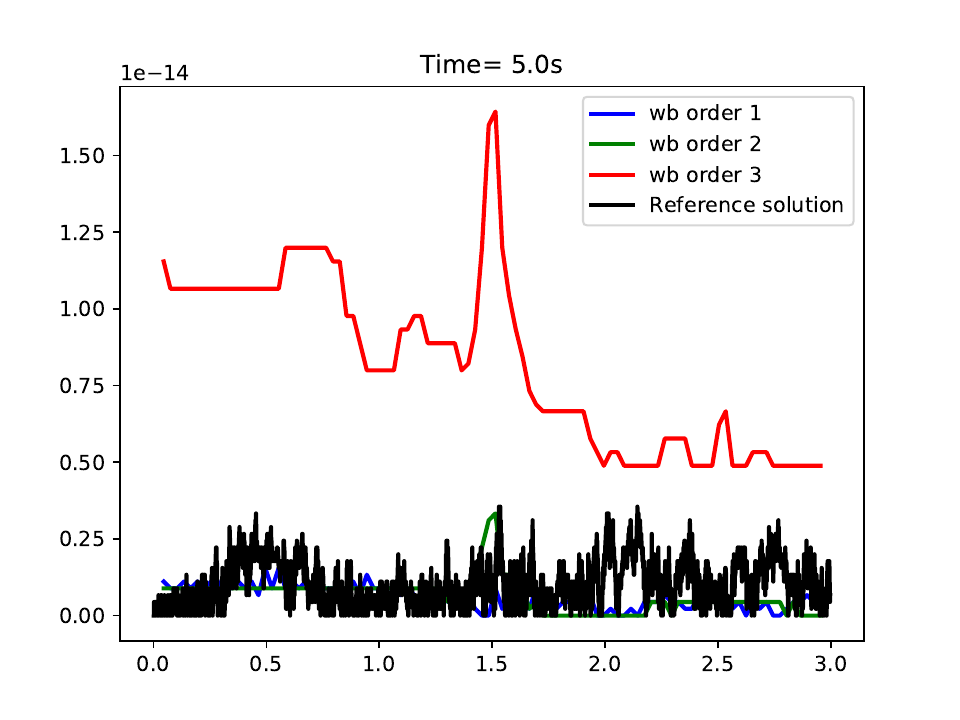}}}
    \caption{\hl{Test 3.3. Reference} and numerical solutions: differences with the stationary solution at times $t = 0.1$~s and $t = 5$~s for $h$ for a $100$-cell~mesh.} \label{fluc_test33}
\end{figure}
\vspace{-20pt}

\begin{figure}[H]  
   \includegraphics[width=0.5\textwidth]{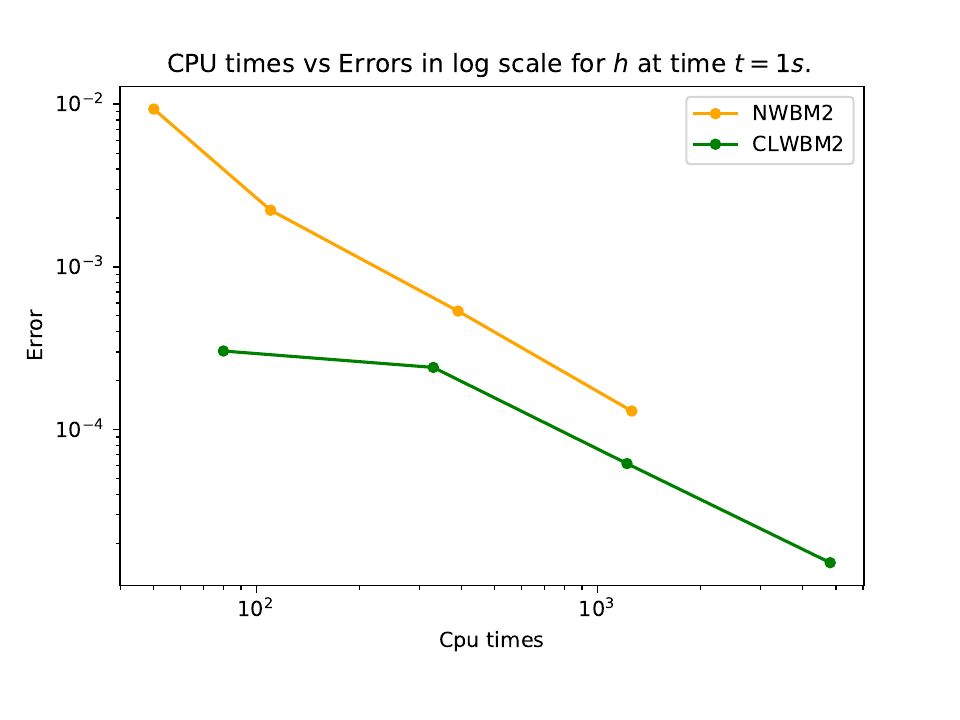}   
   \vspace{-18pt}
    \caption{{Test 3.3. Errors in $L^1$-norm with respect to the reference solution versus CPU times in milliseconds for the NWBM$2$ and CLWBM$2$ at time $t=1$~s.}} \label{test33_eficiencia}
\end{figure}

\subsection{Problem 4: Shallow Water Equations with Manning~Friction}
Let us consider the shallow water equations with Manning friction:
\begin{equation}\label{swf_ecuacion}
    \begin{cases}
    h_t+q_x=0,\\
q_t + \left( \displaystyle \frac{q^2}{h} + \frac{1}{2} g h^2 \right)_x= ghH_x - \displaystyle \frac{k q |q|}{h^{\eta}}.\\
    
    \end{cases}
\end{equation}
In the previous, $k$ is the Manning friction coefficient, and $\eta$ is set to $ \frac{7}{3}$.

The stationary solutions satisfy the ODE system
\begin{equation}\label{u'=Gswf}
\begin{cases}
\left( - {u^2}+ gh \right) h_x= ghH_x- \displaystyle \frac{k q |q|}{h^{\eta}} ,\\
q_x=0.
\end{cases}
\end{equation}

In what follows, we consider some experiments taken from Reference~\cite{michel2017well}.

\subsubsection{Test~4.1}
Let us consider first a moving equilibria test with constant bathymetry. Therefore, the~smooth stationary solutions satisfy \begin{equation}\label{u'=Gswf_fondoplano}
\begin{cases}
\left( - q^2/h^2+ gh \right) h_x=- \displaystyle \frac{k q |q|}{h^{\eta}} ,\\
q_x=0.
\end{cases}
\end{equation}
As in Reference~\cite{michel2017well}, we consider the interval the space interval $[0,1]$, discretized with $200$ cells. Now, we consider a subcritical steady state with $q(0) = -1$ and $h(0)$ a positive root of the non-linear function $\xi(h)$, obtained by integrating \eqref{u'=Gswf_fondoplano}, given by
\begin{equation}
    \xi(h)=- \displaystyle \frac{q_0^2}{\eta -1} \left( h^{\eta -1}-h_0^{\eta -1}\right) + \displaystyle \frac{g}{\eta + 2} \left( h^{\eta +2}-h_0^{\eta +2}\right) + kq_0 |q_0| (x-x_0),
\end{equation}
where $x_0=- \Delta x$, $x=0$, $q_0=q(0)$ and $h_0=h^c$ , with~$h^c$ defined by
$$
h^c=\left( \displaystyle \frac{q_0^2}{g} \right) ^{\frac{1}{3}}.
$$
In this experiment, $k=1$, which enforces the stiff character of the friction term. As~boundary conditions, $q$ is set downstream and $h$ upstream, and~the system is integrated up to $t=1$. Due to the stiff character of the friction term, the~explicit non-well-balanced methods turn out to be unstable for small values of $h$ and/or big values of $k$, unless~a very small time step is used. If~an implicit discretization of the stiff source term is considered, NWBM$i$, $i=1,2,3$ methods do not explode, yet there is not convergence to any stationary solution. The~stability restriction is much less severe for the well-balanced schemes. For~these reasons, only the numerical solutions obtained with well-balanced schemes are~shown.

Figure~\ref{test51_est_dif} shows the discrepancies between the stationary and the numerical solutions at $t = 1$~s with CLWBM$i$, $i=1,2,3$;
Table~\ref{test51_error_est} shows the $L^1$ errors for a mesh with $200$ cells.


\vspace{-20pt}

\begin{figure}[H]
{\captionsetup{position=bottom,justification=centering}
     \subfloat[CLWBM$i$, $i=1,2,3$. $h$.]{
   \includegraphics[width=0.38\textwidth]{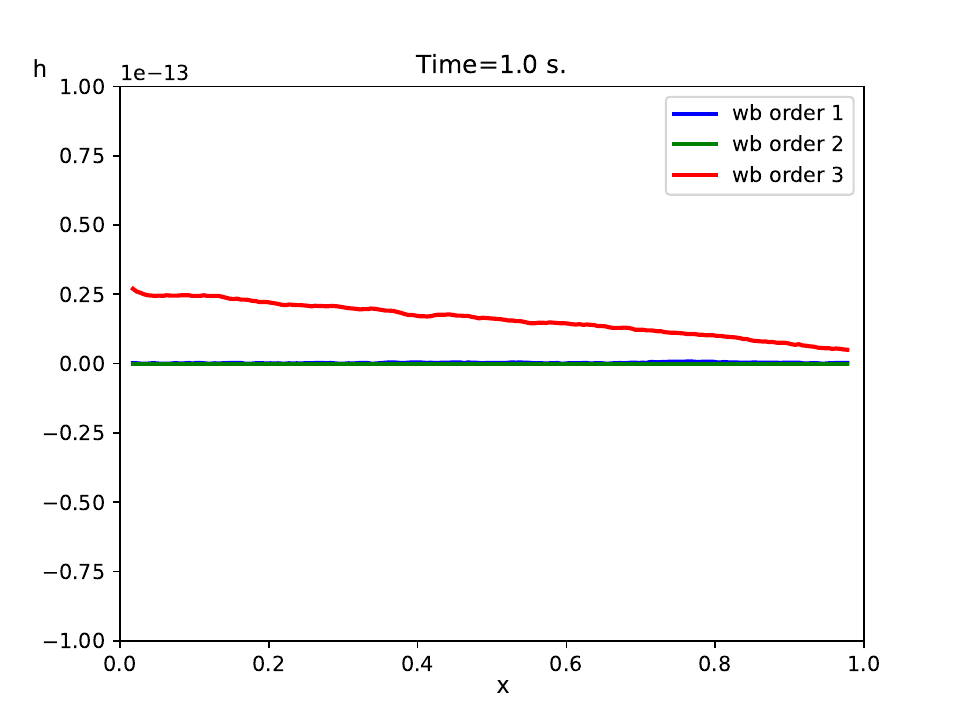}}
   \subfloat[CLWBM$i$, $i=1,2,3$. $q$.]{
   \includegraphics[width=0.38\textwidth]{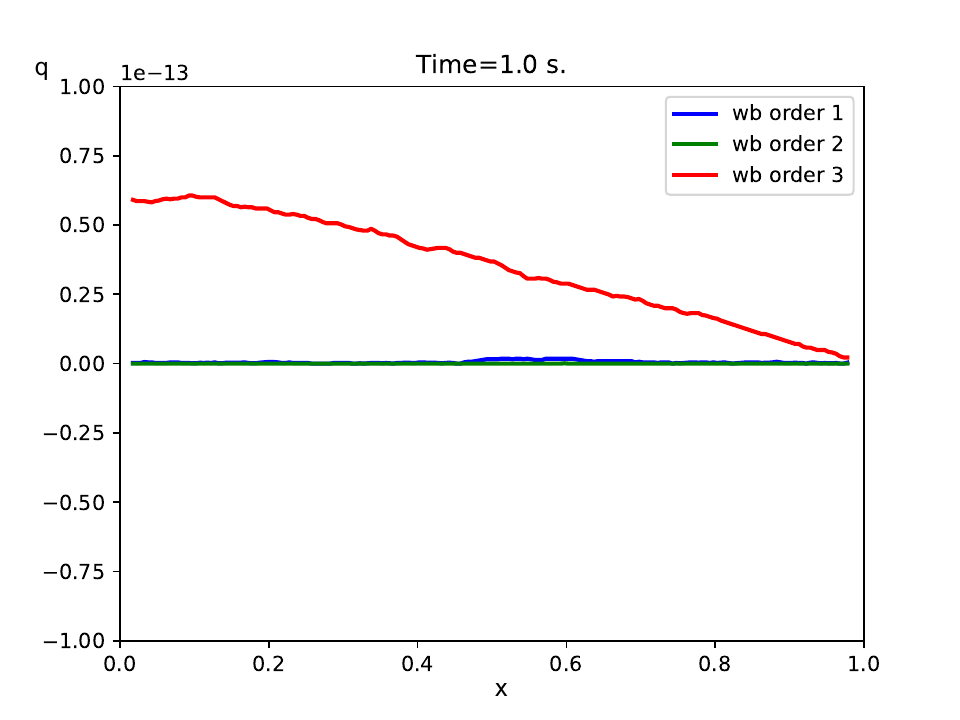}}}
    \caption{\hl{Test 4.1.} Differences between the stationary and the numerical solutions for CLWBM$i$, $i=1,2,3$ at time $t = 1$~s for a $200$-cell~mesh.} \label{test51_est_dif}
\end{figure}
\unskip

\end{paracol}
\nointerlineskip
\begin{specialtable}[H]
\widetable
\caption{Test 4.1. Differences in $L^1$-norm with respect to the stationary solution for CLWBM$i$ ($i=1,2,3$) for the $200$-cell mesh at time $t=1$~s.} \label{test51_error_est}
\setlength{\cellWidtha}{\columnwidth/7-2\tabcolsep-0in}
\setlength{\cellWidthb}{\columnwidth/7-2\tabcolsep+0in}
\setlength{\cellWidthc}{\columnwidth/7-2\tabcolsep-0in}
\setlength{\cellWidthd}{\columnwidth/7-2\tabcolsep+0in}
\setlength{\cellWidthe}{\columnwidth/7-2\tabcolsep-0in}
\setlength{\cellWidthf}{\columnwidth/7-2\tabcolsep+0in}
\setlength{\cellWidthg}{\columnwidth/7-2\tabcolsep-0in}

\scalebox{1}[1]{\begin{tabularx}{\columnwidth}{>{\PreserveBackslash\centering}m{\cellWidtha}>{\PreserveBackslash\centering}m{\cellWidthb}>{\PreserveBackslash\centering}m{\cellWidthc}>{\PreserveBackslash\centering}m{\cellWidthd}>{\PreserveBackslash\centering}m{\cellWidthe}>{\PreserveBackslash\centering}m{\cellWidthf}>{\PreserveBackslash\centering}m{\cellWidthg}}
\toprule
\textbf{Method} & \multicolumn{2}{c}{\boldmath{\textbf{Error ($i=1$)}} }&\multicolumn{2}{c}{\boldmath{\textbf{Error ($i=2$)}}}&\multicolumn{2}{c}{\boldmath{\textbf{Error ($i=3$)}}}\\

  & \boldmath{$h$}&\boldmath{$q$ }&\boldmath{$h$}&\boldmath{ $q$} &\boldmath{$h$}&\boldmath{$q$}\\\midrule
CLWBM$i$& 2.24 $\times$ $10^{-16}$ & 5.06 $\times$ $10^{-16}$ & 0.00 & 5.56 $\times$ $10^{-19}$ & 1.58 $\times$ $10^{-14}$ & 3.38 $\times$ $10^{-14}$\\\bottomrule

\end{tabularx}}

\end{specialtable}
\begin{paracol}{2}
\switchcolumn

\vspace{-8pt}

\subsubsection{Test~4.2}
Now, we consider a numerical experiment that corresponds to a small perturbation of the previous stationary solution $U^*(x) = [h^*(x), q^*(x)]^T$: 
$$
h_0(x)=\begin{cases}
h^*(x)+0.2,
& \mbox{if $\displaystyle \frac{3}{7} \leq x \leq \displaystyle \frac{4}{7}$,}\\
h^*(x), & \mbox{otherwise,}
\end{cases}
$$
$$
q_0(x)=q^*(x).$$
As in Reference~\cite{michel2017well}, we consider a mesh with $100$ cells, and the model is run up to $t= 9$~s.

As in the previous test case, only the well-balanced schemes will be considered. 
\mbox{Figure~\ref{test51_per_dif}} shows the discrepancies between the stationary and the numerical solutions for the water depth computed with CLWBM$i$, $i=1,2,3$ at the time steps $t=0.06$~s and $t=9$~s (the figures are similar for $q$). In~black, we plot a reference solution computed with a 1st-order well-balanced scheme with a $1600$-cell~mesh.
\vspace{-20pt}

\begin{figure}[H]
  {\captionsetup{position=bottom,justification=centering}
 \subfloat[CLWBM$i$, $i=1,2,3$. $t = 0.06$~s.]{
   \includegraphics[width=0.38\textwidth]{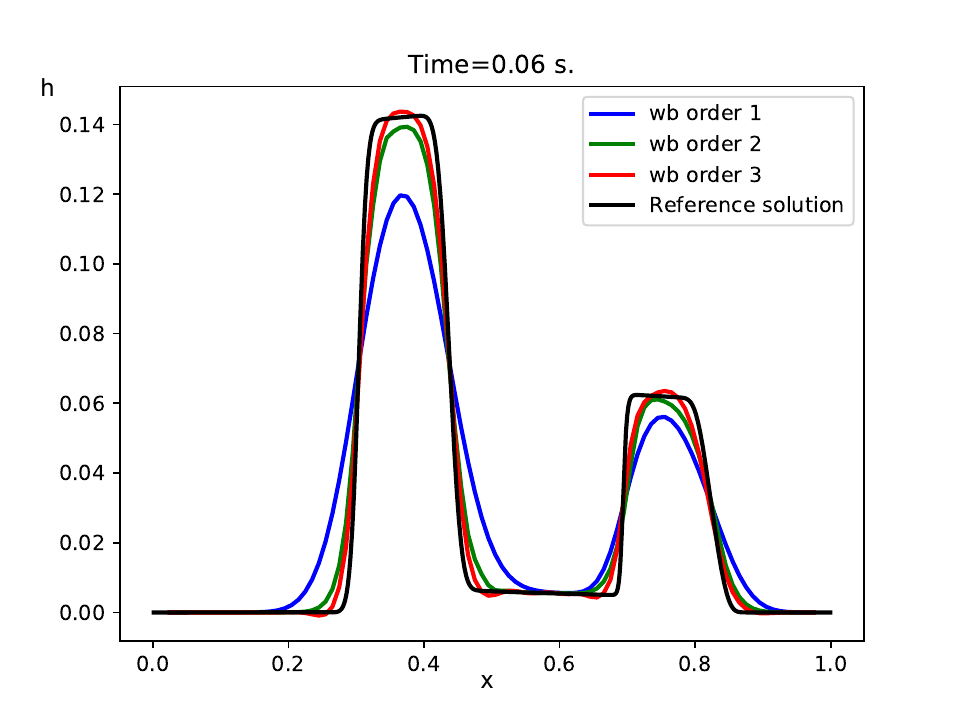}}
   \subfloat[CLWBM$i$, $i=1,2,3$. $t = 9$~s.]{
   \includegraphics[width=0.38\textwidth]{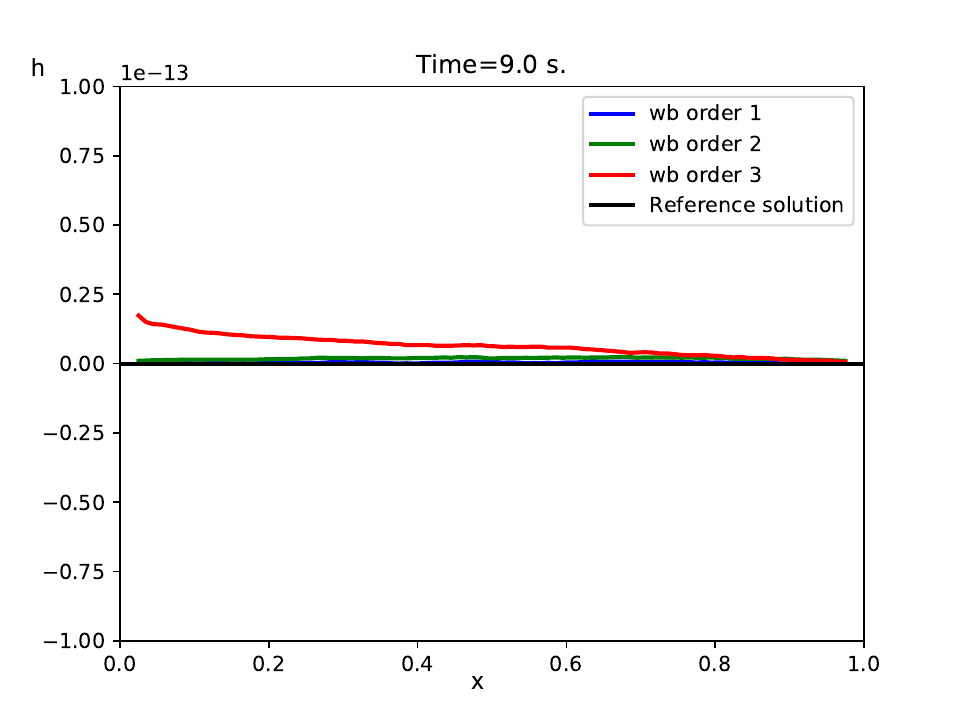}}}
    \caption{\hl{Test 4.2. Reference} and numerical solutions: differences with the stationary solution for CLWBM$i$, $i=1,2,3$, at~times $t = 0.06, 9$~s for $h$ on a $100$-cell~mesh.} \label{test51_per_dif}
\end{figure}

As expected, the~stationary solution is preserved by the well-balanced schemes after the perturbation left the domain. Table~\ref{test51_error_per} shows the $L^1$ errors with respect to the steady state at $t = 9$~s for $100$-cell~mesh.

\end{paracol}
\nointerlineskip
\begin{specialtable}[H]
\widetable
\caption{Test 4.2. Differences in $L^1$-norm with respect to the stationary solution for CLWBM$i$ ($i=1,2,3$) for the $100$-cell mesh at time $t=9$~s.} \label{test51_error_per}
\setlength{\cellWidtha}{\columnwidth/7-2\tabcolsep-0in}
\setlength{\cellWidthb}{\columnwidth/7-2\tabcolsep+0in}
\setlength{\cellWidthc}{\columnwidth/7-2\tabcolsep-0in}
\setlength{\cellWidthd}{\columnwidth/7-2\tabcolsep+0in}
\setlength{\cellWidthe}{\columnwidth/7-2\tabcolsep-0in}
\setlength{\cellWidthf}{\columnwidth/7-2\tabcolsep+0in}
\setlength{\cellWidthg}{\columnwidth/7-2\tabcolsep-0in}

\scalebox{1}[1]{\begin{tabularx}{\columnwidth}{>{\PreserveBackslash\centering}m{\cellWidtha}>{\PreserveBackslash\centering}m{\cellWidthb}>{\PreserveBackslash\centering}m{\cellWidthc}>{\PreserveBackslash\centering}m{\cellWidthd}>{\PreserveBackslash\centering}m{\cellWidthe}>{\PreserveBackslash\centering}m{\cellWidthf}>{\PreserveBackslash\centering}m{\cellWidthg}}
\toprule
\textbf{Method} & \multicolumn{2}{c}{\boldmath{\textbf{Error ($i=1$)}} }&\multicolumn{2}{c}{\boldmath{\textbf{Error ($i=2$)}}}&\multicolumn{2}{c}{\boldmath{\textbf{Error ($i=3$)}}}\\
  & \boldmath{$h$}&\boldmath{$q$ }&\boldmath{$h$}&\boldmath{ $q$} &\boldmath{$h$}&\boldmath{$q$}\\\midrule

CLWBM$i$& 2.99 $\times$ $10^{-16}$ & 3.97 $\times$ $10^{-16}$ & 1.81 $\times$ $10^{-15}$ & 2.76 $\times$ $10^{-15}$ & 6.50 $\times$ $10^{-14}$ & 1.77 $\times$ $10^{-14}$\\\bottomrule

\end{tabularx}}

\end{specialtable}
\begin{paracol}{2}
\switchcolumn

\vspace{-8pt}

\subsubsection{Test~4.3}
We end this subsection with two numerical experiments whose objective is show the efficiency of the schemes to preserve steady states that involve a varying bathymetry and friction. In~this case, $k$ takes the value of $0.01$ and the space domain is again the interval $[0,1]$. The~bathymetry is given by
\begin{equation}\label{swf_fondo}
H(x)= 1 - \displaystyle \frac{1}{2} \frac{e^{\cos(4 \pi x)} -e^{-1}}{e - e^{-1}}.
\end{equation}
The first experiment concerns the preservation of the supercritical stationary solution corresponding to \eqref{u'=Gswf} with initial conditions $q(0)=1$ and $h(0)=0.3$.

Following Reference~\cite{michel2017well}, we take 100 discretization cells, and~the numerical simulation of this experiment is run until $t = 1$~s.

Figure~\ref{test52_est_numsol} shows the numerical solutions at $t=1$~s for NWBM$i$, $i=1,2,3$ and CLWBM$i$, $i=1,2,3$. Notice that NWBM1 is completely wrong: due to the numerical diffusion, the~supercritical regime is lost at the right side of the domain, which causes a shock wave traveling to the left (of course, this behavior disappears when the mesh is fine enough). Therefore, in~Figure~\ref{test52_est_dif}, where the differences with the stationary solution are plotted, only the results for NWBM$i$, $i=2,3$ and for CLWBM$i$, $i=1,2,3$ are shown.
Table~\ref{test52_error_est} shows the errors corresponding to the different methods for the 100-cell mesh. We also consider a reference solution computed with a well-balanced first order scheme on a $1600$-cell~mesh.
\vspace{-20pt}

\begin{figure}[H]
{\captionsetup{position=bottom,justification=centering}
  \subfloat[NWBM$i$, $i=1,2,3$. Free surface.]{
   \includegraphics[width=0.38\textwidth]{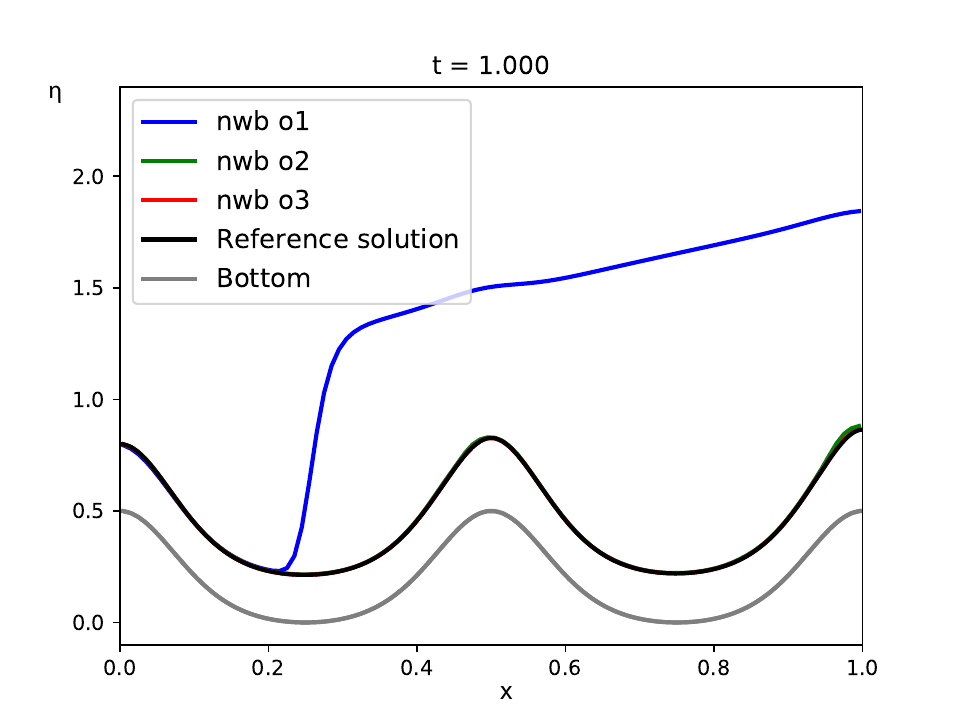}}
   \subfloat[CLWBM$i$, $i=1,2,3$. Free surface.]{
   \includegraphics[width=0.38\textwidth]{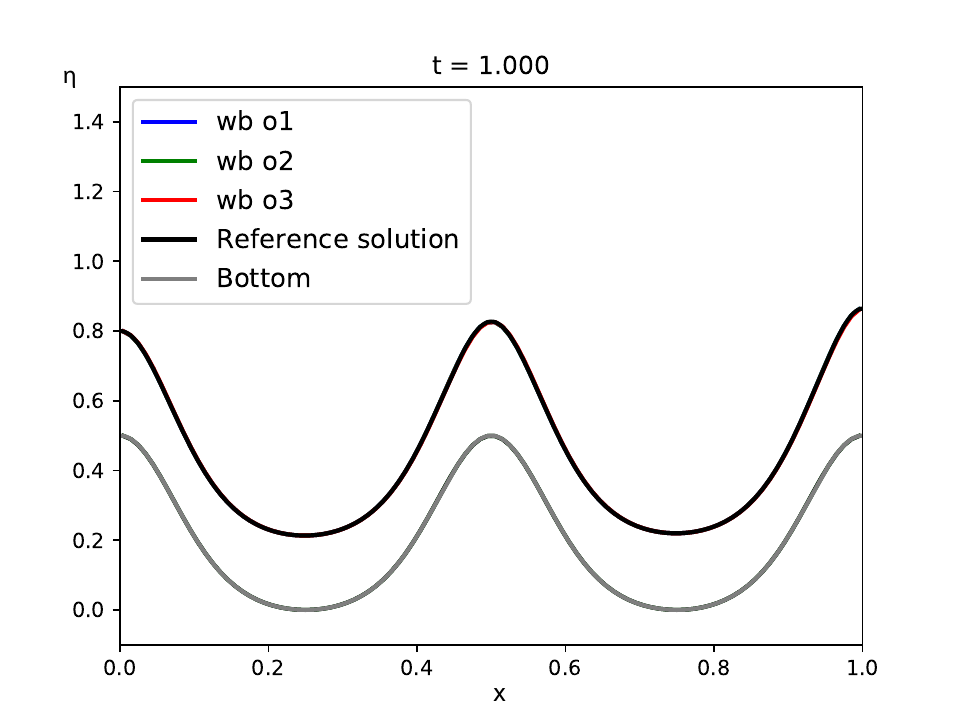}}
   \vspace{-10pt}
 \subfloat[NWBM$i$, $i=1,2,3$. $q$.]{
   \includegraphics[width=0.38\textwidth]{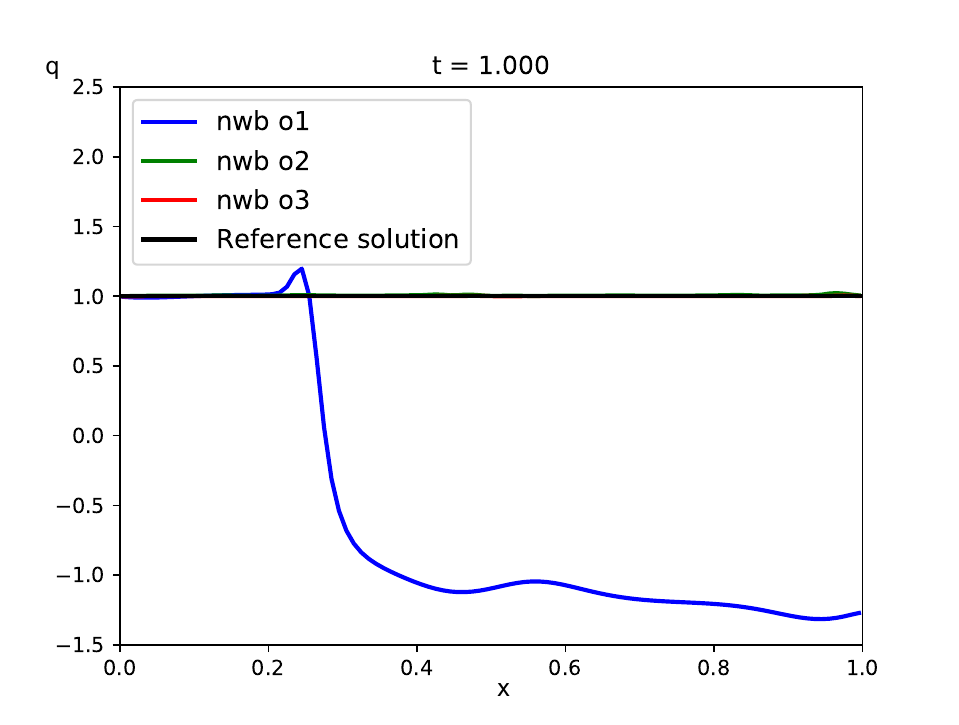}}
   \subfloat[CLWBM$i$, $i=1,2,3$. $q$.]{
   \includegraphics[width=0.38\textwidth]{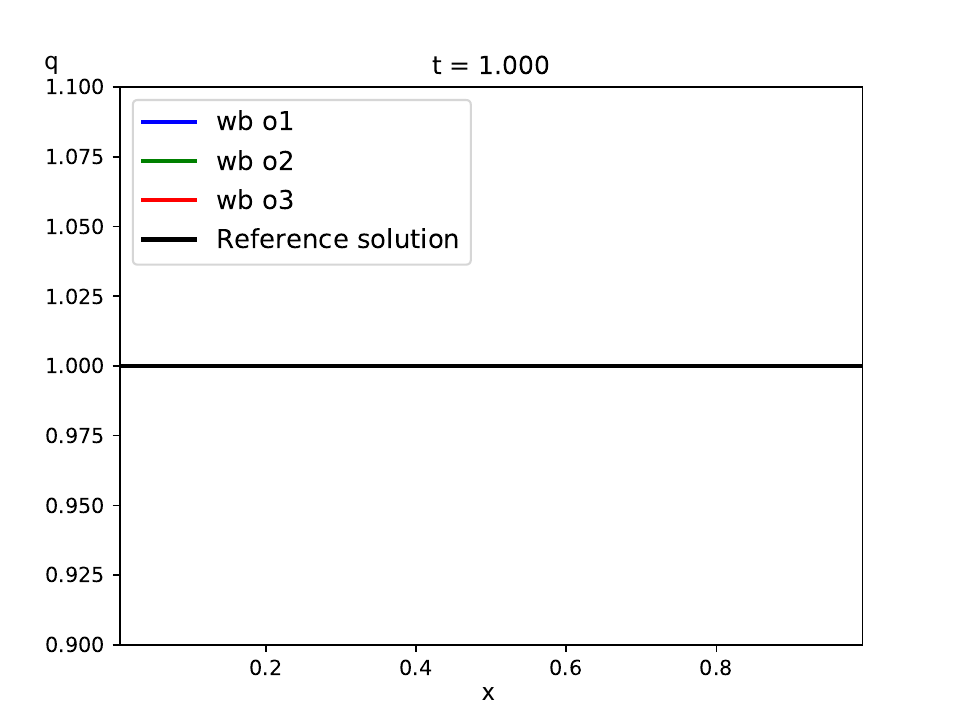}}}
    \caption{Test 4.3. 100-cell numerical solutions and reference solution at time $t = 1$~s. } \label{test52_est_numsol}
\end{figure}
\vspace{-28pt}

\begin{figure}[H]
{\captionsetup{position=bottom,justification=centering}
  \subfloat[NWBM$i$, $i=2,3$. $h$.]{
   \includegraphics[width=0.38\textwidth]{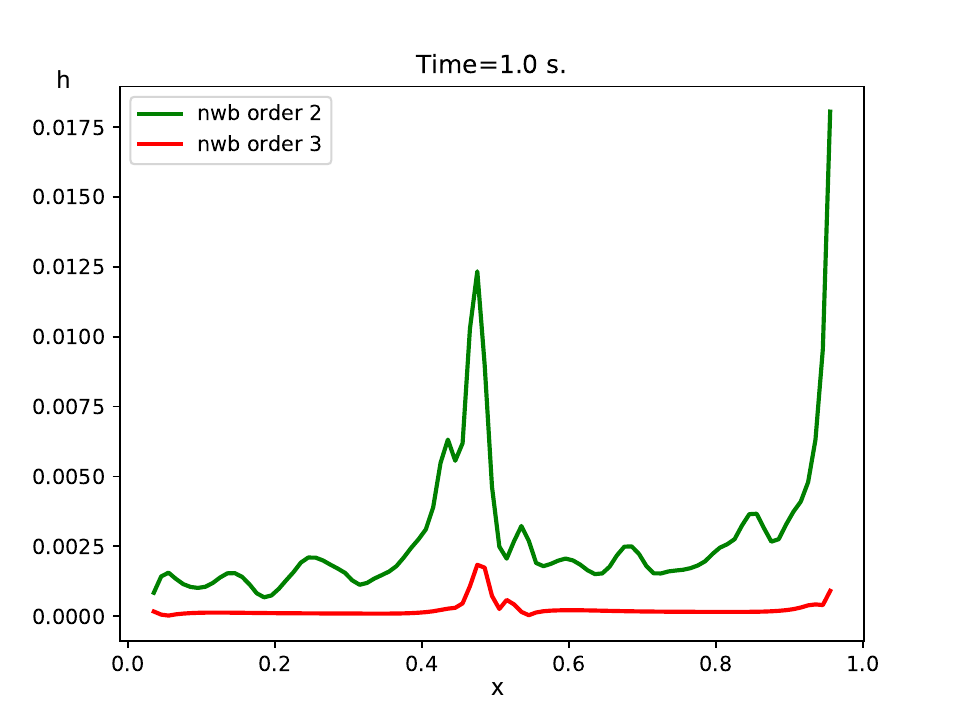}}
   \subfloat[CLWBM$i$, $i=1,2,3$. $h$.]{
   \includegraphics[width=0.38\textwidth]{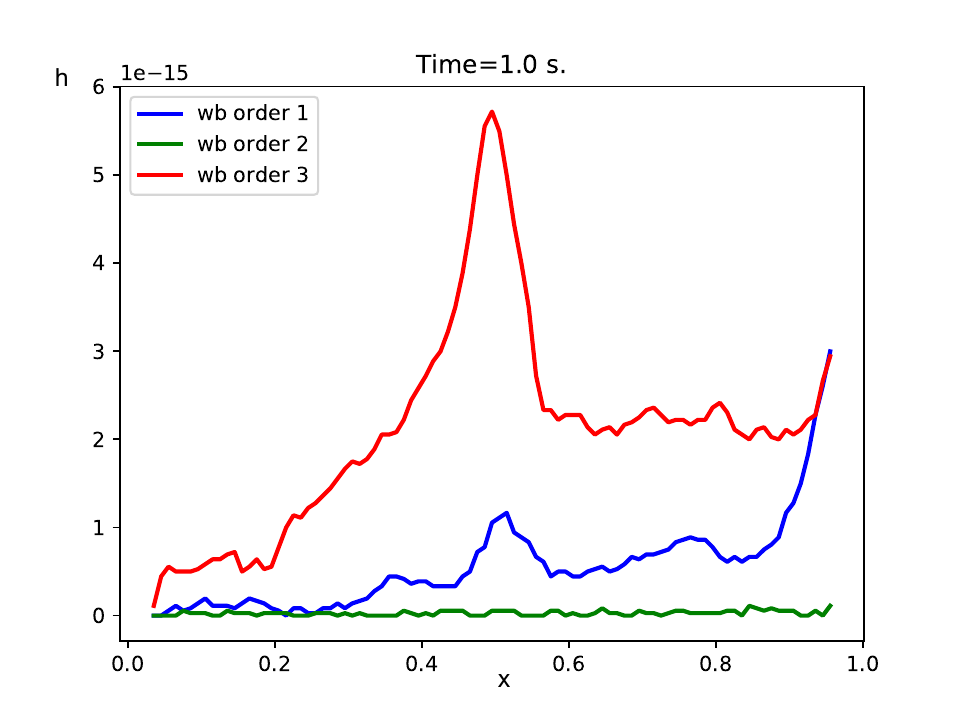}} }
 \caption{\textit{Cont}.}
\end{figure}
 
  \begin{figure}[H]\ContinuedFloat
\centering
\setcounter{subfigure}{3}

{\captionsetup{position=bottom,justification=centering}
 \subfloat[NWBM$i$, $i=2,3$. $q$.]{
   \includegraphics[width=0.38\textwidth]{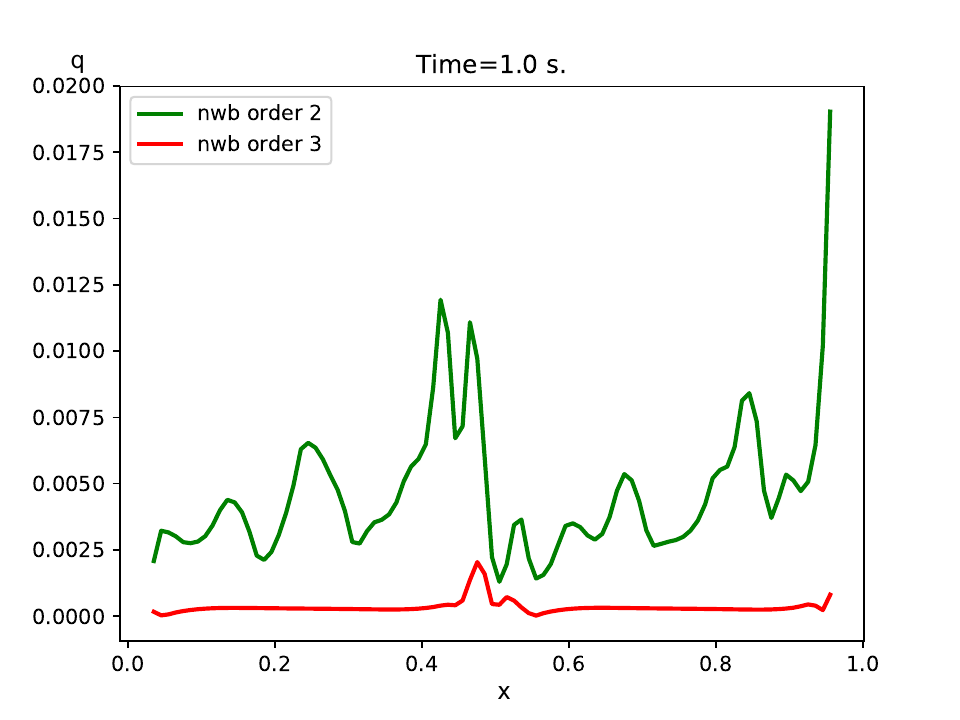}}
   \subfloat[CLWBM$i$, $i=1,2,3$. $q$.]{
   \includegraphics[width=0.38\textwidth]{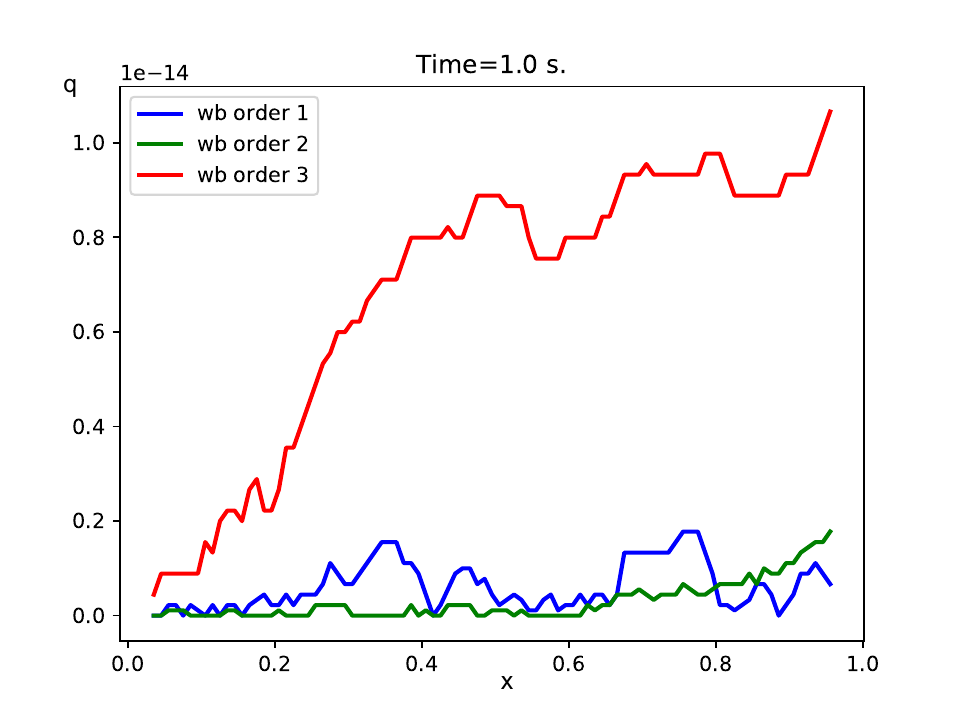}}}
    \caption{\hl{Test 4.3. Differences} between the stationary and the numerical solutions at time $t = 1$~s with a $100$-cell~mesh.} \label{test52_est_dif}
\end{figure}
\unskip

\end{paracol}
\nointerlineskip
\begin{specialtable}[H]
\widetable
\caption{Test 4.3. Differences in $L^1$-norm with respect to the stationary solution for NWBM$i$ and CLWBM$i$ ($i=1,2,3$) for the $100$-cell mesh at time $t=1$~s.} \label{test52_error_est}

\setlength{\cellWidtha}{\columnwidth/7-2\tabcolsep-0in}
\setlength{\cellWidthb}{\columnwidth/7-2\tabcolsep+0in}
\setlength{\cellWidthc}{\columnwidth/7-2\tabcolsep-0in}
\setlength{\cellWidthd}{\columnwidth/7-2\tabcolsep+0in}
\setlength{\cellWidthe}{\columnwidth/7-2\tabcolsep-0in}
\setlength{\cellWidthf}{\columnwidth/7-2\tabcolsep+0in}
\setlength{\cellWidthg}{\columnwidth/7-2\tabcolsep-0in}

\scalebox{1}[1]{\begin{tabularx}{\columnwidth}{>{\PreserveBackslash\centering}m{\cellWidtha}>{\PreserveBackslash\centering}m{\cellWidthb}>{\PreserveBackslash\centering}m{\cellWidthc}>{\PreserveBackslash\centering}m{\cellWidthd}>{\PreserveBackslash\centering}m{\cellWidthe}>{\PreserveBackslash\centering}m{\cellWidthf}>{\PreserveBackslash\centering}m{\cellWidthg}}
\toprule
\textbf{Method} & \multicolumn{2}{c}{\boldmath{\textbf{Error ($i=1$)}} }&\multicolumn{2}{c}{\boldmath{\textbf{Error ($i=2$)}}}&\multicolumn{2}{c}{\boldmath{\textbf{Error ($i=3$)}}}\\
  & \boldmath{$h$}&\boldmath{$q$ }&\boldmath{$h$}&\boldmath{ $q$} &\boldmath{$h$}&\boldmath{$q$}\\\midrule
NWBM$i$& 8.28 $\times$ $10^{-1}$  & 1.54 & 3.61 $\times$ $10^{-3}$ & 4.88 $\times$ $10^{-3}$ & 3.65 $\times$ $10^{-3}$ & 4.30 $\times$ $10^{-4}$\\
CLWBM$i$& 7.03 $\times$ $10^{-16}$ & 5.85 $\times$ $10^{-16}$ & 3.22 $\times$ $10^{-17}$ & 3.75 $\times$ $10^{-16}$ & 2.14 $\times$ $10^{-15}$ & 6.87 $\times$ $10^{-15}$\\\bottomrule

\end{tabularx}}
\end{specialtable}
\begin{paracol}{2}
\switchcolumn

\vspace{-8pt}

\subsubsection{Test~4.4}

Following Reference~\cite{michel2017well}, the~very last experiment in this subsection focuses on a perturbation of the aforementioned steady state. The~initial condition $U_0(x)=[h_0(x),q_0(x)]^T$ is given by: 
$$
h_0(x)=\begin{cases}
h^*(x)+0.05,
& \mbox{if $x \in \left[ \displaystyle \frac{2}{7} , \displaystyle \frac{3}{7} \right] \cup \left[ \displaystyle \frac{4}{7} , \displaystyle \frac{5}{7} \right]$,}\\
h^*(x), & \mbox{otherwise,}
\end{cases}
$$
$$
q_0(x)=\begin{cases}
q^*(x)+0.5,
& \mbox{if $x \in \left[ \displaystyle \frac{2}{7} , \displaystyle \frac{3}{7} \right] \cup \left[ \displaystyle \frac{4}{7} , \displaystyle \frac{5}{7} \right]$,}\\
q^*(x), & \mbox{otherwise,}
\end{cases}
$$
where $U^*(x)=[h^*(x),q^*(x)]^T$ is the stationary solution considered in the previous experiment. As~in Reference~\cite{michel2017well}, we use a 100-cell mesh for the numerical simulation on the domain $[0, 1]$, and the computations are carried out until the final time $t= 2$ s.

Observe that there are not big differences among the solutions during for small times. Nevertheless, for the NNWBM$1$, the supercritical regime is lost on the right of the domain, causing a shock traveling to the left, as~in Test 4.3. Figure~\ref{test52_per_dif} shows then the differences between the numerical solutions and the stationary solution for $h$ with NWBM$i$, $i=1,2,3$ at time $t=0.015$ and $i=2,3$ at time $t=2$, and~with CLWBM$i$, $i=1,2,3$ at times $t=0.015$ and $t=2$ (the plots are similar for $q$). Again, in~black, we plot a reference solution computed on a $1600$-cell mesh with a first order well-balanced~scheme.

\vspace{-20pt}

\begin{figure}[H]
{\captionsetup{position=bottom,justification=centering}
 \subfloat[NWBM$i$, $i =1,2,3$. $t = 0.015$~s. ]{
   \includegraphics[width=0.38\textwidth]{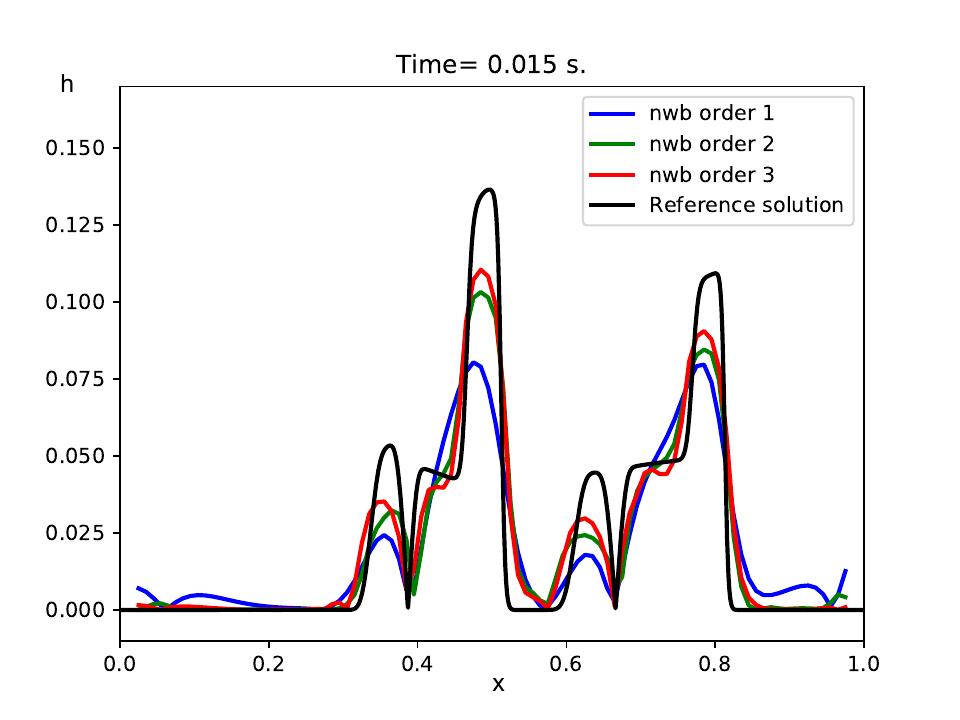}}
 \subfloat[CLWBM$i$, $i=1,2,3$. $t = 0.015$~s. ]{
   \includegraphics[width=0.38\textwidth]{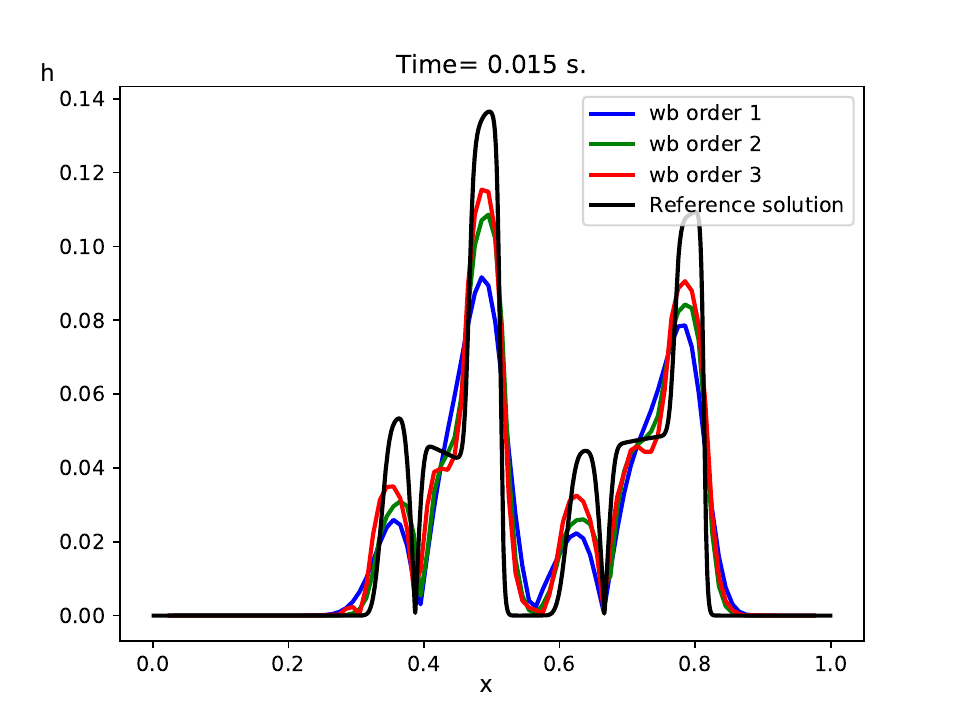}}\vspace{-10pt}
   \subfloat[NWBM$i$, $i =2,3$. $t = 2$~s. ]{
   \includegraphics[width=0.38\textwidth]{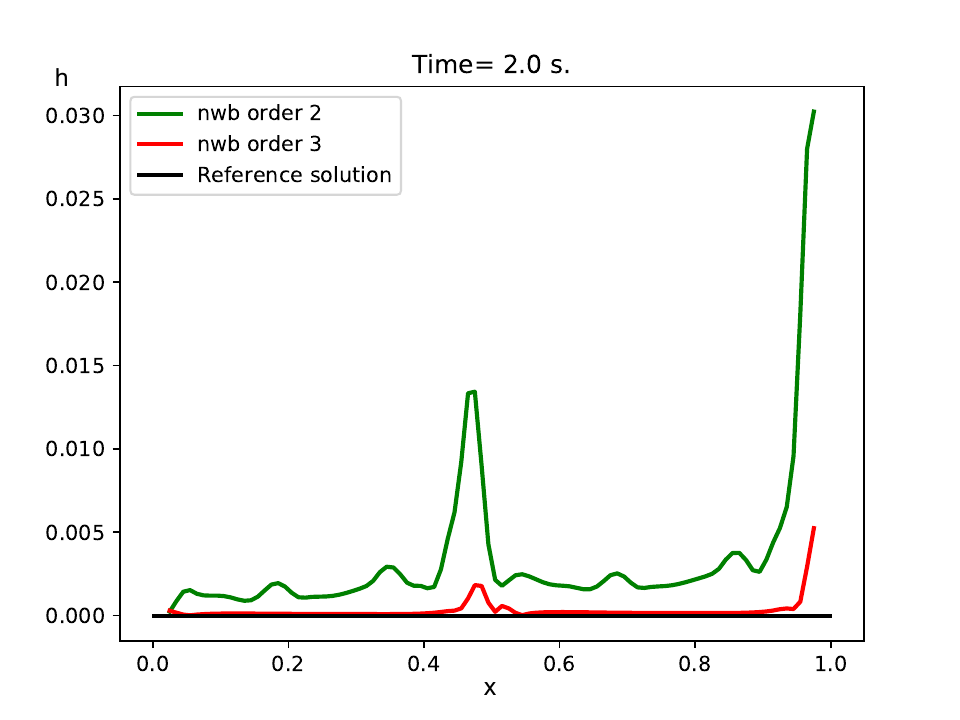}}
   \subfloat[CLWBM$i$, $i=1,2,3$. $t = 2$~s.]{
   \includegraphics[width=0.38\textwidth]{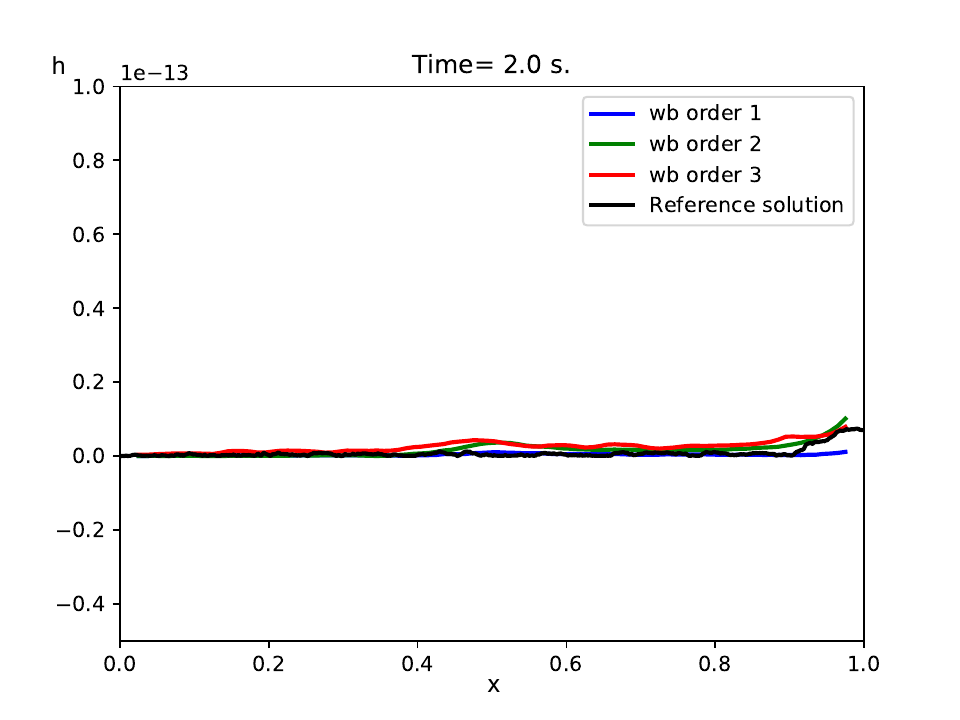}}}
    \caption{\hl{Test 4.4. Reference }and numerical solutions: differences with the stationary solution at times $t = 0.015, 2$ s for $h$, for~a $100$-cell~mesh.} \label{test52_per_dif}

\end{figure}

Once more, well-balanced schemes are able to recover the stationary states after the perturbation left the domain. Table~\ref{test52_error_per} shows the differences in $L^1$-norm with respect to the stationary solution at time $t = 2$ s for a mesh with $100$ cells. 


\end{paracol}
\nointerlineskip
\begin{specialtable}[H]
\widetable
\caption{Test 4.4. Differences in $L^1$-norm with respect to the stationary solution for NWBM$i$ and CLWBM$i$ ($i=1,2,3$) at time $t=2$ s for a $100$-cell~mesh.} \label{test52_error_per}
\setlength{\cellWidtha}{\columnwidth/7-2\tabcolsep-0in}
\setlength{\cellWidthb}{\columnwidth/7-2\tabcolsep+0in}
\setlength{\cellWidthc}{\columnwidth/7-2\tabcolsep-0in}
\setlength{\cellWidthd}{\columnwidth/7-2\tabcolsep+0in}
\setlength{\cellWidthe}{\columnwidth/7-2\tabcolsep-0in}
\setlength{\cellWidthf}{\columnwidth/7-2\tabcolsep+0in}
\setlength{\cellWidthg}{\columnwidth/7-2\tabcolsep-0in}

\scalebox{1}[1]{\begin{tabularx}{\columnwidth}{>{\PreserveBackslash\centering}m{\cellWidtha}>{\PreserveBackslash\centering}m{\cellWidthb}>{\PreserveBackslash\centering}m{\cellWidthc}>{\PreserveBackslash\centering}m{\cellWidthd}>{\PreserveBackslash\centering}m{\cellWidthe}>{\PreserveBackslash\centering}m{\cellWidthf}>{\PreserveBackslash\centering}m{\cellWidthg}}
\toprule
\textbf{Method} & \multicolumn{2}{c}{\boldmath{\textbf{Error ($i=1$)}} }&\multicolumn{2}{c}{\boldmath{\textbf{Error ($i=2$)}}}&\multicolumn{2}{c}{\boldmath{\textbf{Error ($i=3$)}}}\\
  & \boldmath{$h$}&\boldmath{$q$ }&\boldmath{$h$}&\boldmath{ $q$} &\boldmath{$h$}&\boldmath{$q$}\\\midrule
  
NWBM$i$& 2.42  & 6.12 & 3.57 $\times$ $10^{-3}$ & 4.87 $\times$ $10^{-3}$ & 1.39 $\times$ $10^{-3}$ & 4.30 $\times$ $10^{-4}$\\
CLWBM$i$& 3.73 $\times$ $10^{-16}$ & 3.60 $\times$ $10^{-16}$ & 1.80 $\times$ $10^{-15}$ & 1.99 $\times$ $10^{-15}$ & 2.64 $\times$ $10^{-15}$ & 8.93 $\times$ $10^{-15}$\\\bottomrule

\end{tabularx}}

\end{specialtable}
\begin{paracol}{2}
\switchcolumn

\vspace{-8pt}

\subsection{Problem 5: Compressible Euler Equations with Gravitational~Force}
Finally, let us consider the compressible Euler equations for gas dynamics with gravity:
\begin{equation}\label{euler_equations}
\begin{cases}
\rho_t+(\rho u)_x=0,\step
(\rho u)_t +\left( \rho u^2 + p \right)_x= -\rho H_x,\step
E_t +\left( u (E + p) \right)_x= -\rho u H_x.\\
\end{cases}
\end{equation}
Here, $u$ is the velocity, $p\geq 0$ the pressure, $\rho \geq 0$ the density, $q=\rho u$ the momentum, $E$ the total energy per unit volume, and~$H(x)$ is the gravitational potential. The~internal energy $e$ is defined by $\rho e = E - \displaystyle \frac{1}{2} \rho u^2$, and pressure could be determined from the internal energy using the equation of state. In~this work, we suppose an ideal gas; therefore,
\begin{equation}\label{EoS}
    p= (\gamma -1 ) \rho e,
\end{equation}
$\gamma >1$ being the adiabatic constant. Here, it is set to $\gamma = 1.5$. 

Note that \eqref{euler_equations} is a particular case of \eqref{sle} with $N=3$,
$$ U =\begin{pmatrix}
\rho \step
\rho u \step
E \\
\end{pmatrix} , \quad f(U) =\begin{pmatrix}
\rho u \step
 \rho u^2 + p\step
 u(E+p)\\
\end{pmatrix}, \quad S(U) =\begin{pmatrix}
0\step
- \rho\step
- \rho~u\\
\end{pmatrix}.$$

Supposing that the system is strictly hyperbolic and considering \eqref{EoS}, the~stationary solutions satisfy the ODE system:
\begin{equation}\label{u'=Geuler}
\begin{cases}
q_x=0,\step
\displaystyle \frac{d \hat{U}}{dx} = G(x, \hat{U}), \\
\end{cases}
\end{equation}
where
\[ \hat{U} =\begin{pmatrix}
\rho \step
E \\
\end{pmatrix} , \quad G(x, \hat{U}) = - \begin{pmatrix}
\displaystyle \frac{\rho}{c^2 - u^2} \\[4mm]
 \displaystyle \frac{\rho}{\gamma-1} \left( 1 + \frac{3-\gamma}{2} \frac{u^2}{c^2 - u^2} \right)\\
\end{pmatrix} H_x,\]
where $c = \sqrt{ \gamma {p}/{\rho}}$ is the sound speed, which is valid for stationary regular solution. Notice that, for regular solutions of the Euler equations, the entropy is constant along material lines. So that, for~stationary solutions, $u s_x = 0$, where $s$ denotes the entropy density. Therefore, if $u$ does not vanish, it follows that these stationary solutions are~isentropic.

\subsubsection*{Test~5.1.}
In this test, we consider the space domain $[-1,1]$, and the system is integrated up to time $t=5$~s. The~gravity potential is $H(x)=x$. 

As initial condition, a~supersonic stationary solution corresponding to the solution of the Cauchy problem
\begin{equation} \label{test61_cini}
\begin{cases}
q_x=0,\step
\displaystyle \frac{d \hat{U}}{dx} = G(x, \hat{U}), \step
\rho(-1)=1, \, q(-1)=10, \, E(-1)=52
\end{cases}
\end{equation}
is computed with the Gauss-Legendre collocation method. Density, momentum, and energy is set downstream, and open free boundary conditions upstream. $CFL$ parameter is set to~$0.9$.

Figures~\ref{fluc_test61_rho}--\ref{fluc_test61_E} show the differences between the stationary and the numerical solutions for NWBM$i$ and CLWBM$i$, $i=1,2,3$ for $\rho$, $q$, and $E$, respectively; Table~\ref{test61_error} shows the errors for for the 100-cell mesh with NWBM$i$ and CLWBM$i$, $i=1,2,3$ at time $t=5$~s. \mbox{Table~\ref{euler_times}} shows the computational cost for NWBM$i$, CWBM$i$, and CLWBM$i$, $i=1,2,3$. Notice that CLWBM are significantly more efficient than CWBM. We have carried out similar experiments with a perturbation of this stationary solution as initial condition, and~we obtained the expected results, i.e., only the well-balanced methods are able to exactly recover the steady state solution. For~the sake of conciseness, we do not present the tests~here. 
\vspace{-20pt}

\begin{figure}[H]
{\captionsetup{position=bottom,justification=centering}
  \subfloat[NWBM$i$, $i=1,2,3$. Density ]{
   \includegraphics[width=0.38\textwidth]{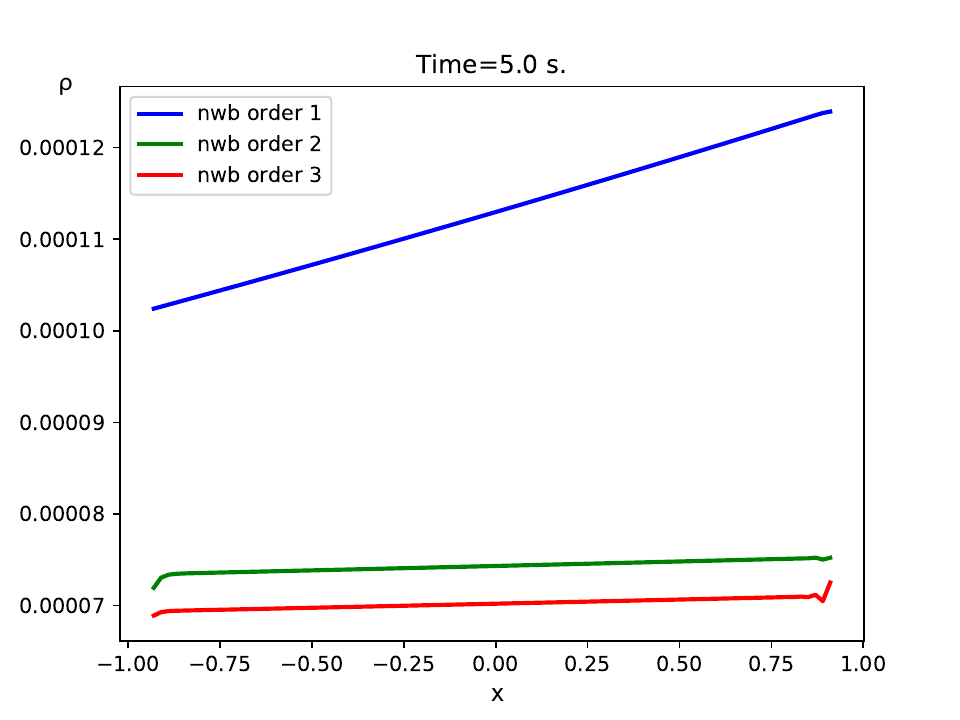}}
 \subfloat[CLWBM$i$, $i=1,2,3$. Density ]{
   \includegraphics[width=0.38\textwidth]{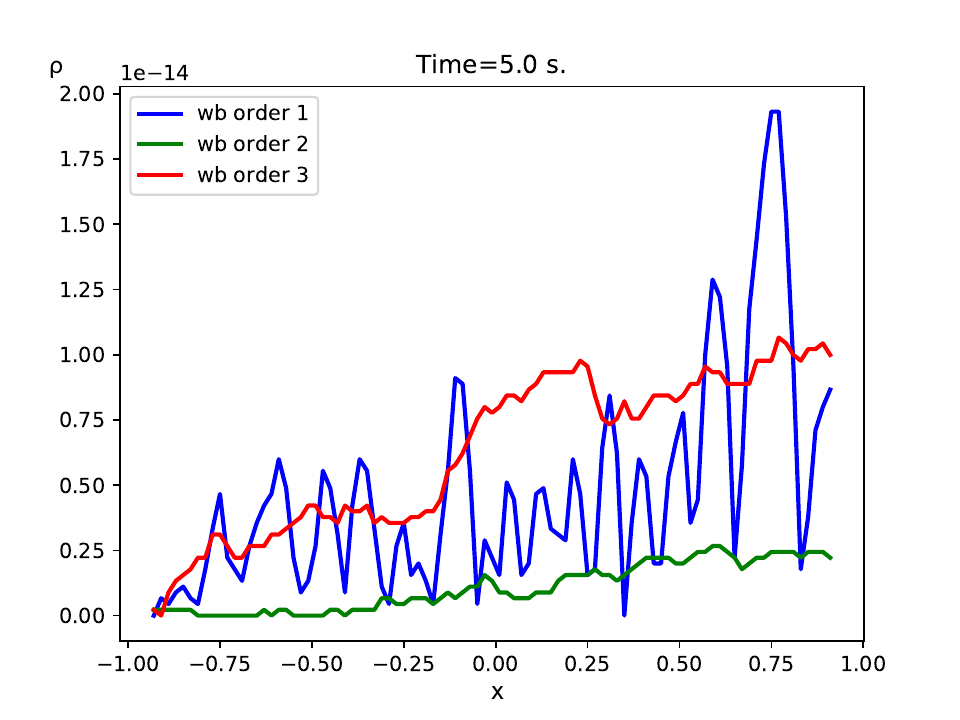}}}
    \caption{\hl{Test 5.1. Differences }between the stationary and the numerical solutions at time $t = 5$~s for $\rho$ with a mesh of $100$ cells.} \label{fluc_test61_rho}

\end{figure}
\vspace{-20pt}

\begin{figure}[H]
{\captionsetup{position=bottom,justification=centering}  
    \subfloat[NWBM$i$, $i=1,2,3$. Momentum]{
   \includegraphics[width=0.38\textwidth]{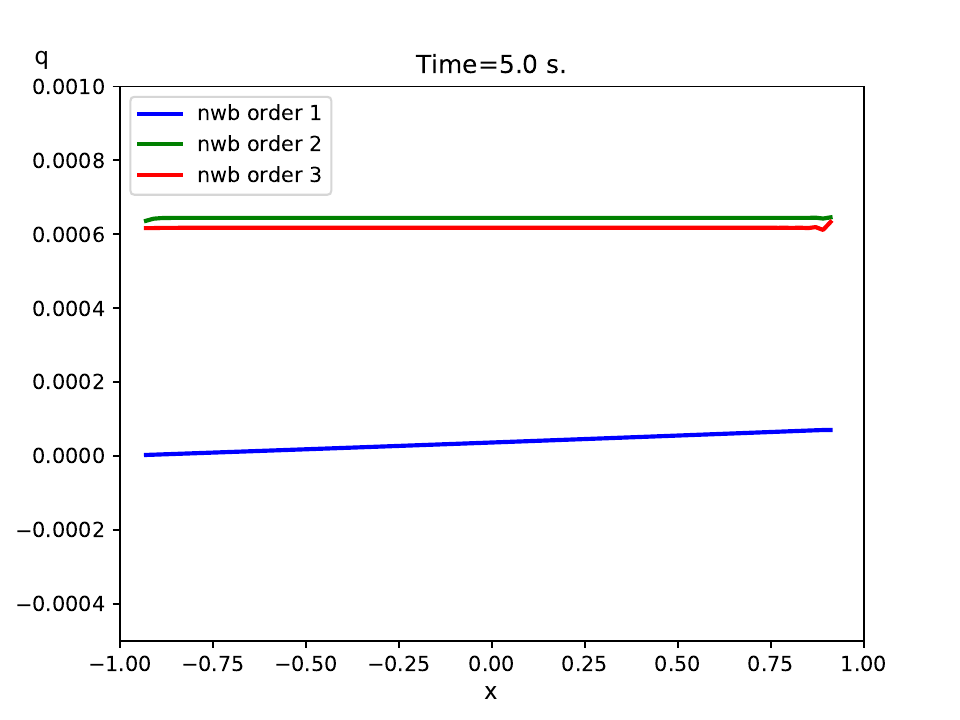}}
 \subfloat[CLWBM$i$, $i=1,2,3$. Momentum ]{
   \includegraphics[width=0.38\textwidth]{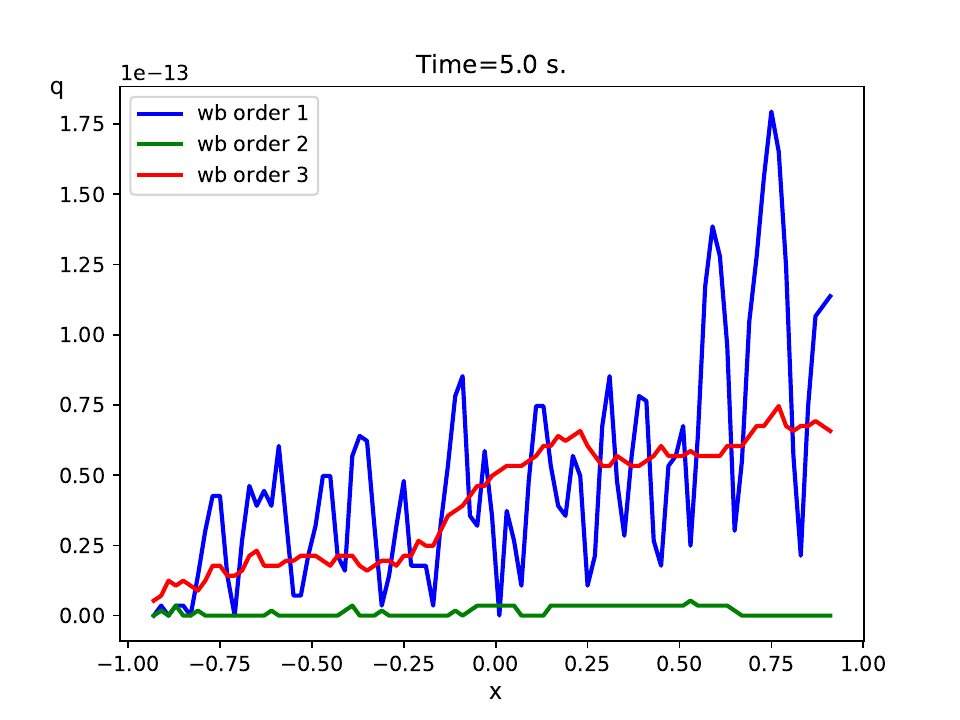}}}   
    \caption{\hl{Test 5.1. Differences }between the stationary and the numerical solutions at time $t = 5$~s for $q$ with a mesh of $100$ cells.} \label{fluc_test61_q}
\end{figure}
\vspace{-20pt}

\begin{figure}[H]
{\captionsetup{position=bottom,justification=centering}  
   \subfloat[NWBM$i$, $i=1,2,3$. Energy]{
   \includegraphics[width=0.38\textwidth]{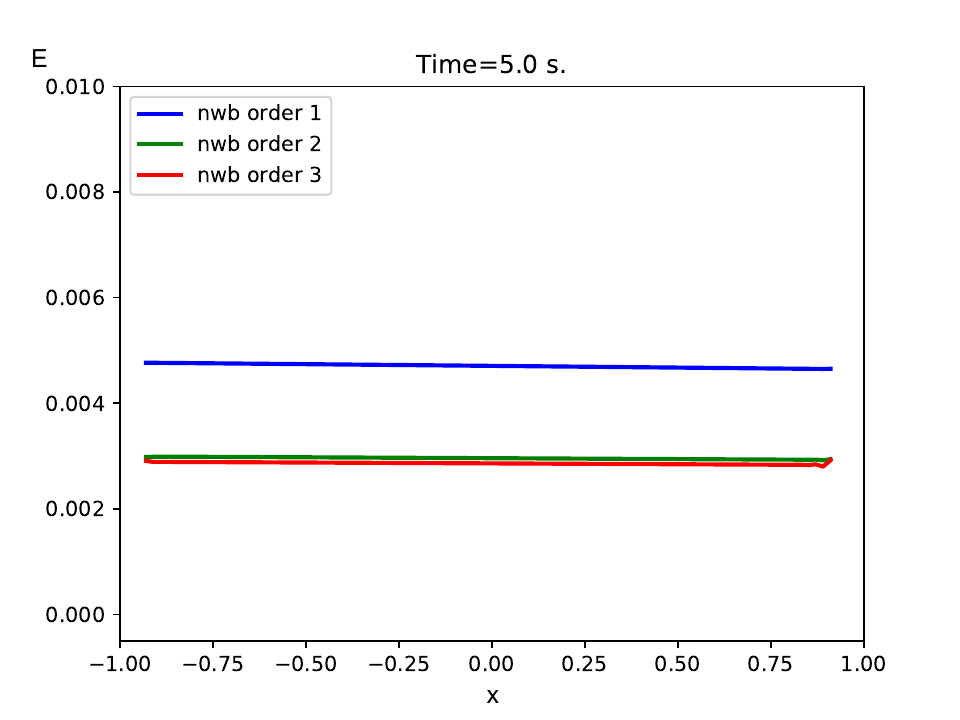}}
 \subfloat[CLWBM$i$, $i=1,2,3$. Energy]{
   \includegraphics[width=0.38\textwidth]{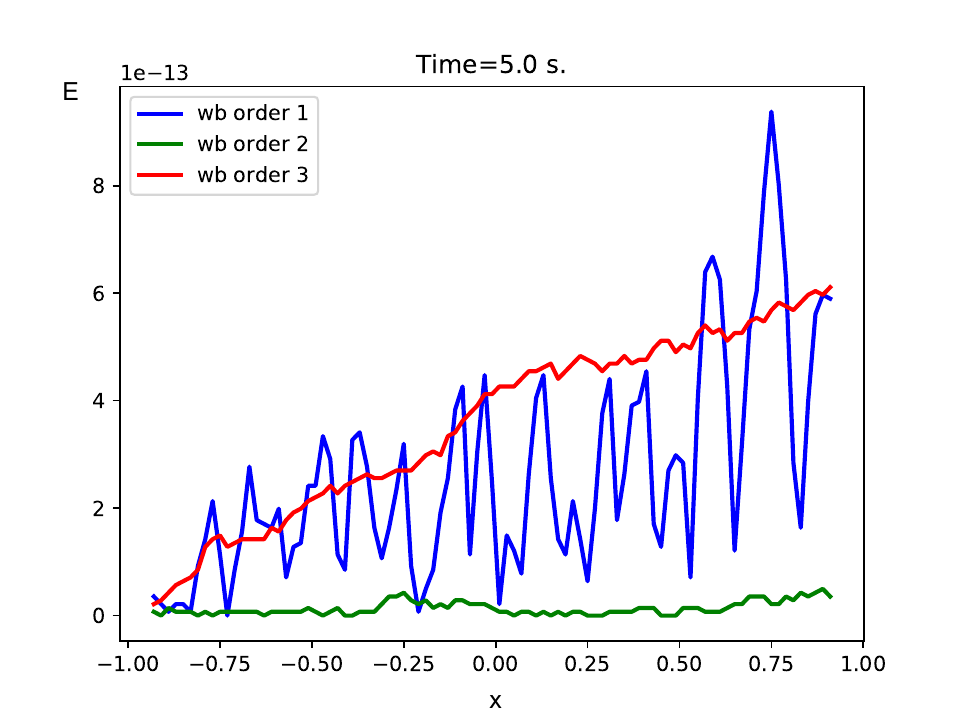}}}
    \caption{Test 5.1. Differences between the stationary and the numerical solutions at time $t = 5$~s for $E$ with a mesh of $100$ cells. } \label{fluc_test61_E}
\end{figure}
\vspace{-20pt}

\begin{specialtable}[H]
\caption{Test 5.1. Differences in $L^1$-norm with respect to the stationary solution for NWBM$i$ and CLWBM$i$ ($i=1,2,3$) for the mesh with $100$ cells at time $t=5$~s.} \label{test61_error}
\setlength{\cellWidtha}{\columnwidth/4-2\tabcolsep-0in}
\setlength{\cellWidthb}{\columnwidth/4-2\tabcolsep+0in}
\setlength{\cellWidthc}{\columnwidth/4-2\tabcolsep-0.0in}
\setlength{\cellWidthd}{\columnwidth/4-2\tabcolsep-0.0in}

\scalebox{1}[1]{\begin{tabularx}{\columnwidth}{>{\PreserveBackslash\centering}m{\cellWidtha}>{\PreserveBackslash\centering}m{\cellWidthb}>{\PreserveBackslash\centering}m{\cellWidthc}>{\PreserveBackslash\centering}m{\cellWidthd}}
\toprule

\textbf{Method} & \boldmath{\textbf{Error ($i=1$)}} &\boldmath{\textbf{Error ($i=2$)}}&\boldmath{\textbf{Error ($i=3$)}}\\\midrule
\multicolumn{4}{c}{$h$}\\\midrule
NWBM$i$& 2.23 $\times$ $10^{-4}$  & 9.41 $\times$ $10^{-5}$ & 9.51 $\times$ $10^{-3}$ \\
CLWBM$i$& 6.97 $\times$ $10^{-15}$ & 6.58 $\times$ $10^{-14}$ & 3.20 $\times$ $10^{-13}$ \\\midrule
\multicolumn{4}{c}{$q$}\\\midrule
NWBM$i$& 1.45 $\times$ $10^{-4}$  & 1.28 $\times$ $10^{-3}$ & 5.95 $\times$ $10^{-3}$ \\
CLWBM$i$&  2.22 $\times$ $10^{-15}$ & 2.81 $\times$ $10^{-15}$ & 2.77 $\times$ $10^{-14}$ \\\midrule
\multicolumn{4}{c}{$E$}\\\midrule
NWBM$i$& 1.38 $\times$ $10^{-4}$  & 1.22 $\times$ $10^{-3}$ & 5.78 $\times$ $10^{-3}$ \\
CLWBM$i$&  1.24 $\times$ $10^{-14}$  & 8.13 $\times$ $10^{-14}$ & 7.15 $\times$ $10^{-13}$\\\bottomrule
\end{tabularx}}
\end{specialtable}
\unskip

\begin{specialtable}[H]
\caption{Test 5.1. Computational cost (milliseconds) for the mesh with $100$ cells. $t=5$~s.}\label{euler_times}
\setlength{\cellWidtha}{\columnwidth/4-2\tabcolsep-0in}
\setlength{\cellWidthb}{\columnwidth/4-2\tabcolsep+0in}
\setlength{\cellWidthc}{\columnwidth/4-2\tabcolsep-0.0in}
\setlength{\cellWidthd}{\columnwidth/4-2\tabcolsep-0.0in}

\scalebox{1}[1]{\begin{tabularx}{\columnwidth}{>{\PreserveBackslash\centering}m{\cellWidtha}>{\PreserveBackslash\centering}m{\cellWidthb}>{\PreserveBackslash\centering}m{\cellWidthc}>{\PreserveBackslash\centering}m{\cellWidthd}}
\toprule

        \boldmath{ \textbf{ Order ($i$)  }   }& \boldmath{  \textbf{ NWBM$i$} }& \boldmath{\textbf{CWBM$i$} }&\boldmath{\textbf{ CLWBM$i$ }}\\
 \midrule
 1& 120 & 4480 & 220 \\ 
                   2 & 250 & 8960 & 1440 \\ 
                   3& 670 & 19430 & 4750 \\\bottomrule

\end{tabularx}}
\end{specialtable}
\unskip

\section{Conclusions}\label{sec5}

Following Reference~\cite{sinum2008}, we propose a family of well-balanced high-order numerical schemes that can be applied to general 1D balance laws. Note that the main difficulty of the method proposed in Reference~\cite{sinum2008} comes from the fact that local stationary solutions must be computed at every cell that may require the knowledge of the stationary states. When solving the ODE \eqref{sblst} for the stationary solution is not possible or too costly, the~first step has to be computed numerically: what we propose here is to build both the stationary solutions and the local solvers for the first step on the basis of the collocation RK methods. The~well-balanced property of these numerical methods is precisely stated, and~we have also introduced a general technique that allows us to deal with resonant problems for any 1D systems of balance~laws.

We have considered a good number of examples: Burgers equation, shallow water system with friction, and Euler equations with gravity to check the well-balanced property of the methods and to test their efficiency. In~some of these test cases, the~stationary solutions are known either in implicit or explicit form, while, in others, the only information comes from the ODE \eqref{sblst} that the stationary solutions satisfy: in the former case, we can compare the efficiency of the new implementation, while, in the latter, we can show the generality of the~methods. 

Although well-balanced schemes are more costly, they are more effective than standard non-well-balanced methods when they are applied to the propagation of small perturbation around an equilibria or when long-time integration is required. Furthermore, the~tests show that the strategy based on the collocation RK methods introduced in this work is more efficient than the one presented in Reference~\cite{GomezCastroPares2020}. 

The methods considered here can be applied in principle to any system of balance laws. For~instance, in~the context of Euler equations, in Reference~\cite{gaburro2018well}, a~first and second order numerical scheme that is well-balanced for particular 1D radial solutions for the Euler equations in cylindrical coordinates has been developed following {the strategy proposed by two of the authors that was later reviewed in} Reference \cite{CastroPares2019}. It is possible to use
the technique described here to extend the methods in Reference~\cite{gaburro2018well} to more general source terms and higher-order. This will be considered in a future~work.

Futher developments will include applications to:

\begin{itemize}
    \item Balance laws \eqref{sle} with non-regular source terms, that is systems for which $H$ has jump discontinuities.
    \item Multidimensional problems.
    \end{itemize}

\authorcontributions{\hl{x}.}

\funding{This research has been partially supported by the Spanish Government and FEDER through the coordinated Research project RTI2018-096064-B-C1, The Junta de Andaluc\'ia research project P18-RT-3163 and the Junta de Andalucia-FEDER-University of M\'alaga Research project UMA18-FEDERJA-161. G. Russo acknowledges partial support from ITN-ETN Horizon 2020 Project ModCompShock, ``Modeling and Computation of Shocks and Interfaces'', Project Reference 642768, and~from the Italian Ministry of University and Research (MIUR), PRIN Project 2017 (No. 2017KKJP4X entitled “Innovative numerical methods for evolutionary partial differential equations and applications”. I. Gómez-Bueno is also supported by a Grant from “El Ministerio de Ciencia, Innovación y Universidades”, Spain (FPU2019/01541).}

\institutionalreview{\hl{x}.}

\informedconsent{\hl{x}.}

\dataavailability{\hl{x}.}
 
 \conflictsofinterest{\hl{x}.}

\end{paracol}

\reftitle{References}

\externalbibliography{yes}

\end{document}